\documentclass{article}
\usepackage[T1]{fontenc}
\usepackage{amsthm, fullpage} 
\usepackage{amsmath}
\usepackage{amssymb}
\usepackage{amsfonts}
\usepackage{eucal}
\usepackage[all]{xy}
\usepackage[frenchb]{babel}
\usepackage{pslatex}
\usepackage[pagebackref]{hyperref}

\newtheorem{prop}{Proposition}[section]
\newtheorem{theo}[prop]{Théor\`eme}
\newtheorem*{theo**}{Théorème}
\newtheorem{coro}[prop]{Corollaire}
\newtheorem*{conj*}{Conjecture}
\newtheorem{lemm}[prop]{Lemme}

\theoremstyle{definition}
\newtheorem{conj}[prop]{Conjecture}
\newtheorem{vide}[prop]{}
\newtheorem{defi}[prop]{Définition}

\theoremstyle{remark}
\newtheorem{rema}[prop]{Remarques}
\newtheorem{exem}[prop]{Exemple}
\newtheorem{nota}[prop]{Notations}

\numberwithin{equation}{prop}

\newcommand{\riso}{ \overset{\sim}{\longrightarrow}\, }
\newcommand{\liso}{ \overset{\sim}{\longleftarrow}\, }

\newcommand{\Spf}{\mathrm{Spf}\,}

\renewcommand{\sp}{\mathrm{sp}}

\renewcommand{\AA}{{\mathcal{A}}}

\newcommand{\FF}{{\mathcal{F}}}
\newcommand{\B}{{\mathcal{B}}}
\newcommand{\CC}{{\mathcal{C}}}
\newcommand{\E}{{\mathcal{E}}}
\newcommand{\G}{{\mathcal{G}}}

\newcommand{\M}{{\mathcal{M}}}
\newcommand{\NN}{{\mathcal{N}}}
\newcommand{\D}{{\mathcal{D}}}
\newcommand{\I}{{\mathcal{I}}}
\newcommand{\PP}{{\mathcal{P}}}

\renewcommand{\O}{{\mathcal{O}}}

\newcommand{\V}{\mathcal{V}}
\renewcommand{\S}{\mathcal{S}}

\newcommand{\Y}{\mathcal{Y}}
\newcommand{\ZZ}{\mathcal{Z}}
\newcommand{\X}{\mathfrak{X}}

\newcommand{\m}{\mathfrak{m}}
\newcommand{\U}{\mathfrak{U}}

\newcommand{\DD}{\mathbb{D}}
\renewcommand{\L}{\mathbb{L}}
\newcommand{\R}{\mathbb{R}}
\newcommand{\Q}{\mathbb{Q}}
\newcommand{\Z}{\mathbb{Z}}
\newcommand{\N}{\mathbb{N}}

\newcommand{\hdag}{  \phantom{}{^{\dag} }    }

\begin{document}
\title{Log-isocristaux surconvergents et holonomie}
\author{Daniel Caro \footnote{L'auteur a bénéficié du soutien du réseau européen TMR \textit{Arithmetic Algebraic Geometry}
(contrat numéro UE MRTN-CT-2003-504917).}}

\maketitle

\begin{abstract}
Soient $\V$ un anneau de valuation discrète complet d'inégales caractéristiques $(0,p)$
de corps résiduel parfait $k$.
On se donne $\X$ un $\V$-schéma formel séparé et lisse,
$\ZZ$ un diviseur à croisements normaux strict de $\X$,
$\Y:=\X \setminus \ZZ$ l'ouvert correspondant,
$X$, $Z$, $Y$ les fibres spéciales respectives,
$T$ un diviseur de $X$,
$M (\ZZ)$ la log-structure sur $\X$ définie par $\ZZ$,
$\X ^\#:= (\X, M (\ZZ))$ le log-$\V$-schéma formel lisse induit et $u$ : $\X ^\# \rightarrow \X$ le morphisme structural.
Soit $\E $ un log-isocristal sur $\X ^\#$ surconvergent le long de $T$, i.e., un
$\D ^\dag _{\X^\#} (\hdag T) _{\Q}$-module à gauche cohérent localement projectif et de type fini sur $\O _{\X} (\hdag T) _{\Q}$.
Nous vérifions d'abord que le complexe
$u _{T,+} (\E)$ est réduit à un terme et que
$u _{T,+} (\E )$ est un $\D ^\dag _{\X} (\hdag T) _{\Q}$-module holonome
(via une extension de la notion d'holonomie sans structure de Frobenius et avec un diviseur).
Nous prouvons ensuite que le complexe
de de Rham logarithmique de $\E $ est isomorphe au complexe de de Rham
de $u _{T,+} (\E )$.
En notant $\E (\hdag Z)$ l'isocristal sur $Y\setminus T$ surconvergent le long de $Z\cup T$ induit par $\E $,
on dispose d'un morphisme canonique $\rho _{\E}$ : $u _{T,+} (\E )\rightarrow \E  (\hdag Z)$.
Nous établissons que lorsque $\E $ est en fait
un isocristal sur $\X $ surconvergent le long de $T$, $\rho _{\E}$ est un isomorphisme.
\end{abstract}

\tableofcontents

\section*{Introduction}
Soient $\V$ un anneau de valuation discrète complet d'inégales caractéristiques $(0,p)$,
d'idéal maximal $\mathfrak{m}$
et de corps résiduel parfait $k$.
Shiho dans \cite[3.1.8]{shihoII} et Tsuzuki dans \cite[1.3.1]{tsumono}
ont énoncé la conjecture suivante,
que l'on peut nommer la {\og conjecture de monodromie génériquement finie\fg} (terminologie de Shiho et Tsuzuki)
ou {\og conjecture de réduction semi-stable \fg} (terminologie de Kedlaya \cite{kedlaya_semi-stable}) :
\begin{conj*}
[*]
\label{semi-stable}
Soient $Y$ une variété lisse sur $k$ (i.e., un $k$-schéma lisse séparé de type fini)
et $E$ un $F$-isocristal surconvergent sur $Y$.
Il existe alors :
\begin{enumerate}
\item un morphisme propre, surjectif, génériquement étale $g$ : $Y _1 \rightarrow Y$,
\item une immersion ouverte $j _1$ : $Y _1 \hookrightarrow X _1$ dans une variété projective lisse telle que
$D _1 := X _1 \setminus Y _1$ soit un diviseur à croisements normaux strict de $X _1$,
\item un log-isocristal convergent $G _1$ sur ($X _1$, $M _1$)
(où $M _1$ est la log structure associée à ($X _1$, $D _1$)) sur $\Spf \V$,
\end{enumerate}
tels que $j _1 ^\dag (G _1)\riso g ^* E $, avec $j _1 ^\dag (G _1)$ signifiant l'isocristal
surconvergent sur $Y _1$ induit par $G _1$.
\end{conj*}

Cette conjecture est établi dans les deux cas particuliers suivants :

$\bullet $ Tsuzuki a obtenu dans \cite{tsumono} un théorème de monodromie génériquement finie et étale
pour les $F$-isocristaux surconvergents {\it unités}. Ce théorème est identique à la conjecture ci-dessous
à la différence près que les pentes nulles par Frobenius permettent d'éviter les structures logarithmiques
(i.e., on peut supposer que $G _1$ est un isocristal convergent sur $X _1$).

$\bullet $ Lorsque $Y$ est une courbe, Kedlaya (voir \cite[1.1]{kedlaya_semi-stable}) a prouvé cette conjecture.
D'après Kedlaya, le cas des surfaces devrait se déduire de ses travaux
\cite{kedlaya-semistableI}, \cite{kedlaya-semistableII}, \cite{kedlaya-semistableIII}.

Voici quelques conséquences de ces deux cas où la conjecture $(*)$ est validée :
nous avions prouvé, grâce au premier cas,
 la surholonomie (et donc l'holonomie) des $F$-isocristaux surconvergents unités sur les variétés lisses
 (voir \cite{caro_unite} et \cite{caro_surholonome}).
En utilisant le second cas, nous avions aussi obtenu
l'holonomie des $F$-isocristaux surconvergents sur les courbes lisses
(voir \cite[4]{caro_courbe-nouveau}).

Les travaux de Kedlaya portent à croire que la conjecture $(*)$
est exacte.
De plus, il est raisonnable de penser
que cette conjecture implique, de façon analogue au cas des courbes ou à celui des $F$-isocristaux surconvergents unités,
la surholonomie des $F$-isocristaux surconvergents sur les variétés lisses. Ce résultat constituerait
une réciproque au dévissage en $F$-isocristaux surconvergents des $F$-complexes surholonomes (voir \cite{caro_devissge_surcoh} ou
\cite{caro-2006-surcoh-surcv}).

Dans cet article, nous nous intéressons à l'holonomie des $F$-isocristaux surconvergents
provenant de log-isocristaux convergents,
ce qui fournit une première étape pour valider l'holonomie et la surholonomie
des $F$-isocristaux surconvergents sur les variétés lisses.
Précisons à présent
les différentes sections de
ce travail.

Soient $\X$ un $\V$-schéma formel séparé et lisse, $\ZZ$ un diviseur à croisements normaux strict de $\X$,
$\Y:=\X \setminus \ZZ$ l'ouvert correspondant,
$X$, $Z$, $Y$ les fibres spéciales respectives,
$M (\ZZ)$ la log-structure sur $\X$ définie par $\ZZ$,
$\X ^\#:= (\X, M (\ZZ))$ le log-$\V$-schéma formel lisse induit,
$T$ un diviseur de $X$ et $u$ : $\X ^\# \rightarrow \X$ le morphisme structural.
On désigne par $\D ^\dag _{\X^\#} (\hdag T) _{\Q}$ le faisceau des opérateurs différentiels de niveau fini sur $\X ^\#$
à singularités surconvergentes le long de $T$.
Soit $\E$ un log isocristal sur $\X ^\#$ surconvergent le long de $T$, i.e.,
  un $\D  ^\dag _{\X ^\#} (\hdag T) _{\Q} $-module cohérent,
  localement projectif et de type fini sur $\O _{\X} (\hdag T) _{\Q}$
(voir \ref{rema-isoctheodef}).

Dans les deux premières parties,
nous donnons quelques compléments sur les log-$\D$-modules arithmétiques qui étendent au cas logarithmique plusieurs
propriétés déjà connues dans le cas non logarithmique.
Nous étudierons notamment les suites de Spencer logarithmiques.
Nous prouvons dans la troisième partie un isomorphisme d'associativité de produits tensoriels
mélangeant $\D$-modules arithmétiques et log $\D$-modules arithmétiques (voir \ref{A2CNiso}).

Notons $\omega _{\X ^\#}$ et $\omega _{\X}$ les faisceaux des formes différentielles de degré maximal
sur respectivement $\X $ et $\X ^\#$.
On vérifie que $\O _{\X} (\ZZ) _{\Q}:= \mathcal{H} om _{\O _{\X,\Q}} ( \omega _{\X,\Q} , \omega _{\X ^\#,\Q})$,
où les indices $\Q$ signifient que l'on applique $-\otimes _\Z \Q $,
est muni d'une structure canonique de
$\D ^\dag _{\X^\#} (\hdag T) _{\Q}$-module à gauche et on pose
$\E (\ZZ) := \E \otimes _{\O _{\X,\Q}} \O _{\X} (\ZZ) _{\Q}$.
Nous étudierons dans la quatrième partie les
log isocristaux sur $\X ^\#$ surconvergents le long de $T$.
 Nous obtiendrons en particulier, grâce aux résultats sur les suites exactes de Spencer du deuxième chapitre, que le morphisme canonique
$\D ^\dag _{\X} (\hdag T) _{\Q} \otimes _{\D ^\dag _{\X^\#} (\hdag T) _{\Q}} ^{\L} \E (\ZZ)
  \rightarrow
  \D ^\dag _{\X} (\hdag T) _{\Q} \otimes _{\D ^\dag _{\X^\#} (\hdag T) _{\Q}} \E (\ZZ)$ est un isomorphisme
  (voir \ref{resspencerextlogdiag4}).

Dans une cinquième partie, nous définissons les foncteurs duaux, images directes
et images directes extraordinaires par $u$ en nous inspirant du cas non-logarithmique.
Grâce à l'isomorphisme ci-dessus (\ref{resspencerextlogdiag4}), on vérifie l'isomorphisme canonique :
$u _{T,+} (\E) \riso \D ^\dag _{\X} (\hdag T) _{\Q} \otimes  _{\D ^\dag _{\X^\#} (\hdag T) _{\Q}} \E (\ZZ)$,
où $u _{T,+}$ désigne l'image directe par $u$ à singularités surconvergentes le long de $T$.
En particulier, $u _{T,+} (\E)$ est réduit à un terme.
On en déduira alors que
$u _{T,+} (\E)$ est un $\D ^\dag _{\X} (\hdag T) _{\Q}$-module holonome
(via une extension de la notion d'holonomie sans structure de Frobenius).
Lorsque $T$ est vide et $\E$ est muni d'une structure de Frobenius, cela
signifie que $u _{+} (\E)$ est holonome au sens de Berthelot.

Nous prouvons dans une dernière partie que le complexe
de de Rham logarithmique de $\E$ est isomorphe au complexe de de Rham
de $u _{T,+} (\E)$.
En notant $\E (\hdag Z)$ l'isocristal sur $Y\setminus T$ surconvergent le long de $Z\cup T$ induit par $\E$,
on dispose d'un morphisme canonique $\rho _\E$ : $u _{T,+} (\E)\rightarrow \E (\hdag Z)$.
Nous établissons que si $\E$ est en fait un
isocristal sur $\X$ surconvergent le long de $T$ (e.g.,
le coefficient constant $\O _{\X} (\hdag T) _\Q$)
alors $\rho _\E$ est un isomorphisme.
La preuve utilise l'isomorphisme d'associativité du troisième chapitre qui ne fonctionne plus
si l'isocristal surconvergent est seulement un log isocristal surconvergent.
Moyennant quelques hypothèses sur les exposants le long des composantes irréductibles de $Z$,
on peut cependant conjecturer que $\rho _\E$ est un isomorphisme. Nous nous attaquerons à cette conjecture
dans un prochain travail.
Notons enfin que le fait que $\rho _\E$ soit un isomorphisme
implique que $\E (\hdag Z)$ est un $\D ^\dag _{\X} (\hdag T) _{\Q}$-module holonome et que
le complexe de de Rham logarithmique de $\E $
est isomorphe au complexe de de Rham de $\E (\hdag Z)$.

\subsection*{Remerciements :} Je remercie Y. Nakkajima et A. Shiho pour la suggestion
de considérer des diviseurs de $\X$ à la place de $X$.
Je remercie N. Tsuzuki pour les discussions
qui m'ont aidé à formuler une bonne version de la conjecture \ref{conj} et
qui ont suscité l'idée
d'étendre au cas des log isocristaux surconvergents ce qui avait été traité dans une première version.
Je remercie C. Noot-Huyghe pour une erreur décelée dans une précédente version.

\section*{Conventions et notations}

Nous conserverons les hypothèses concernant $\V$, $\mathfrak{m}$, $k$.
De plus, soit $\X$ un $\V$-schéma formel séparé et lisse.
Un sous-$\V$-schéma formel fermé $\ZZ$ de $\X$
est un diviseur si son idéal structural est localement principal.
De manière analogue à \cite[2.4]{dejong},
on dit que $\ZZ $ est un diviseur à croisements normaux strict s'il existe des diviseurs $\ZZ _1,\dots, \ZZ _e$ de $\X$ tels que
$\ZZ = \cup _{i=1,\dots, e} \ZZ _i$ (en tant que sous-$\V$-schémas formels fermés),
pour tout sous-ensemble non-vide $I \subset \{1,\dots, e\}$ le sous-$\V$-schéma formel $\cap _{i \in I} \ZZ _i$ soit
un $\V$-schéma formel lisse dont la fibre spéciale est de codimension $\# I$ dans celle de $\X$.
Dans tout cet article, fixons $\ZZ$ un diviseur à croisements normaux strict de $\X$.
Pour tout point $x$ de $\X$, il existe un ouvert $\U$ de $\X$ le contenant et muni des coordonnées locales
$t _1,\dots, t _d$ telles
que $\ZZ \cap \U = V (t_1\cdots t _s)$, avec $s \leq e$ (le nombre $s$ minimal que l'on peut obtenir correspond
au nombre de composantes irréductibles de $\ZZ$ contenant $x$). On remarque alors que le sous-monoïde
$\O _{\U} ^{\times} t _1 ^\N \cdots t _s ^\N\subset \O _{\U}$
est indépendant du choix de telles $t_1, \dots, t_d$ et ne dépend que de $\ZZ$.
Comme les limites projectives de monoïdes sont bien définies, il en résulte alors par
recollement (voir \cite[3.3]{EGAI}) un faisceau de monoïdes sur $\X$ noté $M (\ZZ)$ et une injection canonique $M (\ZZ) \hookrightarrow \O _{\X}$.
Notons $\Y:=\X \setminus \ZZ$ l'ouvert correspondant,
$\X ^\#$ le log $\V$-schéma formel lisse égal à $(\X, M (\ZZ))$.
De plus, soient $d = \dim X$,
$\S =\Spf \V $ et, pour tout entier $i\geq 0$,
$S _i := \mathrm{Spec} \, \V / \mathfrak{m} ^{i+1}$, $X _i := \X \times _{\S} S _i$,
$X ^\# _i := \X ^\# \times _{\S} S _i$,
$Z _i   := \ZZ \times _{\S} S _i$,
$Y _i := \Y \times _{\S} S _i$.
On remarque que la log structure de $X _i ^\#$ correspond à celle définie par
Y. Nakkajima et A. Shiho dans \cite[8.7-8.8]{nakkajimashiho}.

Lorsque l'on ne voudra pas distinguer le cas formel du cas algébrique,
on écrira $X $ (resp. $X ^\#$, $Y$, $Z$) à la place de $\X$ ou $X _i$ (resp. $\X ^\#$ ou $X _i ^\#$, resp. $\Y$ ou $Y _i$, resp. $\ZZ$ ou $Z _i$).
On désigne par $j $ : $ Y \subset X$ l'inclusion canonique.
On remarque que si $t _1,\dots, t _d$ sont des coordonnées locales de $X $
telles que $Z =V (t_1\cdots t _s)$ et $X ^\# = X ^\# \setminus V (t_{s+1} \cdots t _d)$
alors $t _1,\dots, t _d$ fournissent des coordonnées locales logarithmiques
sur $X ^\#$
(i.e., $\{d log (t _1),\dots, d log (t _d)\}$ est une base de $\Omega ^1 _{X ^\#/S}$).
Cette situation est toujours possible au voisinage de tout point de $X ^\#$.

{\it Convention :} Par la suite, les coordonnées locales logarithmiques $t _1,\dots, t _d\in M(\ZZ)$ de $X ^\#$ seront toujours supposées
construites de la façon précédente. En particulier, $t _1,\dots, t _d\in M(\ZZ)$ sont aussi des coordonnées locales de $X$.

Sauf mention explicite du contraire, tous les modules seront des modules à gauche et $m\geq 0$ sera
un entier fixé.
Si $\E$ est un faisceau abélien sur un espace topologique, $\E _\Q$ désignera le faisceau $\E \otimes _\Z \Q$
et $\widehat{\E}$ le complété de $\E$ pour la topologie $p$-adique.

Soit $\AA$ un faisceau d'anneaux sur $X$.
Si $*$ est l'un des symboles $\emptyset$, $+$, $-$, ou $\mathrm{b}$, on note $D ^* ( \AA )$
la catégorie dérivée des complexes de $\AA$-modules (à gauche) vérifiant les conditions correspondantes d'annulation
des faisceaux de cohomologie. Lorsque l'on souhaitera préciser entre droite et gauche, on écrira
$D ^* ( \overset{ ^g}{}\AA )$ ou $D ^* ( \AA \overset{ ^d}{})$.
Conformément aux notations et définitions de
\cite[5.2]{sga6},
on désignera par $D ^{\mathrm{b}} _{\mathrm{coh}} ( \AA )$
(resp.
$D_{\mathrm{parf}} ( \AA )$)
la sous-catégorie pleine de $D  ( \AA )$
dont les objets sont les
complexes à cohomologie cohérente et bornée
(resp.
les complexes parfaits).

Le symbole
$\underline{k}$ désigne le multi-indice
$(k _1,\dots, k_d) \in \N ^d$ et $|\underline{k}|=k_1+\dots + k_d$.
Par convention, si $x _1,\dots, x _d$ sont des éléments de l'idéal structural d'une $m$-PD-algèbre,
on notera
$\underline{x}  ^{\underline{k}}:=x _{1 } ^{k _1}\cdots x _{d} ^{k _d}$,
$\underline{x}  ^{\{\underline{k}\} _{(m)}}:=x _{1 } ^{\{k _1\} _{(m)}}\cdots x _{d} ^{\{k _d\} _{(m)}}$
où $x _i ^{\{k _i\} _{(m)}}$ est
la puissance divisée partielle de niveau $m$ d'ordre $k_i$ de l'élément $x_i$ (voir \cite[1.3]{Be1}).

Si $D$, $D '$ sont deux anneaux,
nous dirons que $E$ est un $(D , D')$-bimodule (resp. bimodule à gauche, resp. bimodule à droite) si $E$ est muni
de deux structures compatibles de $D$ -module à gauche (resp. à gauche, resp. à droite) et $D'$-module à droite
(resp. à gauche, resp. à droite). Si $D =D '$, nous dirons simplement $D$-bimodules (resp. bimodules à gauche,
resp. à droite).

  Si $\phi$ : $\AA \rightarrow \B$ est un homomorphisme de faisceaux d'anneaux, pour tout $\AA$-module $\M$
  (resp. $\B$-module $\NN$),
  on notera $\phi ^\flat (\M) := \mathcal{H} om _{\AA} (\B,\,\M)$
  (resp. $\phi _* (\NN)$ désigne $\NN$ vu comme $\AA$-module).
  La convention est de travailler dans
  les catégories dérivées mais les foncteurs de la forme $\phi ^\flat$ que nous utiliseront seront exacts.

\section{\label{nota1}Log-$\D$-modules arithmétiques}

\begin{nota}
\label{nota11}
Comme aucune confusion n'est possible,
nous omettrons toujours d'indiquer la base $S$ (munie de la structure logarithmique triviale).
Notons ainsi $\PP ^n _{X ^\# (m)} $ (au lieu de $\PP ^n _{X ^\# /S(m)} $)
le {\it voisinage à puissances divisées de niveau $m$ et d'ordre $n$ de $X ^\#$}
(voir \cite[2.2]{these_montagnon} pour le cas algébrique mais
pour le cas formel cela se définit de façon identique).
En dualisant celui-ci via sa structure gauche de $\O _X$-algèbre, on obtient
le {\it faisceau des opérateurs différentiels de niveau $m$ et d'ordre $n$ sur $X ^\#$} et noté
$\D ^{(m)} _{X ^\#,n}$. En prenant l'union sur les entiers $n$, on obtient
le {\it faisceau des opérateurs différentiels de niveau $m$ et d'ordre fini sur $X ^\#$} et noté
$\D ^{(m)} _{X ^\#}$ (voir \cite[2.3]{these_montagnon}).

Si $t _1,\dots, t _d$ sont des coordonnées locales logarithmiques de $X ^\#$ (voir ci-dessus les conventions adoptées à ce sujet) alors,
en notant $\tau _i = 1\otimes t _i - t_i \otimes 1$ et $\tau _{i\# } = \frac{1}{t_i} \tau _i$,
par \cite[2.2.1]{these_montagnon} (et de même pour les log-$\V$-schémas formels),
la famille $\underline{\tau}  ^{\{\underline{k}\} _{(m)}}$ (resp. $\underline{\tau} _{\#} ^{\{\underline{k}\} _{(m)}}$)
pour $|\underline{k}| \leq n$ forme une base
de $\PP ^n _{X  (m)} $ (resp. $\PP ^n _{X ^\# (m)} $) sur $\O _X$.
En dualisant, on obtient une base canonique de $\D ^{(m)} _{X }$ (resp. $\D ^{(m)} _{X ^\#}$) sur $\O _X$ que l'on notera
$\underline{\partial}  ^{<\underline{k}>_{(m)}}$
(resp. $\underline{\partial} _\# ^{<\underline{k}>_{(m)} }$).
Lorsque $\underline{k}$ est de la forme $\underline{k} = (0,\dots, 0,k,0,\dots, 0)$, où $k$ est à la $i$-ième place, nous noterons
$\partial _i ^{<k>_{(m)}}$
(resp. $\partial _{\#i} ^{<k>_{(m)} }$)
à la place de
$\underline{\partial}  ^{<\underline{k}>_{(m)}}$
(resp. $\underline{\partial} _\# ^{<\underline{k}>_{(m)} }$).
Rappelons que d'après \cite[Lemme 2.3.4]{these_montagnon}, pour tout $\underline{k} \in \N ^d$,
on dispose alors de la formule :
\begin{equation}
  \label{prodpartial}
  \underline{\partial} _\# ^{<\underline{k}>_{(m)} }= \prod _{i=1} ^d \partial _{\#i} ^{<k_i>_{(m)} },
\end{equation}
ce qui justifie nos notations.
Lorsque qu'il n'y aura aucune ambiguïté sur le niveau $m$, nous écrirons simplement
$\underline{\partial}  ^{<\underline{k}> }$ et $\underline{\partial} _\# ^{<\underline{k}> }$.
On calcule la relation :
$\underline{\partial} _\# ^{<\underline{k}>}= \underline{t} ^{\underline{k}} \underline{\partial} ^{<\underline{k}>}
(=t _1 ^{k_1 } \cdots t _d ^{k_d} \underline{\partial} ^{<\underline{k}>} ) $.

\end{nota}

Dans la section \cite[2.6]{these_montagnon}, Montagnon définit les $m$-PD-stratifications logarithmiques et montre
qu'une structure de $\D _{X ^\# } ^{(m)}$-module à gauche prolongeant une structure de $\O _{X}$-module
est équivalente à celle d'une $m$-PD-stratification.
Nous allons maintenant étendre cette notion de $m$-PD-stratification en remplaçant le faisceau $\O _X$ par
un faisceau de $\O _X$-algèbres comme suit.

\begin{defi}
\label{defBcomp}
De manière analogue à \cite[2.3]{Be1},
une $\O _X$-algèbre commutative $\B _X$ est dite munie d'une structure de $\D ^{(m)}_{X^\#}$-module à gauche
{\it compatible} à sa structure de $\O _X$-algèbre si :
\begin{enumerate}
  \item \label{defBcompi} La structure de $\O _X$-algèbre est égale à la structure de $\O _X$-module
sous-jacente à la structure de $\D _{X^\#} ^{(m)}$-module.
\item \label{defBcompii} Les isomorphismes $\epsilon _n ^{\B _X}$ : $ \PP ^n _{X ^\#(m)} \otimes _{\O _X } \B _X  \riso \B _X \otimes _{\O _X
}\PP ^n _{X ^\# (m)} $ de la $m$-PD-stratification  correspondante sont
des isomorphismes de $ \PP ^n _{X ^\# (m)}$-algèbres.
\end{enumerate}
\end{defi}

\begin{rema}
\label{remaLeibnitz}
De manière analogue à \cite[2.3.4.1]{Be1},
en coordonnées locales, lorsque la condition \ref{defBcomp}.\ref{defBcompi} est vérifiée,
la condition \ref{defBcomp}.\ref{defBcompii} équivaut à la vérification de la formule de Leibnitz
\begin{equation}
  \label{Leibnitz}
  \forall f,\, g \in \B _X,\ \forall \underline{k} \in \N ^d, \hfill
\underline{\partial} _\# ^{<\underline{k}>} ( f g)=
\sum _{\underline{h}\leq \underline{k}}
\left \{ \begin{smallmatrix}  \underline{k} \\  \underline{h} \\ \end{smallmatrix} \right \}
\underline{\partial} _\# ^{<\underline{k}-\underline{h}>} ( f)
\underline{\partial} _\# ^{<\underline{h}>} ( g),
\end{equation}
où, pour éviter les confusions, on désigne par $\underline{\partial} _\# ^{<\underline{k}>} ( fg)$
l'action comme $\smash{\widetilde{\D}} _{X^\#} ^{(m)}$-module de $\underline{\partial} _\# ^{<\underline{k}>}$
sur $fg$ etc.
\end{rema}

\begin{exem}
\label{exB(T)}
  Soient $T$ un diviseur de $X _0$ et $r$ un multiple de $p ^{m+1}$. Berthelot a construit dans \cite[4.2.3]{Be1} un
faisceau de $\O_X$-algèbres commutatives noté $\B _X  (T,r)$, muni
d'une structure compatible de $\D ^{(m)}_X$-module à gauche. La $m$-PD-stratification
associée induit alors par extension via $\PP ^n _{X  (m)}\rightarrow  \PP ^n _{X ^\# (m)}$
une $m$-PD-stratification logarithmique relativement à $X ^\#$. On obtient ainsi une structure
canonique de $\D ^{(m)}_{X^\#}$-module à gauche sur $\B _X (T,r)$ compatible à sa structure de $\O _X$-algèbre.
\end{exem}

\begin{vide}

Fixons pour toute la suite d'une part $\B _X$ une $\O _X$-algèbre commutative munie d'une structure
de $\D ^{(m)} _{X^\#}$-module à gauche compatible
à sa structure de $\O _X$-algèbre et
d'autre part $\B _{X} '$ une $\B _{X}$-algèbre commutative munie d'une
structure compatible de $\D _{X^\#} ^{(m)}$-module à gauche
telle que $\B _{X} \rightarrow \B _{X}  '$ soit $\D ^{(m)} _{X^\#}$-linéaire.

La construction qui suit (i.e., de $\widetilde{\D} ^{(m)} _{X^\#}$) reprend (de façon plus concise) celle du cas non logarithmique
(voir \cite[1.1]{caro_comparaison}). Le lecteur avide de détails pourra s'y reporter.

On notera par la suite $\widetilde{\PP} ^n _{X ^\# (m)} $ le faisceau de $\PP ^n _{X ^\# (m)} $-algèbres
$\B _X \otimes _{{\O }_X } \PP ^n _{X ^\# (m)} $ (par convention, vu sa position à droite,
rappelons que l'on choisit la structure gauche de $\O _X$-algèbre de
$\PP ^n _{X ^\# (m)}$ pour calculer $\B _X \otimes _{{\O }_X } \PP ^n _{X ^\# (m)}$). Le morphisme canonique $\tilde{d} ^n _0$ : $\B _X \rightarrow\B _X
\otimes _{\O _X }\PP ^n _{X ^\# (m)} $ munit le faisceau $\widetilde{\PP} ^n _{X ^\# (m)} $ d'une structure de
$\B_X$-algèbre que l'on appellera structure {\it gauche}. De plus, le morphisme de $\B _X$-algèbres
$$\tilde{d} _1 ^n \ : \ \B _X \rightarrow
\PP ^n _{X ^\# (m)} \otimes _{\O _X } \B _X  \underset{\epsilon _n ^{\B _X}}{\riso} \B _X \otimes _{\O _X
}\PP ^n _{X ^\# (m)} $$ induit une deuxième structure de $\B _X$-algèbre sur $\widetilde{\PP} ^n _{X ^\# (m)} $ que l'on
appellera structure {\it droite}.

Pour tout couple d'entiers positifs, $n$ et $n'$, le faisceau $\PP ^n _{X ^\# (m)} \otimes _{\O _X} \PP ^{n'} _{X ^\# (m)}$
possède trois structures de $\O _X$-algèbres. La structure de $\O _X$-algèbre de $\PP ^n _{X ^\# (m)}
\otimes _{\O _X} \PP ^{n'} _{X ^\# (m)}$ provenant de la structure gauche de $\PP ^n _{X ^\# (m)} $ sera appelée structure
{\it gauche}, celle provenant de la structure "produit tensoriel" sera dénommée structure du
{\it centre} et enfin celle découlant de la structure droite de $\PP ^{n'} _{X ^\# (m)}$ sera dite structure {\it droite}.

Soit le faisceau $\widetilde{\PP} ^n _{X ^\# (m)} \otimes _{\B _X} \widetilde{\PP} ^{n'} _{X ^\# (m)}$,
où, pour calculer le produit tensoriel, on utilise
la structure droite de $\widetilde{\PP} ^n _{X ^\# (m)} $
et la structure gauche de $\widetilde{\PP} ^{n'} _{X ^\# (m)}$.
On munit $\widetilde{\PP} ^n _{X ^\# (m)} \otimes _{\B _X} \widetilde{\PP} ^{n'} _{X ^\# (m)}$,
de manière analogue à $\PP ^n _{X ^\# (m)} \otimes _{\O _X} \PP ^{n'} _{X ^\# (m)}$,
de structures gauche, du centre et droite de $\B _X$-algèbre.
On désigne par
$\tilde{d} ^{n,n'} _0$,
$\tilde{d} ^{n,n'} _1$,
$\tilde{d} ^{n,n'} _2$
les homomorphismes
$\B _{X } \rightarrow \widetilde{\PP} ^n _{X ^\# (m)} \otimes _{\B _X } \widetilde{\PP} ^{n'} _{X ^\# (m)}$
munissant $\widetilde{\PP} ^n _{X ^\# (m)} \otimes _{\B _X } \widetilde{\PP} ^{n'} _{X ^\# (m)}$ de sa structure de
$\B _X$-algèbre à gauche, du centre, à droite.

On introduit les homomorphismes de $\B _X$-algèbres suivants
\begin{gather}
\widetilde{\delta} ^{n,n'} _{(m)} \ : \
\xymatrix {
{\widetilde{\PP} ^{n+n'} _{X ^\# (m)}}
\ar[rr] ^(0.4){\B _X \otimes \delta ^{n,n'} _{(m)}}
&&
{\widetilde{\PP} ^n _{X ^\# (m)} \otimes _{\B _X} \widetilde{\PP} ^{n'} _{X ^\# (m)}}
},
\notag \\
\widetilde{q _0} ^{n,n'}  \ : \ \widetilde{\PP} ^{n+n'} _{X ^\# (m)} \longrightarrow \widetilde{\PP} ^n _{X ^\# (m)}
\longrightarrow \widetilde{\PP} ^n _{X ^\# (m)} \otimes _{\B _X} \widetilde{\PP} ^{n'} _{X ^\# (m)}, \notag \\
\widetilde{q _1} ^{n,n'} \ : \ \widetilde{\PP} ^{n+n'} _{X ^\# (m)} \longrightarrow \widetilde{\PP} ^{n'} _{X ^\# (m)}
\longrightarrow \widetilde{\PP} ^n _{X ^\# (m)} \otimes _{\B _X} \widetilde{\PP} ^{n'} _{X ^\# (m)}, \notag
\end{gather}
où
$\delta ^{n,n'} _{(m)}$ :
$\PP ^{n+n'} _{X ^\#,\,(m)}\rightarrow \PP ^n _{X ^\#,\,(m)} \otimes _{\O _X} \PP ^{n'} _{X ^\#,\,(m)}$
désigne l'homomorphisme construit par Montagnon dans \cite[2.3.2.A]{these_montagnon}
(dans le cas algébrique, mais la construction formelle est identique).

\end{vide}
\begin{vide}
\label{def1botimesd}
On notera $\smash{\widetilde{\D}} _{X ^\#,\,n} ^{(m)} := \mathcal{H}om _{\B _X} (\tilde{d} ^n _{0*} \widetilde{\PP} ^n _{X ^\# (m)}  , \B _X)$
le dual $\B _X$-linéaire pour la structure gauche de
$\widetilde{\PP} ^n _{X ^\# (m)}$.
Pour $n' \geq n$,
les surjections $\widetilde{\PP}  ^{n'} _{X ^\# (m)} \twoheadrightarrow \widetilde{\PP} ^n _{X ^\# (m)}$
fournissent des injections
$\smash{\widetilde{\D}} _{X ^\#,\,n} ^{(m)} \hookrightarrow \smash{\widetilde{\D}} _{X ^\#,\,n'} ^{(m)}.$
On munit le faisceau $\smash{\widetilde{\D}} _{X^\#} ^{(m)}:=\cup _{n \in \N}\smash{\widetilde{\D}} _{X ^\#,\,n} ^{(m)}$
d'une structure d'anneau grâce aux accouplements
 $ \smash{\widetilde{\D}} _{X ^\#,\,n} ^{(m)} \times \smash{\widetilde{\D}} _{X ^\#,\,n'} ^{(m)}
 \rightarrow \smash{\widetilde{\D}} _{X ^\#,\,n+n'} ^{(m)}$
définis comme suit : si $\widetilde{P} \in \smash{\widetilde{\D}} _{X ^\#,\,n} ^{(m)}$, $\widetilde{P} '\in \smash{\widetilde{\D}} _{X ^\#,\,n'}
^{(m)}$, alors $\widetilde{P} \cdot \widetilde{P}'$ est l'homomorphisme composé
\begin{equation}\label{deltatildeann}
  \widetilde{\PP} ^{n+n'} _{X ^\# (m)} \overset{\widetilde{\delta} ^{n,n'} _{(m)}}{\longrightarrow}
\widetilde{\PP} ^n _{X ^\# (m)} \otimes _{\B _X} \widetilde{\PP} ^{n'} _{X ^\# (m)} \overset{Id \otimes
\widetilde{P}'}{\longrightarrow} \widetilde{\PP} ^n _{X ^\# (m)} \overset{ \widetilde{P}}{\longrightarrow} \B _X.
\end{equation}
Cette loi de composition est bien {\it associative}.

En procédant de manière analogue
à \cite[2.3.5]{Be1}, on munit le faisceau $\B _X \otimes _{\O _X} \D _{X^\#} ^{(m)}$
d'une structure de $\B _X$-algèbre en définissant le produit
$P'':= P\cdot P'$ de $P \in \mathcal{H}om _{\O _X} (\PP ^n _{X ^\# (m)} , \B _X )$ et de
$P '\in \mathcal{H}om _{\O _X} (\PP ^{n'} _{X ^\# (m)} , \B _X )$ comme le composé
$$ \PP ^{n+n'} _{X ^\# (m)} \overset{\delta ^{n,n'} _{(m)}}{\longrightarrow}
{\PP} ^n _{X ^\# (m)} \otimes _{\O _X} {\PP} ^{n'} _{X ^\# (m)}\overset{Id \otimes P '}{\longrightarrow} {\PP} ^n _{X ^\# (m)}
\otimes _{\O _X} \B _X \overset{\epsilon _n ^{\B _X}}{\longrightarrow} \B _X \otimes _{\O _X} {\PP} ^n _{X ^\# (m)}
\overset{ Id \otimes P }{\longrightarrow} \B _X \otimes _{\O _X} \B _X \overset{\mu}{\longrightarrow} \B _X,$$
où $\mu$ est la multiplication canonique de $\B _X$.

Comme pour \cite[1.1.13]{caro_comparaison}, on vérifie que l'isomorphisme $\B _X$-linéaire canonique
$\B _X \otimes _{\O _X} \D _{X^\#} ^{(m)} \riso \smash{\widetilde{\D}} _{X^\#} ^{(m)}$
est un isomorphisme de $\B _X$-algèbres.
Enfin, si $t _1 ,\dots ,t _d$, sont des coordonnées logarithmiques locales,
avec les notations de \ref{nota11},
si aucune confusion n'est à craindre, on notera
encore
$\underline{\tau} _{\#} ^{\{\underline{k}\} _{(m)}}$
(resp. $\underline{\partial} _{\#} ^{<\underline{k}> _{(m)}}$)
l'élément de $\widetilde{\PP} ^n _{X ^\# (m)}$ (resp. $\widetilde{\D} ^{(m)} _{X ^\#}$)
à la place de $1\otimes \underline{\tau} _{\#} ^{\{\underline{k}\} _{(m)}}$
(resp. $1\otimes \underline{\partial} _{\#} ^{<\underline{k}> _{(m)}}$).

\end{vide}

\begin{defi}
\label{definstrat} Soit $\E$ un $\B _X$-module. Une \textit{m-PD-stratification} (ou \textit{PD-stratification de
niveau m}) $\widetilde{\epsilon} $ sur $\E$ {\it relativement à $(X ^\#/S,\,\B _X)$}
(ou {\it relativement à $\B _X$} si aucune confusion n'est à craindre)
est la donnée d'une famille compatible
d'isomorphismes $\widetilde{\PP} ^n _{X ^\# (m)}$-linéaires
$$ \widetilde{\epsilon}^\E  _n \ :\ \widetilde{\PP} ^n _{X ^\# (m)}\otimes _{\B _X} \E \riso
\E \otimes _{\B _X} \widetilde{\PP} ^n _{X ^\# (m)},$$ où les produits tensoriels sont respectivement pris pour les
structures droite et gauche de $\widetilde{\PP} ^n _{X ^\# (m)}$, ces isomorphismes étant astreints aux conditions
suivantes :
\begin{enumerate}
\item $\widetilde{\epsilon} ^\E_0 = Id _{\E}$ ;
\item La condition de cocycle est validée, i.e., pour tous $n$, $n'$,
le diagramme
$$\xymatrix {
{\widetilde{\PP} ^n _{X ^\# (m)} \otimes _{\B _X} \widetilde{\PP} ^{n'} _{X ^\# (m)}\otimes _{\B _X} \E }
\ar[rr]^{\widetilde{\delta} ^{n,n'*} _{(m)}(\widetilde{\epsilon} ^\E _{n +n'})}
_{\widetilde{{}\hspace{0,3cm}  }}
 \ar[rd]^{\tilde{q} _1 ^{n,n'*}(\widetilde{\epsilon} ^\E_{n +n'})} _{\widetilde{{}\hspace{0,3cm}  }}
 &
 {}
 &
{\E \otimes _{\B _X} \widetilde{\PP} ^n _{X ^\# (m)} \otimes _{\B _X} \widetilde{\PP} ^{n'} _{X ^\# (m)}} \\
{}
&
{\widetilde{\PP} ^n _{X ^\# (m)} \otimes _{\B _X} \E \otimes _{\B _X} \widetilde{\PP} ^{n'} _{X ^\# (m)}}
\ar[ur] ^{\tilde{q} _0 ^{n,n'*}(\widetilde{\epsilon} ^\E _{n +n'})}_{\widetilde{{}\hspace{0,3cm}  }}
&
{}
}$$ est commutatif.
\end{enumerate}

\end{defi}
\begin{prop}
\label{eqDtilde}
Pour tout $\B '_X$-module $\E$, il y a équivalence entre les données suivantes :

a) Une structure de $\B '_X \otimes _{\O _X} \D _{X ^\#}^{(m)}$-module à gauche sur $\E$
prolongeant sa structure canonique de $\B '_X$-module ;

b) Une $m$-PD-stratification $\widetilde{\epsilon} ^{\prime \E}= (\widetilde{\epsilon} ^{\prime \E}_n )$ sur $\E$
relative à $\B '_X$.

Notons alors $(\widetilde{\epsilon} _n ^{\B '_X})$ la $m$-PD-stratification relative à $\B _X$ de $\B'_X$.
Les données $a)$ et $b)$ sont équivalentes à la suivante :

c) Une $m$-PD-stratification $\widetilde{\epsilon} ^{\E}= (\widetilde{\epsilon} ^{\E}_n )$ sur $\E$
relative à $\B _X$ dont les
isomorphismes $\widetilde{\epsilon} ^{\E} _n$
sont semi-linéaires par rapport aux isomorphismes $\widetilde{\epsilon} _n ^{\B '_X}$ ;

De plus, un homomorphisme $\B '_X$-linéaire
$\E \rightarrow \FF$ entre deux
$\B '_X \otimes _{\O _X} \D _{X ^\#}^{(m)}$-modules à gauche est
$\B '_X \otimes _{\O _X} \D _{X ^\#}^{(m)}$-linéaire si et seulement s'il commute aux isomorphismes
$\widetilde{\epsilon} ^{\E} _n$ et $\widetilde{\epsilon} ^{\FF} _n$
(resp. $\widetilde{\epsilon} ^{\prime \E} _n$ et $\widetilde{\epsilon} ^{\prime \FF} _n$).
Le morphisme sera dit {\it horizontal}.

Enfin, en coordonnées locales, pour toute section $x$ de $\E$, on a la formule :
\begin{equation}
\label{taylor}
\widetilde{\epsilon} ^{\prime \E} _n ( 1\otimes x )= \sum _{|\underline{k}|\leq n}\underline{\partial}
_\# ^{<\underline{k}>} x \otimes \underline{\tau} _\# ^{\{ \underline{k}\} }.
\end{equation}

\end{prop}

\begin{proof}
  Identique à celle de \cite[1.1.16]{caro_comparaison}.
\end{proof}

\begin{exem}
  \label{BXDmod}
On munit $\B _X$ de la $m$-PD-stratification triviale, i.e., les isomorphismes
de la $m$-PD-stratification sont les identités de $\widetilde{\PP} ^n _{X ^\# (m)}$.
La structure canonique de $\smash{\widetilde{\D}} _{X^\#} ^{(m)}$ induite sur $\B _X$ est
alors décrite de la façon suivante :
si $P \in \smash{\widetilde{\D}} _{X^\#} ^{(m)}$ et $b \in \B _X$, l'action de $P$ sur $b$
est l'image de $b$ via le morphisme composé :
$\B _X \overset{\tilde{d} _1 ^n}{\longrightarrow} \widetilde{\PP} ^n _{X ^\# (m)}
\overset{P}{\longrightarrow} \B _X $ (on le voit par exemple par $\B _X$-linéarité via \ref{taylor}).
On obtient ainsi, pour tout $b \in \B _X$, la formule analogue à \cite[2.2.4.(i)]{Be1} :
$\tilde{d} _1 ^n (b) = \sum _{|\underline{k}|\leq n}
(\underline{\partial} _\# ^{<\underline{k}>} ( b) ) \underline{\tau} _\# ^{\{ \underline{k}\} }$,
où, pour éviter les confusions, on désigne par $\underline{\partial} _\# ^{<\underline{k}>} ( b)$
l'action comme $\smash{\widetilde{\D}} _{X^\#} ^{(m)}$-module de $\underline{\partial} _\# ^{<\underline{k}>}$
sur $b$
(et non la multiplication via la structure droite de $\B _X$-algèbre de $\smash{\widetilde{\D}} _{X^\#} ^{(m)}$).
Il en résulte la formule analogue à \cite[2.2.4.(iv)]{Be1} ou \cite[2.3.5.1]{Be1} :
\begin{equation}
  \label{be1224iv}
\underline{\partial} _\# ^{<\underline{k}>}  b =
\sum _{\underline{h}\leq \underline{k}}
\left \{ \begin{smallmatrix}  \underline{k} \\  \underline{h} \\ \end{smallmatrix} \right \}
\underline{\partial} _\# ^{<\underline{k}-\underline{h}>} ( b)
\underline{\partial} _\# ^{<\underline{h}>} .
\end{equation}

\end{exem}

De manière analogue à \cite[1.1.3]{Be2} (ou \cite[1.2.22]{caro_comparaison}), on définit
les $m$-PD-{\it costratifications relativement à $\B _X$} :

\begin{defi}\label{defincostrat}
  Soit $\M$ un $\B _X$-module. Une $m$-PD-{\it costratification sur $\M$ relativement}
  à $(X ^\# /S,\,\B _X)$ (ou $\B _X$ si aucune ambiguïté n'est à craindre)
  est la donnée d'une famille {\it compatible} d'isomorphismes $\widetilde{\PP} ^n _{X ^\# (m)}$-linéaires
  $$ \widetilde{\epsilon} ^\M _n \ :\
  \mathcal{H} om _{\B _X} ( \tilde{d} ^{n} _{0*} \widetilde{\PP} ^n _{X ^\# (m)},\,\M)
  \riso
  \mathcal{H} om _{\B _X} ( \tilde{d} ^{n} _{1*} \widetilde{\PP} ^n _{X ^\# (m)},\,\M),$$
  ceux-ci vérifiant les conditions suivantes :
  \begin{enumerate}
    \item $\widetilde{\epsilon} ^\M _0 = Id _{\M}$ ;
    \item Pour tous $n$, $n'$, le diagramme
    \begin{equation}
      \label{definstratdiag1}
      \xymatrix@C=-1cm {
 {  \mathcal{H} om _{\B _X}
  ( \tilde{d} ^{n,n'} _{0*}( \widetilde{\PP} ^n _{X ^\# (m)} \otimes _{\B _X} \widetilde{\PP} ^{n'} _{X ^\# (m)}),\,\M)}
  \ar[rr]^{\widetilde{\delta} ^{n,n'\flat} _{(m)}(\widetilde{\epsilon} ^\M _{n +n'})}
  _{\widetilde{{}\hspace{0,3cm}  }}
   \ar[rd]^{\tilde{q} _0 ^{n,n'\flat}(\widetilde{\epsilon} ^\M _{n +n'})}
   _{\widetilde{{}\hspace{0,3cm}  }}
   &&
 {  \mathcal{H} om _{\B _X}
  ( \tilde{d} ^{n,n'} _{2*}( \widetilde{\PP} ^n _{X ^\# (m)} \otimes _{\B _X} \widetilde{\PP} ^{n'} _{X ^\# (m)}),\,\M)}
\\
&
 {  \mathcal{H} om _{\B _X}
  ( \tilde{d} ^{n,n'} _{1*} (\widetilde{\PP} ^n _{X ^\# (m)} \otimes _{\B _X} \widetilde{\PP} ^{n'} _{X ^\# (m)}),\,\M)}
  \ar[ur] ^{\tilde{q} _1 ^{n,n'\flat}(\widetilde{\epsilon} ^\M _{n +n'})}
  _{\widetilde{{}\hspace{0,3cm}  }}
  }
    \end{equation}
    est commutatif.
  \end{enumerate}
\end{defi}

\begin{prop}
\label{Dmodstrat}
Pour tout $\B '_X$-module $\M$, il y a équivalence entre les données suivantes :

  a) Une structure de $\B '_X \otimes _{\O _X} \D _{X ^\#}^{(m)}$-module à droite sur $\M$
  prolongeant sa structure de $\B '_X$-module ;

  b) Une $m$-PD-costratification $( \widetilde{\epsilon} _n ^{\M})$ relativement à $\B _X$ sur $\M$
telle que les isomorphismes $\widetilde{\epsilon} _n ^{\M}$ soient
semi-linéaires par rapport à $(\widetilde{\epsilon} _n ^{\B '_X})^{-1}$ ;

  c) Une $m$-PD-costratification $( \widetilde{\epsilon} _n ^{\prime \M})$ relativement à $\B '_X$ sur $\M$.

Un homomorphisme $\B ' _X$-linéaire $\M \rightarrow \NN$ entre deux
$\B '_X \otimes _{\O _X} \D _{X ^\#}^{(m)}$-modules à droite est
$\B '_X \otimes _{\O _X} \D _{X ^\#}^{(m)}$-linéaire si et seulement s'il commute aux isomorphismes
$\widetilde{\epsilon} _n ^{\M}$ et $\widetilde{\epsilon} _n ^{\NN}$
(resp. $\widetilde{\epsilon} _n ^{\prime \M}$ et $\widetilde{\epsilon} _n ^{\prime \NN}$).

\end{prop}

\begin{proof}
    Identique à celle de \cite[1.1.23]{caro_comparaison}.
\end{proof}

\begin{prop}\label{MotimesEetc}
  Soient $\E$, $\FF$ deux $\widetilde{\D} ^{(m)} _{X ^\#}$-modules à gauche,
  $\M$, $\NN$ deux $\widetilde{\D} ^{(m)} _{X ^\#}$-modules à droite.
La structure de
  $\B _X$-module de $\mathcal{H} om _{\B _X} (\NN,\M )$
  (resp. $\E \otimes _{\B _X} \FF$, resp. $\mathcal{H} om _{\B _X} (\E,\FF )$)
  se prolonge en une structure canonique de
$\widetilde{\D} ^{(m)} _{X ^\#}$-module à gauche.
  La structure de
  $\B _X$-module de $\M \otimes _{\B _X} \E$
  (resp. $\mathcal{H} om _{\B _X} (\E, \M)$)
  se prolonge en une structure canonique de
$\widetilde{\D} ^{(m)} _{X ^\#}$-module à droite.
\end{prop}
\begin{proof}
Cela se vérifie comme dans le cas non logarithmique (\cite[1.1.18 et 1.1.24]{caro_comparaison}).
Pour mémoire, via pour $i=0,1$ les identifications
\begin{gather}
(\E \otimes _{\B _X} \FF ) \otimes _{\B _X} {d} ^{n} _{i*} \widetilde{\PP} ^n _{X ^\#,\,(m)}
\riso
({d} ^{n} _{i*} \widetilde{\PP} ^n _{X ^\#,\,(m)}  \otimes _{\B _X} \E )
\otimes _{\widetilde{\PP} ^n _{X ^\#,\,(m)} }
({d} ^{n} _{i*} \widetilde{\PP} ^n _{X ^\#,\,(m)}  \otimes _{\B _X} \FF ),
\\
\mathcal{H} om _{\B _X} (\E,\FF ) \otimes _{\B _X} {d} ^{n} _{i*} \widetilde{\PP} ^n _{X ^\#,\,(m)}
\riso
\mathcal{H} om _{\widetilde{\PP} ^n _{X ^\#,\,(m)}}
(\E \otimes _{\B _X} {d} ^{n} _{i*} \widetilde{\PP} ^n _{X ^\#,\,(m)},
\FF \otimes _{\B _X} {d} ^{n} _{i*} \widetilde{\PP} ^n _{X ^\#,\,(m)}),
\\
\mathcal{H} om _{\B _X} (\NN,\M ) \otimes _{\B _X} {d} ^{n} _{i*} \widetilde{\PP} ^n _{X ^\#,\,(m)}
\riso
\mathcal{H} om _{\widetilde{\PP} ^n _{X ^\#,\,(m)}}
(\mathcal{H} om _{\B _{X}} ( {d} ^{n} _{i*} \widetilde{\PP} ^n _{X ^\#,\,(m)}, \NN),
\mathcal{H} om _{\B _{X}} ( {d} ^{n} _{i*} \widetilde{\PP} ^n _{X ^\#,\,(m)}, \M)),
\\
 \mathcal{H} om _{\B _{X}} ( {d} ^{n} _{i*} \widetilde{\PP} ^n _{X ^\#,\,(m)},\,\M \otimes _{\B _X} \E )
 \riso
 \mathcal{H} om _{\B _{X}} ( {d} ^{n} _{i*} \widetilde{\PP} ^n _{X ^\#,\,(m)},\,\M )
 \otimes _{\widetilde{\PP} ^n _{X ^\#,\,(m)} }
({d} ^{n} _{i*} \widetilde{\PP} ^n _{X ^\#,\,(m)}  \otimes _{\B _X} \E ), \\
  \mathcal{H} om _{\B _{X}} ( {d} ^{n} _{i*} \widetilde{\PP} ^n _{X ^\#,\,(m)},
\mathcal{H} om _{\B _X} (\E, \M))
\riso
\mathcal{H} om _{\widetilde{\PP} ^n _{X ^\#,\,(m)}} (\E \otimes _{\B _X} {d} ^{n} _{i*} \widetilde{\PP} ^n _{X ^\#,\,(m)},
\mathcal{H} om _{\B _{X}} ( {d} ^{n} _{i*} \widetilde{\PP} ^n _{X ^\#,\,(m)},  \M)),
\end{gather}
on obtient les $m$-PD-stratifications cherchées en posant
$\epsilon _n ^{\E \otimes \FF} : =
\epsilon _n ^{\E } \otimes _{\widetilde{\PP} ^n _{X ^\#,\,(m)} } \epsilon _n ^{\FF}$,
$\epsilon _n ^{ \mathcal{H} om _{\B _X} (\E, \FF) }: =
\mathcal{H} om _{\widetilde{\PP} ^n _{X ^\#,\,(m)} } ((\epsilon _n ^{\E })^{-1}, \epsilon _n ^{\FF })$,
$\epsilon _n ^{ \mathcal{H} om _{\B _X} (\NN, \M) }: =
\mathcal{H} om _{\widetilde{\PP} ^n _{X ^\#,\,(m)} } (\epsilon _n ^{\NN }, (\epsilon _n ^{\M })^{-1})$,
$\epsilon _n ^{\M \otimes \E} : =
\epsilon _n ^{\M } \otimes _{\widetilde{\PP} ^n _{X ^\#,\,(m)} } (\epsilon _n ^{\E}) ^{-1}$,
$\epsilon _n ^{\mathcal{H} om _{\B _X} (\E, \M) } :=
\mathcal{H} om _{\widetilde{\PP} ^n _{X ^\#,\,(m)} } (\epsilon _n ^{\E }, \epsilon _n ^{\M })$.

\end{proof}

\begin{rema}\label{rem261mon}
Avec les notations de \ref{MotimesEetc},
contrairement au cas non logarithmique, le calcul
de l'inverse de $\epsilon _n ^{\E}$ à partir de la formule tautologique dit du développement de Taylor \ref{taylor} n'est pas aisé.
En particulier, la formule de \cite[2.3.2.3]{Be1} :
\begin{equation}
  \label{rem261mon=}
(\epsilon _n  ) ^{-1} (x \otimes 1) =
\sum _{\underline{k} \leq n} (-1) ^{|\underline{k}|} \underline{\tau} ^{\{\underline{k}\}}
\otimes \underline{\partial}  ^{<\underline{k}>}x,
\end{equation}
est (en général) fausse si on remplace respectivement {\og $\underline{\tau} ^{\{\underline{k}\}}$\fg}
par {\og $\underline{\tau} ^{\{\underline{k}\}}_\#$\fg} et
{\og $\underline{\partial} ^{<\underline{k}>}$\fg}
par {\og $\underline{\partial} ^{<\underline{k}>}_\#$\fg}.
En effet, la formule
\cite[Lemme 2.3.4]{these_montagnon} implique que l'égalité
$\underline{\partial}  ^{<\underline{k}>} \underline{\partial}  ^{<\underline{h}>}
=
\left <\begin{smallmatrix}  \underline{k} + \underline{h} \\  \underline{k} \\ \end{smallmatrix} \right >
\underline{\partial}  ^{<\underline{k}+\underline{h}>}$
(voir \cite[2.2.4.(iii)]{Be1}) devient fausse avec des dièses (mais elle devient vraie dans
$\mathrm{gr} \widetilde{\D} ^{(m)} _{X ^\#}$ : voir
\ref{lemm-wildetildeDcoh}.\ref{lemm-wildetildeDcohi}).
Ainsi, les analogues logarithmiques des formules \cite[2.3.3.2]{Be1},
\cite[1.1.7.1-3]{Be2} décrivant l'action de $\underline{\partial} ^{<\underline{k}>}_\#$ sur
$\mathcal{H} om _{\B _X} (\E,\FF )$, $\mathcal{H} om _{\B _X} (\E,\M )$, $\mathcal{H} om _{\B _X} (\NN,\M )$
sont faux (i.e., il ne suffit pas d'ajouter des dièses). Néanmoins on dispose des formules \ref{MotimesEetcegpre},
\ref{MotimesEetceg} et \ref{MotimesEetcegbis} ci-dessous (cette dernière nous servira pour obtenir le lemme \ref{MotimesEetcegcons}
du théorème \ref{A2CN}).
\end{rema}

\begin{rema}
  \label{prodtens2alg}
  Soient $\CC _{X}$ (resp. $\CC _{X} '$) une $\B _{X}$-algèbre commutative munie d'une
structure compatible de $\D _{X^\#} ^{(m)}$-module à gauche
telle que $\B _{X} \rightarrow \CC _{X}$ (resp. $\B _{X} \rightarrow \CC _{X}  '$)
soit $\D ^{(m)} _{X^\#}$-linéaire.
La structure produit tensoriel de $\D _{X^\#} ^{(m)}$-module à gauche sur
la $\B_X$-algèbre $\CC _X \otimes _{\B _X} \CC ' _X$ est compatible. On bénéficie de plus
des morphismes d'algèbres $\widetilde{\D} _{X^\#} ^{(m)}$-linéaires
$\CC _X \rightarrow \CC _X \otimes _{\B _X} \CC ' _X$ et
$\CC _X '\rightarrow \CC _X \otimes _{\B _X} \CC ' _X$.
\end{rema}

\begin{lemm}
\label{comm-assoc}
  Soient $\E$, $\FF$, $\G$ des $\widetilde{\D} ^{(m)} _{X ^\#}$-modules à gauche et
  $\M$ un $\widetilde{\D} ^{(m)} _{X ^\#}$-module à droite.
  Les isomorphismes canoniques de commutation et d'associativité
  \begin{gather}
    (\E \otimes _{\B _{X}} \FF )\otimes _{\B _{X}} \G
    \riso
    \E \otimes _{\B _{X}} (\FF \otimes _{\B _{X}} \G),
    \\ \label{comm-assoc2}
    (\M \otimes _{\B _{X}} \FF )\otimes _{\B _{X}} \G
    \riso
    \M \otimes _{\B _{X}} (\FF \otimes _{\B _{X}} \G),
    \\
    \E \otimes _{\B _{X}} \FF \riso \FF \otimes _{\B _{X}} \E
  \end{gather}
  sont $\widetilde{\D} ^{(m)} _{X ^\#}$-linéaires.
\end{lemm}
\begin{proof}
  On vérifie que les morphismes sont horizontaux.
\end{proof}

\begin{prop}
\label{extg-td}
Soient $\E $, $\FF$ deux $\smash{\widetilde{\D}} _{X^\#} ^{(m)}$-modules à gauche,
$\E ' $ un $\B _{X } ' \otimes _{\O _{X }} \D _{X ^\#} ^{(m)}$-module à gauche,
$\M $ un $\smash{\widetilde{\D}} _{X ^\#} ^{(m)}$-module à droite et $\M '  $
un $\B _{X } ' \otimes _{\O _{X }} \D _{X ^\#} ^{(m)}$-module à droite.

La structure canonique (voir \ref{MotimesEetc}) de
$\smash{\widetilde{\D}} _{X^\#} ^{(m)}$-module
à gauche sur
$\B _{X } ' \otimes _{\B _{X }} \E  $
se prolonge en une structure de $\B _{X } ' \otimes _{\B _{X }} \D _{X^\#} ^{(m)}$-module à gauche.
De plus, les isomorphismes canoniques
\begin{gather}
\label{extg}
(\B _{X } ' \otimes _{\B _{X }} \E  ) \otimes _{\B _{X } '} \E '  \riso \E  \otimes _{\B _{X }} \E ' , \\
\label{hombb'}
\mathcal{H} om _{ \B _{X}} (\E, \E')
\riso
\mathcal{H} om _{ \B _{X}'} (\B _{X } ' \otimes _{\B _{X }} \E, \E'),
\\
\label{botimeshom}
\B _{X } ' \otimes _{\B _{X }} \mathcal{H} om _{ \B _{X}} (\E, \FF)
\riso \mathcal{H} om _{ \B _{X} '} (\B _{X } ' \otimes _{\B _{X }} \E, \B _{X } ' \otimes _{\B _{X }}  \FF)
\end{gather}
sont $\widetilde{\D} _{X ^\#} ^{(m)}$-linéaires
(et donc $\B _{X } '\otimes _{\O _{X }} \D _{X ^\#} ^{(m)}$-linéaire pour le dernier).
Par transport de structure, on munit ainsi $\E  \otimes _{\B _{X }} \E ' $ et $\mathcal{H} om _{ \B _{X}} (\E, \E') $
d'une structure canonique de $\B _{X } ' \otimes _{\O _{X }} \D _{X ^\#}^{(m)}$-module à gauche.

De même, on dispose d'une structure canonique $\B _{X } ' \otimes _{\O _{X }} \D _{X ^\#}^{(m)}$-module à gauche (resp. à droite)
sur $\mathcal{H} om _{ \B _{X}} (\M, \M')$ (resp. $\M  \otimes _{\B _{X }} \E ' $, $\M  ' \otimes _{\B _{X }} \E  $,
$\mathcal{H} om _{ \B _{X}} (\E, \M')$) ainsi que des isomorphismes canoniques analogues $\widetilde{\D} _{X ^\#} ^{(m)}$-linéaires.
\end{prop}
\begin{proof}
Concernant le prolongement en une structure de $\B _{X } ' \otimes _{\B _{X }} \D _{X^\#} ^{(m)}$-module, on vérifie que les
$m$-PD-(co)stratifications relatives à $\B _X$ sont
semi-linéaires par rapport à $(\widetilde{\epsilon} _n ^{\B '_X})^{-1}$ (voir \ref{eqDtilde} et \ref{Dmodstrat}).
La $\widetilde{\D} _{X ^\#} ^{(m)}$-linéarité découle tautologiquement des constructions respectives
des $m$-PD-(co)stratifications relatives à $\B _X$.

\end{proof}

\begin{lemm}\label{evalDlin}
Soient $\M$, $\NN$ deux $\widetilde{\D} ^{(m)} _{X ^\#}$-modules à droite. Le morphisme d'évaluation
$$ \mathrm{ev}\ : \ \M \otimes _{\B _X} \mathcal{H} om _{ \B _{X}} (\M, \NN) \rightarrow \NN$$
est $\widetilde{\D} ^{(m)} _{X ^\#}$-linéaire.
\end{lemm}
\begin{proof}
Analogue à celle de \ref{comm-assoc}.
\end{proof}

\begin{lemm}
\label{Hom-otimes}
  Soient $\E$, $\FF$, $\G$ trois $\widetilde{\D} ^{(m)} _{X ^\#}$-modules à gauche et
$\M$, $\NN$ deux $\widetilde{\D} ^{(m)} _{X ^\#}$-modules à droite. Les morphismes canoniques :
\begin{gather}
  \E \otimes _{\B _X} \mathcal{H} om _{\B _X} (\FF, \G) \rightarrow
  \mathcal{H} om _{\B _X} (\FF, \E \otimes _{\B _X} \G),\
  \E \otimes _{\B _X} \mathcal{H} om _{\B _X} (\M, \NN) \rightarrow
  \mathcal{H} om _{\B _X} (\M, \NN \otimes _{\B _X} \E ),\\
   \E \otimes _{\B _X} \mathcal{H} om _{\B _X} (\FF, \M) \rightarrow
  \mathcal{H} om _{\B _X} (\FF, \M \otimes _{\B _X} \E),\
   \M \otimes _{\B _X} \mathcal{H} om _{\B _X} (\FF, \E) \rightarrow
  \mathcal{H} om _{\B _X} (\FF, \M \otimes _{\B _X} \E)
\end{gather}
sont $\widetilde{\D} ^{(m)} _{X ^\#}$-linéaires.
\end{lemm}

\begin{proof}
  Cela s'établit en vérifiant leur commutation aux $m$-PD-(co)stratifications.
\end{proof}

\begin{lemm}
  \label{HomD}
Soient $\M$ un $\widetilde{\D} ^{(m)} _{X ^\#}$-module à droite,
$\NN$ un $\widetilde{\D} ^{(m)} _{X ^\#}$-bimodule,
$\E$ un $\widetilde{\D} ^{(m)} _{X ^\#}$-module à gauche.
On dispose d'un isomorphisme canonique de $\widetilde{\D} ^{(m)} _{X ^\#}$-modules à droite :
\begin{gather}\label{HomD1}
  (\M \otimes _{\B _X} \NN ) \otimes _{\widetilde{\D} ^{(m)} _{X ^\#}} \E
\riso
\M \otimes _{\B _X} (\NN  \otimes _{\widetilde{\D} ^{(m)} _{X ^\#}} \E),
\\ \label{HomD2}
(\M \otimes _{\widetilde{\D} ^{(m)} _{X ^\#}} \NN ) \otimes _{\B _X} \E
\riso
\M \otimes _{\widetilde{\D} ^{(m)} _{X ^\#}} (\NN  \otimes _{\B _X} \E).
\end{gather}
L'isomorphisme \ref{HomD1} reste valable lorsque
$\M$ est un $\widetilde{\D} ^{(m)} _{X ^\#}$-module à gauche.
\end{lemm}

\begin{proof}
On se contente par analogie de vérifier \ref{HomD1}.
  Dans un premier temps, établissons que l'isomorphisme canonique
$$\theta \ :\ (\M \otimes _{\B _X} \widetilde{\D} ^{(m)} _{X ^\#}) \otimes _{\widetilde{\D} ^{(m)} _{X ^\#}} \E
\riso
\M \otimes _{\B _X} \E$$
  est un isomorphisme de $\widetilde{\D} ^{(m)} _{X ^\#}$-modules à droite.
Soient $x$, $y$, $P$  des sections locales de respectivement $\M$, $\E$, $\widetilde{\D} ^{(m)} _{X ^\#}$. Il suffit de prouver la relation
 $\theta (((x \otimes 1) \otimes y)\cdot P) = (x \otimes y)\cdot P$. L'homomorphisme
de $\widetilde{\D} ^{(m)} _{X ^\#}$-modules à gauche $\widetilde{\D} ^{(m)} _{X ^\#} \rightarrow \E$ envoyant $1$ sur $y$ induit l'homomorphisme de
$\widetilde{\D} ^{(m)} _{X ^\#}$-modules à droite $\phi$ : $\M \otimes _{\B _X} \widetilde{\D} ^{(m)} _{X ^\#}\rightarrow \M \otimes _{\B _X} \E$.
Or, on vérifie par un calcul que l'image par $\phi$ de l'action de $P$ pour la structure gauche de $\widetilde{\D} ^{(m)} _{X ^\#}$-modules à droite
sur $(x \otimes 1)$ est égale à $\theta (((x \otimes 1) \otimes y)\cdot P)$. D'un autre côté, on obtient par
$\widetilde{\D} ^{(m)} _{X ^\#}$-linéarité de $\phi$ que cette image est $(x \otimes y)\cdot P$.

Il s'ensuit le cas général via les isomorphismes :
$$(\M \otimes _{\B _X} \NN ) \otimes _{\widetilde{\D} ^{(m)} _{X ^\#}} \E
\liso
(\M \otimes _{\B _X} \widetilde{\D} ^{(m)} _{X ^\#}) \otimes _{\widetilde{\D} ^{(m)} _{X ^\#}} \NN  \otimes _{\widetilde{\D} ^{(m)} _{X ^\#}} \E
\riso
\M \otimes _{\B _X} (\NN  \otimes _{\widetilde{\D} ^{(m)} _{X ^\#}} \E).$$
\end{proof}

\begin{nota}
\label{widetildeomega}
Notons $\widetilde{\omega} _{X ^\#}:= \B _X \otimes _{\O _X}\omega _{X ^\#}$,
$\B _Y:= \B _X  |Y$,
$\widetilde{\D} ^{(m)} _{Y}:= \B _Y \otimes _{\O _Y} \D ^{(m)} _{Y}$
et
$\widetilde{\omega}  _{Y}:= \B _Y \otimes _{\O _Y} \omega _Y$.
Le faisceau
$\widetilde{\omega}_{X ^\#}$ est muni d'une structure canonique de
$\widetilde{\D} ^{(m)} _{X ^\#}$-module à droite.
En effet, Montagnon l'a déjà établi pour
$\omega_{X ^\#}$ de la manière suivante :
on vérifie via les formules \cite[2.3.3.1]{Be1} et \cite[1.2.3]{Be2}
que $\omega _{X ^\#} $ est un sous-$\D ^{(m)} _{X ^\#}$-module à droite de
$j _* \omega _{Y}$.
On en déduit alors une structure de $\widetilde{\D} ^{(m)} _{X ^\#}$-module à droite sur
$\widetilde{\omega}_{X ^\#}$ (voir \ref{extg-td}).

Soit $\delta$ :
  $\omega_Y \otimes _{\O_Y} \D ^{(m)} _Y \riso
  \omega_Y \otimes _{\O_Y} \D ^{(m)} _Y $ l'isomorphisme de transposition (voir \cite[1.3.3]{Be2}).
Par un calcul en coordonnées locales (voir la formule \cite[1.3.1.1]{Be2}), on obtient de plus les factorisations :
\begin{equation}
  \label{deltadiese}
  \xymatrix {
  {\omega_{X } \otimes _{\O_X} \D ^{(m)} _{X ^\#}}
  \ar@{.>}[d] ^-\sim _-{\delta}
  \ar@{^{(}->}[r]
  &
  {\omega_{X ^\#} \otimes _{\O_X} \D ^{(m)} _{X ^\#}}
  \ar@{.>}[d] ^-\sim _-{\delta _\#}
  \ar@{^{(}->}[r]
  &
{j _* (\omega_Y \otimes _{\O_Y} \D ^{(m)} _Y)}
\ar[d] ^-\sim _-{\delta}
\\
{\omega_{X } \otimes _{\O_X} \D ^{(m)} _{X ^\#}}
\ar@{^{(}->}[r]
  &
{\omega_{X ^\#} \otimes _{\O_X} \D ^{(m)} _{X ^\#}}
\ar @{^{(}->}[r]
  &
{j _* ( \omega_Y \otimes _{\O_Y} \D ^{(m)} _Y).}
}
\end{equation}
Il découle de \cite[1.3.3]{Be2} que
ces deux factorisations sont aussi les uniques involutions échangeant les deux structures
de $ \D ^{(m)} _{X ^\#}$-modules à droite respectives et vérifiant,
pour toute section $x$ de $\omega_{\X }$ ou $\omega_{\X ^\#}$,
$x \otimes 1 \mapsto x \otimes1$ (cela implique en particulier que ce sont bien des isomorphismes).

Si $t _1,\dots, t _d$ sont des coordonnées locales logarithmiques de $X ^\#$,
en identifiant $\omega_{X }$ et $\O _X$ (à $a \in \O _X$ correspond
$a d t _1 \wedge \cdots \wedge d t _d$),
l'isomorphisme $\delta $ induit une flèche :
$\D ^{(m)} _{X ^\#} \rightarrow \D ^{(m)} _{X ^\#}$,
envoyant un opérateur différentiel $P $ de $\D ^{(m)} _{X ^\#}$ sur son {\it adjoint}
$  \overset{^\mathrm{t}}{} P  $ (voir \cite[1.2.2]{Be2}).
Pour tous $P , Q $ de $\D ^{(m)} _{X ^\#}$,
$ \overset{^\mathrm{t}}{} (P  \cdot Q  ) =  \overset{^\mathrm{t}}{} Q  \cdot \overset{^\mathrm{t}}{} P  $,
$ \overset{^\mathrm{t}}{} (\overset{^\mathrm{t}}{} P  ) =P$
et la structure gauche (resp. droite) de $\D ^{(m)} _{X ^\#}$-module à droite (via l'identification
de $\omega_{X }$ et $\O _X$) sur $\omega_{X } \otimes _{\O_X} \D ^{(m)} _{X ^\#}$
est donnée par la multiplication à gauche par l'adjoint (resp. par la multiplication à droite).

De même, en identifiant $\omega_{X ^\#}$ et $\O _X$ (via
$a \frac{d t _1 \wedge \cdots \wedge dt _d}{t _1 \cdots t_d}\leftrightarrow a \in \O _X$),
l'isomorphisme $\delta _\#$ induit une flèche :
$\D ^{(m)} _{X ^\#} \rightarrow \D ^{(m)} _{X ^\#}$, $ P  \mapsto \smash{\widetilde{P }} $.
Pour le différencier de l'adjoint,
$\smash{\widetilde{P }} $
sera appelé {\it adjoint logarithmique} de $P $.
On calcule que
$\smash{\widetilde{P }}  = t_1 \cdots t_d \cdot \overset{^\mathrm{t}}{} P  \cdot\frac{1}{t_1 \cdots t_d } $
(on remarque que ce dernier élément appartient bien
à $\D ^{(m)} _{X ^\#}$).

Pour tout $P   = \sum _{\underline{k}} b _{\underline{k}} \underline{\partial} ^{<\underline{k}>}_\#
\in \widetilde{\D} ^{(m)} _{X ^\#}$ avec $b _{\underline{k}}  \in  \B _X$,
on définit plus généralement les adjoints et adjoints logarithmiques de $P$ en posant
$\overset{^\mathrm{t}}{} P  : = \sum _{\underline{k}} \overset{^\mathrm{t}}{}
\underline{\partial} ^{<\underline{k}>}_\#  b _{\underline{k}}$ et
$\smash{\widetilde{P }} : = t_1 \cdots t_d \cdot \overset{^\mathrm{t}}{} P  \cdot\frac{1}{t_1 \cdots t_d } $.
\end{nota}

\begin{vide}
\label{dr-ga}
Soient $\E$ (resp. $\M$) un $\widetilde{\D} ^{(m)} _{X ^\#}$-module à gauche (resp. à droite).
Si $t _1,\dots, t _d$ sont des coordonnées locales logarithmiques de $X^\#$,
dans l'expression $\widetilde{\omega}_{X ^\#} \otimes _{\B _X} \E$, en identifiant
$\widetilde{\omega}_{X ^\#}$ et $\B _X$ (à $b \in \B _X$ correspond
$b \frac{d t _1 \wedge \cdots \wedge dt _d}{t _1 \cdots t_d}$), l'action à droite
de $P  \in \widetilde{\D} ^{(m)} _{X ^\#}$ sur $e \in \E$ est
$\smash{\widetilde{P }}  .e$.
En effet, par $\B_X$-linéarité, il suffit de le vérifier pour
$P  \in \D ^{(m)} _{X ^\#}$.
Or, via \ref{extg}, $\widetilde{\omega}_{X ^\#} \otimes _{\B _X} \E \riso \omega_{X ^\#} \otimes _{\O _X} \E$.
On utilise ensuite l'isomorphisme canonique
$(\omega_{X ^\#} \otimes _{\O _X} \D ^{(m)} _{X ^\#}) \otimes _{\D ^{(m)} _{X ^\#}} \E
\riso \omega_{X ^\#} \otimes _{\O _X} \E$ de \ref{HomD1}.
On obtient la description analogue pour
$\widetilde{\omega}_{X } \otimes _{\B _{X}} \E$
en échangeant $\smash{\widetilde{P }}$ par $\overset{^\mathrm{t}}{}P  $.

En utilisant la formule \cite[1.1.7.3]{Be2}, on
calcule que l'isomorphisme canonique
$\mathcal{H} om _{\O _X } (\omega_{X}, \omega_{X })
\riso \O _X$
est $\D ^{(m)}  _{X}$-linéaire (et donc $\D ^{(m)}  _{X ^\#}$-linéaire).
Pour obtenir l'isomorphisme analogue avec des dièses, on utilise
les inclusions
$\mathcal{H} om _{\O _X} (\omega_{X ^\# }, \omega_{X ^\# })
\hookrightarrow j_* \mathcal{H} om _{\O _Y } (\omega_{Y}, \omega_{Y })$
et
$\O _X \hookrightarrow j_* \O _Y$.
Puis, on calcule que cela induit
la
$\D ^{(m)}  _{X ^\#}$-linéarité de
$\mathcal{H} om _{\O _X} (\omega_{X ^\# }, \omega_{X ^\# })
\riso
\O _X$.
Par \ref{botimeshom}, il en découle que l'isomorphisme canonique
$\mathcal{H} om _{\B _X} (\widetilde{\omega}_{X ^\# }, \widetilde{\omega}_{X ^\# })
\riso
\B _X$
est $\widetilde{\D} ^{(m)}  _{X ^\#}$-linéaire.

Par \ref{evalDlin},
l'isomorphisme canonique
$\widetilde{\omega}_{X ^\# }\otimes _{\B _X} \mathcal{H} om _{\B _X} (\widetilde{\omega}_{X ^\# }, \M)
\riso \M$
est $\widetilde{\D} ^{(m)}  _{X ^\#}$-linéaire.
De plus, via \ref{Hom-otimes}, on dispose de l'isomorphisme canonique
$\widetilde{\D} ^{(m)}  _{X ^\#}$-linéaire :
$\mathcal{H} om _{\B _X} (\widetilde{\omega}_{X ^\# }, \widetilde{\omega}_{X ^\# }\otimes _{\B _X}  \E )
\riso \mathcal{H} om _{\B _X} (\widetilde{\omega}_{X ^\# }, \widetilde{\omega}_{X ^\# }) \otimes _{\B _X}  \E
\riso \E$. On bénéficie des isomorphismes similaires en remplaçant $\widetilde{\omega}_{X ^\# }$ par
$\widetilde{\omega}_{X }$.
Il en dérive que
les foncteurs $-\otimes _{\B _{X}} \widetilde{\omega}_{X ^\# } ^{-1}$
et $\widetilde{\omega}_{X ^\# }  \otimes _{\B _{X}} -$
(resp. $-\otimes _{\B _{X}} \widetilde{\omega}_{X  } ^{-1}$
et $\widetilde{\omega}_{X }  \otimes _{\B _{X}} -$)
induisent des équivalences quasi-inverses entre
la catégorie des $\widetilde{\D} ^{(m)}  _{X ^\#}$-modules à gauche et celle
des $\widetilde{\D} ^{(m)}  _{X ^\#}$-modules à droite.

On en déduit que, dans $\M \otimes _{\B _{X}} \widetilde{\omega}_{X ^\# } ^{-1}$,
en identifiant comme précédemment $\widetilde{\omega}_{X ^\#}$ et $\B _X$
l'action à gauche
de $P  \in \widetilde{\D} ^{(m)} _{X ^\#}$ sur $m \in \M$ est égale à
$m\cdot \smash{\widetilde{P }}$.
On obtient une description analogue pour
$\M \otimes _{\B _{X}} \widetilde{\omega}_{X } ^{-1}$
en remplaçant $\smash{\widetilde{P }}$ par $\overset{^\mathrm{t}}{}P  $.
Ces structures sont appelés structures {\og tordues\fg}.
\end{vide}

\begin{lemm}
\label{dr-gabis}
Les foncteurs quasi-inverses $-\otimes _{\B _{X}} \widetilde{\omega}_{X ^\# } ^{-1}$
et $\widetilde{\omega}_{X ^\# }  \otimes _{\B _{X}} -$
(resp. $-\otimes _{\B _{X}} \widetilde{\omega}_{X  } ^{-1}$
et $\widetilde{\omega}_{X }  \otimes _{\B _{X}} -$) de \ref{dr-ga}
sont exacts et induisent des équivalences entre
la catégorie des $\widetilde{\D} ^{(m)}  _{X ^\#}$-modules
(resp. cohérents, resp. plats, resp. localement projectifs de type fini) à gauche et celle
des $\widetilde{\D} ^{(m)}  _{X ^\#}$-modules
(resp. cohérents, resp. plats, resp. localement projectifs de type fini) à droite.
\end{lemm}

\begin{proof}
  L'exactitude est triviale. En ce qui concerne la cohérence et la projectivité locale de type fini, il suffit
d'utiliser les descriptions de \ref{dr-ga} sur les structures tordues
(et de se rappeler que les flèches qui associent à un opérateur différentiel son adjoint
ou adjoint logarithmique sont des isomorphismes).
Enfin, pour la platitude, cela découle, pour tous
$\widetilde{\D}  _{X ^\#}$-module à gauche $\E$ et
$\widetilde{\D}  _{X ^\#}$-module à droite $\M$,
de l'isomorphisme canonique (voir \cite[I.2.2]{virrion})
$$ \M \otimes _{\widetilde{\D} _{X ^\# } }  \E
\riso
(\omega _{X ^\#} \otimes _{\O _{X}} \E )
\otimes _{\widetilde{\D} _{X ^\# } }
(\M\otimes _{\O _{X}} \omega _{X ^\#} ^{-1} ),$$
de même en remplaçant {\og $\omega _{X ^\#}$\fg} par {\og $\omega _{X }$\fg}.

\end{proof}

\begin{prop}
\label{propMotimesEetcegpre}
Soient $\E$, $\FF$ deux $\widetilde{\D} ^{(m)} _{X ^\#}$-modules à gauche et
$\M$ un $\widetilde{\D} ^{(m)} _{X ^\#}$-module à droite.
Si $t _1,\dots, t _d$ sont des coordonnées locales logarithmiques de $X^\#$, 
pour tout $\underline{k} \in \N ^d$, pour toutes sections
$e \in \E $, $f \in \FF$, $m \in \M$, la structure canonique de $\widetilde{\D} ^{(m)} _{X ^\#}$-modules à gauche
(resp. à droite) sur $\E \otimes _{\B _X} \FF$ (resp. $\M \otimes _{\B _X} \E$)
est caractérisée par la formule \ref{MotimesEetcegpre} (resp. \ref{MotimesEetceg} ou \ref{MotimesEetcegbis}) ci-dessous :
\begin{gather}
  \label{MotimesEetcegpre}
  \underline{\partial} ^{<\underline{k}>}_\# (e\otimes f)
=
\sum _{\underline{h} \leq \underline{k}}
\left \{ \begin{smallmatrix}   \underline{k} \\   \underline{h} \\ \end{smallmatrix} \right \}
\underline{\partial} ^{<\underline{k}-\underline{h}>}_\# e
\otimes
\underline{\partial} ^{<\underline{h}>}_\# f,
\\
  \label{MotimesEetceg}
  (m \otimes e) \smash{\widetilde{\underline{\partial}}} ^{<\underline{k}>}_\#
=
\sum _{\underline{h} \leq \underline{k}}
\left \{ \begin{smallmatrix}   \underline{k} \\   \underline{h} \\ \end{smallmatrix} \right \}
m \smash{\widetilde{\underline{\partial}}} ^{<\underline{k}-\underline{h}>}_\#
\otimes
\underline{\partial} ^{<\underline{h}>}_\# e,
\\
  \label{MotimesEetcegbis}
  (m \otimes e) \overset{^\mathrm{t}}{} \underline{\partial} ^{<\underline{k}>}_\#
=
\sum _{\underline{h} \leq \underline{k}}
\left \{ \begin{smallmatrix}   \underline{k} \\   \underline{h} \\ \end{smallmatrix} \right \}
m \overset{^\mathrm{t}}{} \underline{\partial} ^{<\underline{k}-\underline{h}>}_\#
\otimes
\underline{\partial} ^{<\underline{h}>}_\# e.
\end{gather}
\end{prop}

\begin{proof}
  La preuve de \ref{MotimesEetcegpre} est analogue à \cite[2.3.3.1]{Be1}.
  On en déduit \ref{MotimesEetceg} et \ref{MotimesEetcegbis} par passage de gauche à droite (i.e., via les équivalences de catégories
  de \ref{widetildeomega}) et via \ref{comm-assoc2}.
\end{proof}

Nous aurons besoin de la proposition suivante pour obtenir \ref{md(eof)=} qui nous permettra de prouver
\ref{theo-compDR}.
\begin{prop}
\label{isotransp}
 Soient $\E$ un $\widetilde{\D} ^{(m)} _{X ^\#}$-module à gauche, $\E \otimes _{\B_X} \widetilde{\D} ^{(m)} _{X ^\#}$ et
$ \widetilde{\D} ^{(m)} _{X ^\#} \otimes _{\B_X} \E$ les faisceaux obtenus en calculant le produit tensoriel via
la structure gauche et respectivement droite de $\B _X$-algèbre de $ \widetilde{\D} ^{(m)} _{X ^\#} $. Il existe un unique isomorphisme de
$ \widetilde{\D} ^{(m)} _{X ^\#} $-bimodules :
$$\gamma _{\E}\ :\ \widetilde{\D} ^{(m)} _{X ^\#} \otimes _{\B_X} \E \riso
\E \otimes _{\B_X} \widetilde{\D} ^{(m)} _{X ^\#}$$
tel que, pour toute section $e$ de $\E$,
$\gamma _{\E} ( 1\otimes e) = e \otimes 1$.
\end{prop}

\begin{proof}
  La preuve est analogue à celle de \cite[1.3.1 et 1.3.2]{Be2} : dans un premier temps, on suppose
  $\B _X = \O _X$. L'homomorphisme $\gamma _{\E}$ est défini de manière unique via
  la formule $\gamma _{\E} ( P \otimes e) := P (e \otimes 1)$, où $P \in \D ^{(m)} _{X ^\#}$ et $e \in \E$.
  Vérifions à présent la $\D ^{(m)} _{X ^\#}$-linéarité à droite.
  Il suffit d'établir, pour tout $\underline{k} \in \N ^d$,
  $\gamma _{\E} (( 1 \otimes e)\overset{^\mathrm{t}}{} \underline{\partial} ^{<\underline{k}>}_\#)=
  ( e \otimes 1)\overset{^\mathrm{t}}{} \underline{\partial} ^{<\underline{k}>}_\#
  (=e \otimes \overset{^\mathrm{t}}{} \underline{\partial} ^{<\underline{k}>}_\#)$.

  Par \ref{MotimesEetcegbis},
$$
\gamma _{\E} (( 1 \otimes e)\overset{^\mathrm{t}}{} \underline{\partial} ^{<\underline{k}>}_\#)=
\gamma _{\E} (\sum _{\underline{h} \leq \underline{k}}
\left \{ \begin{smallmatrix}   \underline{k} \\   \underline{h} \\ \end{smallmatrix} \right \}
\overset{^\mathrm{t}}{} \underline{\partial} ^{<\underline{h}>}_\#
\otimes
\underline{\partial} ^{<\underline{k}-\underline{h}>}_\# e)
=
\sum _{\underline{h} \leq \underline{k}}
\left \{ \begin{smallmatrix}   \underline{k} \\   \underline{h} \\ \end{smallmatrix} \right \}
\overset{^\mathrm{t}}{} \underline{\partial} ^{<\underline{h}>}_\#
(\underline{\partial} ^{<\underline{k}-\underline{h}>}_\# e \otimes 1).
$$
Pour tout $\underline{h} \in \N ^d$, on note $q _i ^{(h _i)}$ le quotient de la division euclidienne de
$h_i$ par $p^m$ et $\underline{q} ^{(\underline{h})}!:=q _1 ^{(h _1)}!\cdots q _d ^{(h _d)}!$.
Grâce à \cite[2.2.4.(ii) et (iv)]{Be1}, on calcule :
$(-1) ^{|\underline{h}|}\overset{^\mathrm{t}}{} \underline{\partial} ^{<\underline{h}>}_\# =
\underline{\partial} ^{<\underline{h}>}  \underline{t} ^{\underline{h}}  =
\sum _{\underline{i} \leq \underline{h}} \underline{q} ^{(\underline{h}-\underline{i})}!
\left \{ \begin{smallmatrix}   \underline{h} \\   \underline{i} \\ \end{smallmatrix} \right \}
\left ( \begin{smallmatrix}   \underline{h} \\   \underline{i} \\ \end{smallmatrix} \right )
\underline{t} ^{\underline{i}} \underline{\partial} ^{<\underline{i}>}
=
\sum _{\underline{i} \leq \underline{h}} \underline{q} ^{(\underline{h}-\underline{i})}!
\left \{ \begin{smallmatrix}   \underline{h} \\   \underline{i} \\ \end{smallmatrix} \right \}
\left ( \begin{smallmatrix}   \underline{h} \\   \underline{i} \\ \end{smallmatrix} \right )
\underline{\partial} ^{<\underline{i}>} _\#  $.
D'où :
\begin{align}
  \notag
\gamma _{\E} (( 1 \otimes e)\overset{^\mathrm{t}}{} \underline{\partial} ^{<\underline{k}>}_\#)
&=
\sum _{\underline{h} \leq \underline{k}}
\sum _{\underline{i} \leq \underline{h}}
(-1) ^{|\underline{h}|}
\underline{q} ^{(\underline{h}-\underline{i})}!
\left \{ \begin{smallmatrix}   \underline{k} \\   \underline{h} \\ \end{smallmatrix} \right \}
\left \{ \begin{smallmatrix}   \underline{h} \\   \underline{i} \\ \end{smallmatrix} \right \}
\left ( \begin{smallmatrix}   \underline{h} \\   \underline{i} \\ \end{smallmatrix} \right )
\underline{\partial} ^{<\underline{i}>} _\#
(\underline{\partial} ^{<\underline{k}-\underline{h}>}_\# e \otimes 1)
        \\  \notag
 &       =
\sum _{\underline{h} \leq \underline{k}}
\sum _{\underline{i} \leq \underline{h}}
\sum _{\underline{j} \leq \underline{i}}
(-1) ^{|\underline{h}|}
\underline{q} ^{(\underline{h}-\underline{i})}!
\left \{ \begin{smallmatrix}   \underline{k} \\   \underline{h} \\ \end{smallmatrix} \right \}
\left \{ \begin{smallmatrix}   \underline{h} \\   \underline{i} \\ \end{smallmatrix} \right \}
\left ( \begin{smallmatrix}   \underline{h} \\   \underline{i} \\ \end{smallmatrix} \right )
\left \{ \begin{smallmatrix}   \underline{i} \\   \underline{j} \\ \end{smallmatrix} \right \}
\underline{\partial} ^{<\underline{i}-\underline{j}>} _\# \underline{\partial} ^{<\underline{k}-\underline{h}>}_\# e
\otimes
\underline{\partial} ^{<\underline{j}>} _\#)
        \\  \notag
  &      =
\sum _{\underline{j} \leq \underline{k}}
\left (
\sum _{\underline{j} \leq \underline{i}\leq \underline{h} \leq \underline{k}}
(-1) ^{|\underline{h}|}
\underline{q} ^{(\underline{h}-\underline{i})}!
\left \{ \begin{smallmatrix}   \underline{k} \\   \underline{h} \\ \end{smallmatrix} \right \}
\left \{ \begin{smallmatrix}   \underline{h} \\   \underline{i} \\ \end{smallmatrix} \right \}
\left ( \begin{smallmatrix}   \underline{h} \\   \underline{i} \\ \end{smallmatrix} \right )
\left \{ \begin{smallmatrix}   \underline{i} \\   \underline{j} \\ \end{smallmatrix} \right \}
\underline{\partial} ^{<\underline{i}-\underline{j}>} _\# \underline{\partial} ^{<\underline{k}-\underline{h}>}_\#
\right ) e
\otimes
\underline{\partial} ^{<\underline{j}>} _\#).
\end{align}
Pour conclure, il suffit alors d'établir la formule
\begin{equation}
  \label{cekilreste}
\sum _{\underline{j} \leq \underline{i}\leq \underline{h} \leq \underline{k}}
(-1) ^{|\underline{h}|}
\underline{q} ^{(\underline{h}-\underline{i})}!
\left \{ \begin{smallmatrix}   \underline{k} \\   \underline{h} \\ \end{smallmatrix} \right \}
\left \{ \begin{smallmatrix}   \underline{h} \\   \underline{i} \\ \end{smallmatrix} \right \}
\left ( \begin{smallmatrix}   \underline{h} \\   \underline{i} \\ \end{smallmatrix} \right )
\left \{ \begin{smallmatrix}   \underline{i} \\   \underline{j} \\ \end{smallmatrix} \right \}
\underline{\partial} ^{<\underline{i}-\underline{j}>} _\# \underline{\partial} ^{<\underline{k}-\underline{h}>}_\#
=
(-1) ^{|\underline{k}|}
\underline{q} ^{(\underline{k}-\underline{j})}!
\left \{ \begin{smallmatrix}   \underline{k} \\   \underline{j} \\ \end{smallmatrix} \right \}
\left ( \begin{smallmatrix}   \underline{k} \\   \underline{j} \\ \end{smallmatrix} \right ) .
\end{equation}
À cette fin, procédons comme suit.
Lorsque $\E$ est égal à $\D ^{(m)} _{X ^\#}$ (pour éviter les confusions, on le note toujours $\E$),
la vérification de la proposition est plus aisée car il suffit de faire un calcul local (via la formule \cite[1.3.1.1]{Be2})
pour constater que le morphisme de \cite[1.3.1]{Be2} se factorise de la manière suivante :
\begin{equation}
  \label{diaggammadiese}
  \xymatrix {
{\D ^{(m)} _{X ^\#} \otimes _{\O _X} \E }
  \ar@{.>}[d] _-{\gamma }
  \ar@{^{(}->}[r]
  &
{j _* (\D ^{(m)} _Y \otimes _{\O _Y} \E |Y)}
\ar[d] ^-\sim _-{\gamma }
\\
  {\E \otimes _{\O _X} \D ^{(m)} _{X ^\#}}
\ar @{^{(}->}[r]
  &
{j _* (\E |Y \otimes _{\O _Y} \D ^{(m)} _Y).}
}
\end{equation}
Cette factorisation envoie bien $1\otimes e$ sur $e \otimes 1$ et correspond donc (par unicité) au
morphisme $\gamma _\E$ que l'on a défini en début de preuve. En reprenant les calculs déjà faits,
comme $\D ^{(m)} _{X ^\#} $ est un $\O _X$-module libre, on obtient alors la formule
\ref{cekilreste} recherchée.

Il reste à vérifier que $\gamma _\E$ est un isomorphisme. Il suffit pour cela de construire
de manière analogue l'unique morphisme de $\D ^{(m)} _{X ^\#} $-bimodules
$\E \otimes _{\O_X} \D ^{(m)} _{X ^\#} \rightarrow
\D ^{(m)} _{X ^\#} \otimes _{\O_X} \E $
qui envoie, pour toute section $e$ de $\E$,
$e \otimes 1 $ sur $ 1\otimes e $.

Traitons à présent le cas où $\B _X$ est quelconque. On construit l'isomorphisme $\gamma _\E$ de la façon suivante :
\begin{equation}
  \notag
\widetilde{\D} ^{(m)} _{X ^\#} \otimes _{\B_X} \E
\underset{\gamma_{\B _X}}{\liso}
(\D ^{(m)} _{X ^\#} \otimes _{\O _X} \B_X)\otimes _{\B_X} \E
\riso
\D ^{(m)} _{X ^\#} \otimes _{\O_X} \E
\underset{\gamma}{\riso}
\E \otimes _{\O_X} \D ^{(m)} _{X ^\#}
\riso
\E \otimes _{\B_X} \widetilde{\D} ^{(m)} _{X ^\#}.
\end{equation}
\end{proof}

\begin{prop}
\label{isotranspd}
  Soient $\M$ un $\widetilde{\D} ^{(m)} _{X ^\#}$-module à droite, $\M \otimes _{\B _X} \widetilde{\D} ^{(m)} _{X ^\#}$ le faisceau
  obtenu en calculant le produit tensoriel via la structure gauche de $\B _X$-algèbre de $\widetilde{\D} ^{(m)} _{X ^\#}$. Il existe une unique
  involution de $\widetilde{\D} ^{(m)} _{X ^\#}$-bimodules à droite
  $\widetilde{\delta} _{\M}$ :
  $\M \otimes _{\B _X} \widetilde{\D} ^{(m)} _{X ^\#} \riso \M \otimes _{\B _X} \widetilde{\D} ^{(m)} _{X ^\#}$ échangeant
  les deux structures de $\widetilde{\D} ^{(m)} _{X ^\#}$-modules à droite et telle que, pour toute section $m$ de $\M$,
$  \widetilde{\delta} _{\M} ( m \otimes 1) =m \otimes 1$.
\end{prop}

\begin{proof}
  La preuve est analogue à celle de \ref{isotransp} où on remplace \cite[1.3.1]{Be2} par
  \cite[1.3.3]{Be2}.
\end{proof}

\begin{vide}
\label{tranpDotimesB}
  De manière analogue à \cite[2.3.5]{Be1}, on bénéficie des formules caractérisant
  la structure d'anneau de $\B _X \otimes _{\O _X} \D ^{(m)} _{X ^\#}$ :
  $(b \otimes 1)(1\otimes P) = b\otimes P$ pour tous $b\in \B _X$ et $P \in \D  ^{(m)} _{X ^\#}$ et
\begin{equation}
\label{prodtensdgann}
 (1\otimes \underline{\partial} ^{<\underline{k}>}_\#) (b \otimes 1) =\sum _{\underline{h} \leq \underline{k}}
\left \{ \begin{smallmatrix}   \underline{k} \\   \underline{h} \\ \end{smallmatrix} \right \}
\underline{\partial} ^{<\underline{k}-\underline{h}>}_\# b
\otimes
\underline{\partial} ^{<\underline{h}>}_\#.
\end{equation}

On munit le faisceau $\D ^{(m)} _{X ^\#}\otimes _{\O _X} \B _X $ d'une structure
d'anneau via l'isomorphisme de transposition
$\gamma _{\B _X}$ : $\D ^{(m)} _{X^\#} \otimes
_{\O _X} \B _X \riso \smash{\widetilde{\D}}  ^{(m)} _{X^\#}$.
Pour tous $b\in \B _X$ et $P \in \D  ^{(m)} _{X ^\#}$, on obtient les formules :
$(P \otimes 1)(1\otimes b) = P\otimes b$ et
\begin{gather}
\label{prodtensdgannbis}
(1 \otimes b)(\overset{^\mathrm{t}}{} \underline{\partial} _\# ^{<\underline{k}>}\otimes 1)
= \sum _{\underline{h} \leq
\underline{k}} \left\{ \begin{smallmatrix} \underline{k} \\ \underline{h} \end{smallmatrix} \right\}
 \overset{^\mathrm{t}}{} \underline{\partial} _\# ^{<\underline{k}-\underline{h}>}\otimes \underline{\partial} _\# ^{<\underline{h}>}b .
\end{gather}

Le morphisme canonique
$\B  _X \otimes _{\B _X}\D  _{X^\#} ^{(m)}
\rightarrow
\B ' _X \otimes _{\B _X}\D  _{X^\#} ^{(m)} $
est aussi un homomorphisme d'anneaux (cela se voit par exemple via la formule \ref{prodtensdgann}).
Si $\E ^{\prime}$ est un $\smash{\widetilde{\D}}  _{X^\#} ^{(m)}$-module,
on vérifie (par exemple via les formules \ref{MotimesEetcegpre} et \ref{prodtensdgann})
que l'homomorphisme canonique
\begin{equation}
\label{extbg} \B '_X \otimes _{\B _X} \E ^{\prime} \rightarrow
(\B ' _X \otimes _{\B _X}\D  _{X^\#} ^{(m)} )\otimes _{\smash{\widetilde{\D}}  _{X^\#} ^{(m)}} \E ^{\prime},
\end{equation}
est un isomorphisme de $\B ' _X \otimes _{\O _X}\D  _{X^\#} ^{(m)} $-modules à gauche.

Si $\M$ est un $\D _{X^\#} ^{(m)} \otimes _{\O _X} \B _X$-module à droite,
on établit (par exemple via les formules \ref{MotimesEetcegbis} et \ref{prodtensdgannbis})
que l'homomorphisme canonique
\begin{equation}
\label{extbd} \M \otimes _{\B _X} \B _X '\rightarrow \M \otimes _{\D _{X^\#} ^{(m)}  \otimes _{\O _X} \B _X } (\D
_{X^\#} ^{(m)} \otimes _{\O _X} \B _X '),
\end{equation}
est un isomorphisme de $\D _{X^\#} ^{(m)} \otimes _{\O _X} \B _X '$-modules à droite.

\end{vide}

\section{\label{chap-coh}Cohérence et résolutions de Spencer}
\begin{vide}
\label{chap-coh-21}
Soit $\B _\X$ une $\O _\X$-algèbre commutative munie d'une structure de $\D ^{(m)} _{\X^\#}$-module à gauche compatible
à sa structure de $\O _\X$-algèbre. Pour tout $i\geq 0$, on pose $\B _{X _i} = \B _\X / \m ^{i+1} \B _\X $.
De plus, sauf mention explicite du contraire,
on supposera qu'il existe une base d'ouverts affines $\mathfrak{B}$ de $X$ telle que
\begin{enumerate}
  \item [(a)] Pour tout $U \in \mathfrak{B} $, l'anneau $\Gamma (U,\B _\X)$ est noethérien ;
  \item [(b')]  Pour tous $U,V \in \mathfrak{B}$ tels que $V \subset U$, l'homomorphisme
  $\Gamma (U,\B _\X )  \rightarrow \Gamma (V,\B _\X ) $ est plat.
  \item [(b)]  Pour tout $i\geq 0$, $\B _{X _i}$ est un $\O _{X _i}$-module quasi-cohérent et l'homomorphisme
 $\B _\X \rightarrow \underset{\longleftarrow}{\mathrm{lim}} \B _{X _i}$ est un isomorphisme.
\end{enumerate}
On remarque que la condition $(b)$ implique $(b')$.
D'après \cite[3.1 et 3.3]{Be1}, on dispose alors de théorèmes de type $A$ et $B$ pour les faisceaux d'anneaux $\B _X$.
Avec ces hypothèses supplémentaires, nous conservons les notations du chapitre \ref{nota1}.

\end{vide}

\begin{exem}
  \label{exB(T)bis}
  Soient $T$ un diviseur de $X _0$ et $r$ un multiple de $p ^{m+1}$.
  On dispose par \ref{exB(T)} du faisceau de $\O _\X$-algèbre $\B _\X  (T,r)$
  muni d'une structure canonique de $\D ^{(m)}_{\X^\#}$-module à gauche
  compatible à sa structure de $\O _\X$-algèbre.
  Son complété $p$-adique $\widehat{\B} _\X  (T,r)$ est aussi muni d'une structure compatible de
$\D ^{(m)}_{\X^\#}$-module à gauche (on le voit par exemple en s'assurant que la formule de Leibnitz \ref{Leibnitz}
reste valable).
Le faisceau $\widehat{\B} _\X  (T,r)$ satisfait toutes les conditions de \ref{chap-coh-21} (voir \cite[4.3.2]{Be1}).
Enfin, on pose $\widehat{\B} _\X ^{(m)} (T )= \widehat{\B} _\X  (T ,p^{m+1})$.
\end{exem}

\begin{prop}
\label{lemm-wildetildeDcoh}
Soit $\CC _X$ une $\O _X$-algèbre commutative munie d'une structure de $\D ^{(m)} _{X^\#}$-module à gauche compatible
à sa structure de $\O _X$-algèbre.
\begin{enumerate}
\item \label{lemm-wildetildeDcohi} L'anneau gradué (associé à la filtration par l'ordre)
$\mathrm{gr} \widetilde{\D} ^{(m)} _{X ^\#}$ est un anneau commutatif. Si $X ^\#$ est muni de coordonnées logarithmiques,
la relation
$\underline{\partial}  _\# ^{<\underline{k}>} \underline{\partial}  ^{<\underline{h}>} _\#
=
\left <\begin{smallmatrix}  \underline{k} + \underline{h} \\  \underline{k} \\ \end{smallmatrix} \right >
\underline{\partial}  ^{<\underline{k}+\underline{h}>} _\#$ devient exacte dans
$\mathrm{gr} \widetilde{\D} ^{(m)} _{X ^\#}$.
  \item \label{lemm-wildetildeDcohii} Si $X ^\#$ est muni de coordonnées logarithmiques,
le faisceau $\CC _X \otimes _{\O _X}  \D ^{(m)} _{X^\#}$
est engendré comme $\O _\X$-algèbre par
les opérateurs $\partial_{\#i} ^{<p^j>_{(m)} }$, où $1\leq i\leq d$, $0\leq j\leq m$,
ces derniers commutant deux à deux.
\item Pour tout ouvert affine $U\subset X$, l'homomorphisme canonique
  \begin{equation}
    \label{lemm-wildetildeDcoh2361}
    \Gamma (U, \CC _X) \otimes _{\Gamma (U,\O _X) } \Gamma (U,\D ^{(m)} _{X^\#})
    \rightarrow \Gamma (U, \CC _X \otimes _{\O _X}  \D ^{(m)} _{X^\#})
  \end{equation}
est un isomorphisme.
\item Si $\CC _X$ satisfait la condition $(a)$ de \ref{chap-coh-21} alors,
  pour tout ouvert affine $U\in \mathfrak{B}$, l'anneau $\Gamma (X ^\#, \smash{\widetilde{\D}} ^{(m)} _{X ^\#})$
  (resp. $\smash{\widetilde{\D}} ^{(m)} _{X ^\# ,x}$ pour tout $x \in X ^\#$) est noethérien à droite et à gauche.

\item Si $\CC _X$ satisfait les conditions $(a)$ et $(b')$ de \ref{chap-coh-21} alors
  le faisceau d'anneaux $\smash{\widetilde{\D}} ^{(m)} _{X ^\#}$ est cohérent à droite et à gauche.
\end{enumerate}
\end{prop}
\begin{proof}
Lorsque $\CC _X =\O _X$, cela correspond à \cite[Propositions 2.3.1-2]{these_montagnon}.
Autrement, on procède de manière identique à \cite[2.2.5]{Be1}, \cite[2.3.6]{Be1} et \cite[3.1.2]{Be1}.
\end{proof}

\begin{defi}
\label{defbonnefilt}
Soit $\M$ un $\smash{\widetilde{\D}} ^{(m)} _{X ^\#}$-module.
Une {\it bonne filtration} de $\M$ est une famille croissante exhaustive
$(\M _r) _{r\in \N}$ de sous-$\B _{X }$-modules cohérents de $\M$ telle que :
\begin{enumerate}
 \item  Pour tous $r,\ s\in \N$,
 $\smash{\widetilde{\D}} ^{(m)} _{X ^\# , r} \cdot \M _s \subset \M _{r+s}$ ;
 \item \label{defbonnefiltii}Il existe un entier $r_1 \in \N$ tel que pour tout entier $r\geq r_1$, on ait :\newline
$\M _r =\sum_{j=0} ^{p^m -1} \smash{\widetilde{\D}} ^{(m)} _{X ^\#, r-r_1 +j}\cdot \M_{r _1 -j} .$
\end{enumerate}
Lorsque la condition \ref{defbonnefiltii}) n'est pas validée, la famille $(\M _r) _{r\in \N}$
définit seulement une {\it filtration}.
\end{defi}

\begin{exem}
On vérifie par un calcul analogue à \cite[2.2.4]{caro_devissge_surcoh}
que la famille $(\smash{\widetilde{\D}} ^{(m)} _{X ^\#, r}) _{r\in \N}$ est une bonne filtration de
$\smash{\widetilde{\D}} ^{(m)} _{X ^\#}$ vérifiant la condition $\ref{defbonnefilt}.\ref{defbonnefiltii}$ pour tout
entier $r_1 \geq p ^m -1$.

\end{exem}
\begin{prop}
\label{globpfbfil}
Un $\smash{\widetilde{\D}} ^{(m)} _{X ^\#}$-module globalement de présentation finie admet une bonne filtration.
\end{prop}
\begin{proof}
  Analogue à \cite[2.2.5]{caro_devissge_surcoh}.
 \end{proof}

\begin{theo}[Théorème B]\label{theoB}
  On suppose $X ^\#$ affine et soit $\M$ un $\smash{\widetilde{\D}} ^{(m)} _{X ^\#}$-module globalement de présentation finie.
  Alors, pour tout entier $i \neq 0$, $H ^i (X ^\# , \M) =0$.
\end{theo}
\begin{proof}
  Cela résulte de \ref{globpfbfil} et du fait que,
  pour tout entier $i \neq 0$, le foncteur $H ^i (X ^\# , -)$ commute aux limites inductives filtrantes et
  du théorème de type $B$ pour les $\B _{X }$-modules cohérents.
 \end{proof}

\begin{theo}[Théorème A]\label{theoA}
  Lorsque $X ^\#$ est affine, les foncteurs $\M \mapsto \Gamma (X ^\#, \M)$ et
  $M \mapsto \smash{\widetilde{\D}} ^{(m)} _{X ^\#} \otimes _{\Gamma (X ^\#,\smash{\widetilde{\D}} ^{(m)} _{X ^\#})} M$ établissent des équivalences
  quasi-inverses entre la catégorie des $\smash{\widetilde{\D}} ^{(m)} _{X ^\#} $-modules globalement de présentation finie et
  celle des $\Gamma (X ^\#,\smash{\widetilde{\D}} ^{(m)} _{X ^\#} )$-modules de type fini.
\end{theo}
\begin{proof}
  Il s'agit de calquer \cite[2.2.7]{caro_devissge_surcoh}.
 \end{proof}

\begin{theo}\label{cohlocbfilt}
Soit $\M$ un $\smash{\widetilde{\D}} ^{(m)} _{X ^\# } $-module.
Alors $\M$ est cohérent si et seulement s'il admet localement de bonnes filtrations.
\end{theo}
\begin{proof}
  Similaire à \cite[2.2.8]{caro_devissge_surcoh}.
 \end{proof}

\begin{rema}

\begin{enumerate}
  \item Dans le cas où $X ^\#$ est un log-schéma, les assertions \ref{globpfbfil}, \ref{theoB}, \ref{theoA}
restent valables en remplaçant l'hypothèse
{\og globalement de présentation finie\fg}
par
{\og cohérent\fg}.
\item Ce chapitre reste valable en remplaçant les modules (à gauche) par des modules à droite.
\end{enumerate}
\end{rema}

  Nous nous restreignons dans la suite de cette section au niveau $0$ car
  les énoncés analogues ne sont plus valables pour un niveau $m$ quelconque, e.g.,
les suites de Spencer ne sont plus exactes.

\begin{prop}\label{resspencerlog}
 Notons $\mathcal{T} _{X ^\#}: = ( \Omega ^1 _{X ^\#}) ^\vee $ le faisceau tangent de $X ^\#$
 et $\widetilde{\mathcal{T}} _{X ^\#} =\B _X \otimes _{\O _X} \mathcal{T} _{X ^\#}$.
Il existe un homomorphisme canonique $\widetilde{\mathcal{T}} _{X ^\#} \rightarrow \widetilde{\D} ^{(0)} _{X ^\#}$
induisant un isomorphisme
$\mathbb{S} (\widetilde{\mathcal{T}} _{X ^\#}) \riso \mathrm{gr} \widetilde{\D} ^{(0)} _{X ^\#}$.
\end{prop}
\begin{proof}
  Le cas où $X^\#$ est un log-schéma et $\B _X = \O _X$ a été vérifié par Montagnon
  dans \cite[5.1.1]{these_montagnon}. Le cas général se traite de la même façon.
\end{proof}

\begin{prop}\label{logstrdcn2}
  Supposons $X$ irréductible et
  $X ^\#$ muni de coordonnées locales logarithmiques
  $t _1, \dots, t _d \in M (Z)$.
  En notant $\sigma $ le morphisme canonique
  $\widetilde{\D} ^{(0)} _{X} \rightarrow \mathrm{Gr} \widetilde{\D} ^{(0)} _{X}$, on pose
  $\sigma _i := \sigma (\smash{\partial _\#} _i)$
  (on rappelle que $\widetilde{\D} ^{(0)} _{X ^\#}  \subset \widetilde{\D} ^{(0)} _{X } $).
   La suite $\{ \sigma _1, \dots ,\sigma _d \}$ est $\mathrm{Gr} \widetilde{\D} ^{(0)} _{X}$-régulière.
   On dispose de plus de l'égalité :
\begin{equation}
  \label{logstrdcn2eg}
  \sigma ( \widetilde{\D} ^{(0)} _{X} ( \smash{\partial _\#}  _1, \dots , \smash{\partial _\#} _d ))
  =
  \mathrm{Gr} \widetilde{\D} ^{(0)} _{X} (\sigma _1, \dots , \sigma _d ).
\end{equation}
\end{prop}

\begin{proof}
La première assertion découle du fait que $\mathrm{Gr} \widetilde{\D} ^{(0)} _{X}$ est une algèbre de polynôme en les
$\sigma (\smash{\partial } _1),\dots, \sigma (\smash{\partial } _d)$ et
que $t _1,\dots, t_d$ est une suite régulière de $\O _X$.
On vérifie alors l'égalité \ref{logstrdcn2eg} par un calcul identique à \cite[4.1.2]{calderon-logDderham}.
\end{proof}

\begin{defi}\label{firstspencer}
Soit $\E$ un $\widetilde{\D} ^{(0)} _{X ^\#}$-module à gauche
muni d'une filtration $\E = \cup _{s\in \N} \E _s$ (voir \ref{defbonnefilt}).
De manière analogue à \cite[1.6]{kashiwarathesis}, on définit un homomorphisme $\widetilde{\D} ^{(0)} _{X ^\#}$-linéaire
\begin{equation}
  \notag
  \epsilon\ :\
\widetilde{\D} ^{(0)} _{X ^\#} \otimes _{\B _{X}} \wedge ^{r} \widetilde{\mathcal{T}} _{X ^\#} \otimes _{\B _{X}}  \E _{s-1}
\rightarrow
\widetilde{\D} ^{(0)} _{X ^\#} \otimes _{\B _{X}} \wedge ^{r-1} \widetilde{\mathcal{T}} _{X ^\#} \otimes _{\B _{X}}  \E _s
\end{equation}
en posant, pour $i =1,\dots , r$,
$P \in \widetilde{\D} ^{(0)} _{X ^\#} $, $\delta _i \in \widetilde{\mathcal{T}} _{X ^\#}$,
$e \in \E _{s-1}$ :
\begin{align}
\notag
  \epsilon ( P \otimes (\delta _1 \wedge \cdots \wedge \delta _r)  \otimes e)
&  =  \sum _{i=1} ^{r} (-1) ^{i-1} P \delta _i \otimes
     (\delta _1
     \wedge \cdots \wedge
     \widehat{\delta _i}
     \wedge \cdots \wedge
     \delta _r )
     \otimes e
     \\
     \notag
 & -  \sum _{i=1} ^{r} (-1) ^{i-1} P  \otimes
     (\delta _1
     \wedge \cdots \wedge
     \widehat{\delta _i}
     \wedge \cdots \wedge
     \delta _r )
     \otimes \delta _i e
     \\
  & +  \sum _{1 \leq i < j\leq r} (-1) ^{i+j} P \otimes ([\delta _i , \delta _j] \wedge \delta _1
     \wedge \cdots \wedge
     \widehat{\delta _i}
     \wedge \cdots \wedge
     \widehat{\delta _j}
     \wedge \cdots \wedge
     \delta _r )
     \otimes e.
\end{align}
On vérifie par un calcul que l'on obtient le complexe
\begin{equation}
\label{sspencer}
  0 \rightarrow
\widetilde{\D} ^{(0)} _{X ^\#} \otimes _{\B _{X}} \wedge ^{d} \widetilde{\mathcal{T}} _{X ^\#} \otimes _{\B _{X}}  \E _{s-d}
\overset{\epsilon}{\rightarrow}
\cdots
\overset{\epsilon}{\rightarrow}
\widetilde{\D} ^{(0)} _{X ^\#} \otimes _{\B _{X}} \wedge ^{1} \widetilde{\mathcal{T}} _{X ^\#} \otimes _{\B _{X}}  \E _{s-1}
\overset{\epsilon}{\rightarrow}
\widetilde{\D} ^{(0)} _{X ^\#}  \otimes _{\B _{X}}  \E _s
\rightarrow
\E
\rightarrow
0
\end{equation}
que l'on écrira $Sp ^\bullet _{s,\widetilde{\D} ^{(0)} _{X ^\#}} (\E)$
et que l'on appellera {\og première suite de Spencer de degré $s$ de $\E$\fg}.

Lorsque $\E$ est $\B _{X}$-cohérent, on le munit de la filtration constante et
$Sp ^\bullet _{\widetilde{\D} ^{(0)} _{X ^\#}} (\E)$ indique la première suite de Spencer associée.

\end{defi}

\begin{rema}
Avec les notations de \ref{firstspencer},
on pose
$\M:= \widetilde{\omega} _{X ^\#} \otimes _{\B _{X}} \E$
et
$\M_s:= \widetilde{\omega} _{X ^\#} \otimes _{\B _{X}} \E_s$.
  On définit un homomorphisme de $\widetilde{\D} ^{(0)} _{X ^\#}$-modules à droite
$\M_{s-1}
\otimes _{\B _{X}}   \wedge ^{r} \widetilde{\mathcal{T}} _{X ^\#}
\otimes _{\B _{X}} \widetilde{\D} ^{(0)} _{X ^\#}
\rightarrow
\M_s \otimes _{\B _{X}} \wedge ^{r-1} \widetilde{\mathcal{T}} _{X ^\#}
\otimes _{\B _{X}} \widetilde{\D} ^{(0)} _{X ^\#}$
  via la commutativité du diagramme
  \begin{equation}
\xymatrix {
{\widetilde{\omega} _{X ^\#} \otimes _{\B _{X}} ( \widetilde{\D} ^{(0)} _{X ^\#} \otimes _{\B _{X}} \wedge ^{r} \widetilde{\mathcal{T}} _{X ^\#}
\otimes _{\B _{X}} \E _{s-1})}
\ar[r] ^{\epsilon} \ar[d] ^-\sim
&
{\widetilde{\omega} _{X ^\#} \otimes _{\B _{X}} (\widetilde{\D} ^{(0)} _{X ^\#} \otimes _{\B _{X}} \wedge ^{r-1} \widetilde{\mathcal{T}} _{X ^\#}
\otimes _{\B _{X}}  \E _s)} \ar[d] ^-\sim
\\
{\M_{s-1}
\otimes _{\B _{X}}   \wedge ^{r} \widetilde{\mathcal{T}} _{X ^\#}
\otimes _{\B _{X}} \widetilde{\D} ^{(0)} _{X ^\#}  }
\ar@{.>}[r] ^{\epsilon}
&
{\M_s \otimes _{\B _{X}} \wedge ^{r-1} \widetilde{\mathcal{T}} _{X ^\#}
\otimes _{\B _{X}} \widetilde{\D} ^{(0)} _{X ^\#}  ,}
}
  \end{equation}
  où les isomorphismes $\widetilde{\D} ^{(0)} _{X ^\#}$-linéaires verticaux se déduisent
  par fonctorialité de l'isomorphisme de transposition
  $\delta _\#$ : $\widetilde{\omega} _{X ^\#} \otimes _{\B _{X}} \widetilde{\D} ^{(0)} _{X ^\#} \riso
  \widetilde{\omega} _{X ^\#} \otimes _{\B _{X}} \widetilde{\D} ^{(0)} _{X ^\#}$ (voir \ref{isotranspd}).
On obtient alors le complexe
\begin{equation}
\label{sspencerd}
  0
  \rightarrow
 \M _{s-d} \otimes _{\B _{X}} \wedge ^{d} \widetilde{\mathcal{T}} _{X ^\#} \otimes _{\B _{X}} \widetilde{\D} ^{(0)} _{X ^\#}
\overset{\epsilon}{\rightarrow}
\cdots
\overset{\epsilon}{\rightarrow}
\M _{s-1} \otimes _{\B _{X}} \wedge ^{1} \widetilde{\mathcal{T}} _{X ^\#} \otimes _{\B _{X}} \widetilde{\D} ^{(0)} _{X ^\#}
\overset{\epsilon}{\rightarrow}
\M _s    \otimes _{\B _{X}}  \widetilde{\D} ^{(0)} _{X ^\#}
\rightarrow \M
\rightarrow
0.
\end{equation}
En le notant $Sp ^\bullet _{s,\widetilde{\D} ^{(0)} _{X ^\#}} (\M)$, on dispose par construction de l'isomorphisme
  $\widetilde{\omega} _{X ^\#} \otimes _{\B _{X}} Sp ^\bullet _{s,\widetilde{\D} ^{(0)} _{X ^\#}} (\E)
  \riso
  Sp ^\bullet _{s,\widetilde{\D} ^{(0)} _{X ^\#}} (\M)$.
\end{rema}

\begin{theo}
\label{theofirstspencer}
  Avec les notations de \ref{firstspencer},
  supposons de plus que la filtration de $\E$
  soit bonne.
   Alors, pour $s$ suffisamment grand,
  $Sp ^\bullet _{s,\widetilde{\D} ^{(0)} _{X ^\#}} (\E)$
  est exacte.
En particulier, si $\E$
est $\B _{X}$-cohérent,
le complexe $Sp ^\bullet _{\widetilde{\D} ^{(0)} _{X ^\#}} (\E)$
est acyclique.
\end{theo}
\begin{proof}
  De manière analogue au début de la preuve de \cite[1.6.1]{kashiwarathesis}
  (on utilise pour cela \ref{lemm-wildetildeDcoh}.\ref{lemm-wildetildeDcohi} pour
  vérifier que l'on obtient un complexe de Koszul),
on établit par récurrence sur $s\geq 0$ l'exactitude de la suite :
  \begin{equation}
  \label{theofirstspencerDpre}
   0 \rightarrow
\widetilde{\D} ^{(0)} _{X ^\#} \otimes _{\B _{X}} \wedge ^{d} \widetilde{\mathcal{T}} _{X ^\#} \otimes _{\B _{X}} \widetilde{\D} ^{(0)} _{X ^\#, s-d}
\overset{\epsilon}{\rightarrow} \cdots \overset{\epsilon}{\rightarrow}
\widetilde{\D} ^{(0)} _{X^\#} \otimes _{\B _{X}} \wedge ^{1} \widetilde{\mathcal{T}} _{X ^\#} \otimes _{\B _{X}}  \widetilde{\D} ^{(0)} _{X ^\#, s-1}
\overset{\epsilon}{\rightarrow}
  \widetilde{\D} ^{(0)} _{X^\#}  \otimes _{\B _{X}}  \widetilde{\D} ^{(0)} _{X ^\#,s}
\rightarrow \widetilde{\D} ^{(0)} _{X^\#}
\rightarrow
0.
\end{equation}
On conclut alors de manière analogue à la fin de la preuve de \cite[1.6.1]{kashiwarathesis}.
\end{proof}

\begin{vide}
  En appliquant $\widetilde{\D} ^{(0)} _{X} \otimes _{\widetilde{\D} ^{(0)} _{X ^\#} } -$ à \ref{sspencer},
  on obtient :
  \begin{gather}
\label{sspencerext}
  0 \rightarrow
\widetilde{\D} ^{(0)} _{X } \otimes _{\B _{X}} \wedge ^{d} \widetilde{\mathcal{T}} _{X ^\#} \otimes _{\B _{X}}  \E _{s-d}
\overset{\epsilon}{\rightarrow}
\cdots
\overset{\epsilon}{\rightarrow}
\widetilde{\D} ^{(0)} _{X} \otimes _{\B _{X}} \wedge ^{1} \widetilde{\mathcal{T}} _{X ^\#} \otimes _{\B _{X}}  \E _{s-1}
\overset{\epsilon}{\rightarrow}
\widetilde{\D} ^{(0)} _{X}  \otimes _{\B _{X}}  \E _s
\rightarrow \widetilde{\D} ^{(0)} _{X} \otimes _{\widetilde{\D} ^{(0)} _{X ^\#} } \E
\rightarrow
0.
\end{gather}
\end{vide}

\begin{theo}\label{theospencerext}
Soit $\E $ un $\widetilde{\D} ^{(0)} _{X ^\# } $-module à gauche cohérent qui soit cohérent et plat sur $\B _X$.
Le complexe de \ref{sspencerext} pour $\E$, $\widetilde{\D} ^{(0)} _{X} \otimes _{\widetilde{\D} ^{(0)} _{X ^\#} }
Sp ^\bullet _{\widetilde{\D} ^{(0)} _{X ^\#}} (\E)$,
 est acyclique.
\end{theo}

\begin{proof}
Il suffit de vérifier l'exactitude de la suite :
\begin{equation}
  \label{sspencerext-}
  0 \rightarrow
\widetilde{\D} ^{(0)} _{X } \otimes _{\B _{X}} \wedge ^{d} \widetilde{\mathcal{T}} _{X ^\#} \otimes _{\B _{X}}  \E
\overset{\epsilon}{\rightarrow}
\cdots
\overset{\epsilon}{\rightarrow}
\widetilde{\D} ^{(0)} _{X} \otimes _{\B _{X}} \wedge ^{1} \widetilde{\mathcal{T}} _{X ^\#} \otimes _{\B _{X}}  \E
\overset{\epsilon}{\rightarrow}
\widetilde{\D} ^{(0)} _{X}  \otimes _{\B _{X}}  \E .
\end{equation}
Comme $\E$ est plat sur $\B _X$,
on obtient la filtration de \ref{sspencerext-} suivante pour $n \in \N $ :
\begin{equation}
  \label{sspencerext-r}
  0 \rightarrow
\widetilde{\D} ^{(0)} _{X , n-d} \otimes _{\B _{X}} \wedge ^{d} \widetilde{\mathcal{T}} _{X ^\#} \otimes _{\B _{X}}  \E
\overset{\epsilon}{\rightarrow}
\cdots
\overset{\epsilon}{\rightarrow}
\widetilde{\D} ^{(0)} _{X,n-1} \otimes _{\B _{X}} \wedge ^{1} \widetilde{\mathcal{T}} _{X ^\#} \otimes _{\B _{X}}  \E
\overset{\epsilon}{\rightarrow}
\widetilde{\D} ^{(0)} _{X,n}  \otimes _{\B _{X}}  \E .
\end{equation}
D'après la preuve de \cite[4.1.3]{calderon-logDderham},
lorsque $\E$ est égal à $\B_X$, le gradué de la filtration \ref{sspencerext-r} donne une suite exacte.
De plus, on vérifie par un calcul immédiat que le gradué de la filtration \ref{sspencerext-r} est canoniquement isomorphe
au gradué de la filtration \ref{sspencerext-r} lorsque $\E$ est égal à $\B _X$ tensorisé par $\E$ au-dessus de $\B _X$.
Comme $\E$ est plat sur $\B _X$, ce gradué est donc une suite exacte. D'où l'exactitude de \ref{sspencerext-}.

\end{proof}

\begin{coro}
\label{corologdiag3}
Soit $\E $ un $\widetilde{\D} ^{(0)} _{X ^\# } $-module à gauche cohérent qui soit cohérent et plat sur $\B _X$.
L'homomorphisme canonique
\begin{gather}
  \label{resspencerextlogdiag3}
  \widetilde{\D} ^{(0)} _{X} \otimes _{\widetilde{\D} ^{(0)} _{X ^\#}} ^{\L} \E
  \rightarrow
  \widetilde{\D} ^{(0)} _{X} \otimes _{\widetilde{\D} ^{(0)} _{X ^\#}} \E
\end{gather}
est alors un isomorphisme.
\end{coro}

\begin{proof}
Cela découle de \ref{theofirstspencer} et \ref{theospencerext}.
\end{proof}

\begin{vide}
\label{remaQ}
Lorsque $X =\X$, toutes les définitions et tous les résultats de cette section s'étendent en remplaçant
$\widetilde{\D} ^{(0)} _{\X ^\#}$ (resp. $\widetilde{\D} ^{(0)} _{\X }$) par
$\widetilde{\D}  _{\X ^\#,\Q }:=\B _\X \otimes _{\O _\X} \D  _{\X ^\#,\Q}$
(resp. $\widetilde{\D}  _{\X ,\Q }:= \B _\X \otimes _{\O _\X} \D _{\X,\Q }$).
On obtient en particulier la proposition ci-après.
\end{vide}

\begin{prop}\label{corologdiag3q}
Soit $\E $ un $\widetilde{\D}  _{\X ^\#,\Q } $-module à gauche cohérent qui soit cohérent et plat sur $\B _{\X,\Q}$.
L'homomorphisme canonique
\begin{gather}
  \label{resspencerextlogdiag3q}
\widetilde{\D}  _{\X ,\Q }\otimes _{\widetilde{\D}  _{\X ^\#,\Q }} ^{\L} \E
  \rightarrow
\widetilde{\D}  _{\X ,\Q }\otimes _{\widetilde{\D}  _{\X ^\#,\Q }}  \E
\end{gather}
est alors un isomorphisme.
\end{prop}

\section{Un isomorphisme d'associativité}

Nous conservons les notations et hypothèses du chapitre \ref{nota1}.
Nous prouvons dans cette section l'isomorphisme d'{\og associativité\fg} \ref{A2CNiso}
puis par passage de gauche à droite celui de \ref{A2CNbis} (on remarque que l'appellation {\og associativité \fg}
est plus adéquate pour \ref{A2CNbis}).
L'isomorphisme \ref{A2CNiso} se généralise pour donner \ref{omegaEDF}, ce qui permettra d'établir
\ref{u+u!}.

\begin{theo}\label{A2CN}
  Soient $\E ^\#$ un $\widetilde{\D} ^{(m)} _{X^\#}$-module à gauche et $\M $ un $\widetilde{\D} ^{(m)} _{X}$-module à droite.
On dispose du morphisme de $\widetilde{\D} ^{(m)} _{X^\#}$-modules à droite :
  $\M\otimes _{\B _{X}} \E ^\# \rightarrow
\M\otimes _{\B _{X}} (\widetilde{\D} ^{(m)} _{X} \otimes _{\widetilde{\D} ^{(m)} _{X^\#}} \E ^\#)$,
envoyant, pour $m\in \M$ et $e \in \E ^\#$, $m \otimes e $ sur $m \otimes (1 \otimes e)$.
  Le morphisme $\widetilde{\D} ^{(m)} _{X}$-linéaire induit par extension :
\begin{equation}
  \label{A2CNiso}
(\M\otimes _{ \B _{X}} \E ^\#) \otimes _{\widetilde{\D} ^{(m)} _{X^\#}} \widetilde{\D} ^{(m)} _{X}
  \rightarrow
\M\otimes _{\B _{X}} ( \widetilde{\D} ^{(m)} _{X} \otimes _{\widetilde{\D} ^{(m)} _{X^\#}} \E ^\# )
\end{equation}
est un isomorphisme de $\widetilde{\D} ^{(m)} _{X}$-modules à droite.
\end{theo}

\begin{proof}On pourra comparer avec la preuve de \cite[A.1]{calderon-narvaez-dual}.
Le fait que la flèche \ref{A2CNiso} soit un isomorphisme est local.
Supposons donc $X ^\#$ muni de coordonnées locales logarithmiques $t _1,\dots, t_d$ et conservons
les notations de \ref{nota11}.
  Il s'agit de prouver que, pour tout $\widetilde{\D} ^{(m)} _{X}$-module à droite $\NN$,
  pour tout morphisme $\widetilde{\D} ^{(m)} _{X ^\#}$-linéaire $\alpha$ : $\M\otimes _{\B _{X}} \E ^\# \rightarrow \NN$,
  il existe un unique morphisme de $\widetilde{\D} ^{(m)} _{X}$-modules à droite
  $\beta$ : $\M\otimes _{\B _{X}} ( \widetilde{\D} ^{(m)} _{X} \otimes _{\widetilde{\D} ^{(m)} _{X^\#}} \E ^\# ) \rightarrow \NN$
  rendant commutatif le diagramme
\begin{equation}
  \notag
  \xymatrix {
  {\M\otimes _{\B _{X}} \E ^\#  } \ar[r]
  \ar@/_0,8cm/ [rr] ^-\alpha
  &
  { \M\otimes _{\B _{X}} (\widetilde{\D} ^{(m)} _{X} \otimes _{\widetilde{\D} ^{(m)} _{X^\#}} \E ^\#) } \ar@{.>}[r] ^-{\exists!\beta}
  &
  {\NN .}
  }
\end{equation}

Traitons d'abord l'unicité de $\beta$. Pour cela, on prouve par récurrence sur $N$ que, pour tous
$m \in \M$, $e \in \E ^\#$, $P \in \widetilde{\D} ^{(m)} _{X,N} $, $\alpha$ détermine de façon unique
l'élément $\beta ( m \otimes (P \otimes e))$.
Lorsque $N =0$, on a forcément $\beta (m \otimes (b \otimes e))= \alpha (m b \otimes e)$,
pour $b \in \B _{X}$.
Par linéarité, on peut supposer $P$ de la forme $\underline{\partial} ^{<\underline{k}>}$.
La formule suivante que doit vérifier $\beta$
(car $\beta$ est $\widetilde{\D} ^{(m)} _{X}$-linéaire et on utilise \cite[1.1.24.1]{caro_comparaison})
\begin{equation}
\notag
\beta (m \otimes ( \underline{\partial} ^{<\underline{k}>} \otimes  e)) =
(-1) ^{|\underline{k}|}
\left (
\beta ( m \otimes ( 1\otimes  e))\underline{\partial} ^{<\underline{k}>}
-
\sum _{\underline{h} <\underline{k}}
(-1) ^{|\underline{h}|}
\left \{ \begin{smallmatrix}   \underline{k} \\   \underline{h} \\ \end{smallmatrix} \right \}
\beta (m \underline{\partial} ^{<\underline{k}-\underline{h}>} \otimes ( \underline{\partial} ^{<\underline{h}>} \otimes  e ))
\right )
\end{equation}
nous permet de conclure la récurrence.

\'Etablissons à présent l'existence de $\beta$.
Par récurrence sur $N$, on construit un morphisme de groupes
$v _N$ : $ \M\times \widetilde{\D} ^{(m)} _{X,N} \times \E ^\#  \rightarrow \NN $ induisant $v _{N -1}$ de la manière suivante :
pour tous $m \in \M$, $e \in \E ^\#$, $b \in \B _{X}$, on pose
$v _0 (m, b,e) :=  \alpha (mb \otimes e) $.
En supposant défini $v _N$, pour tout
$\underline{k} \in \N ^d$ tel que $|\underline{k}| \leq N +1$,
pour tous $m \in \M$, $e \in \E ^\#$, $b \in \B _{X}$, on pose
\begin{equation}\label{A2CNeg1}
  v _{N+1} (m, \underline{\partial} ^{<\underline{k}>}, e)
:=
(-1) ^{|\underline{k}|}
\left (
v _N ( m, 1, e)\underline{\partial} ^{<\underline{k}>}
-
\sum _{\underline{h} <\underline{k}}
(-1) ^{|\underline{h}|}
\left \{ \begin{smallmatrix}   \underline{k} \\   \underline{h} \\ \end{smallmatrix} \right \}
v_N (m \underline{\partial} ^{<\underline{k}-\underline{h}>},\underline{\partial} ^{<\underline{h}>} , e )
\right ).
\end{equation}
Puis, pour tout $P = \sum _{\underline{r} }
b _{\underline{r}} \underline{\partial} ^{<\underline{r}>}
\in
\widetilde{\D} ^{(m)} _{X,N+1} $ où
$b _{\underline{r}}  \in \B _X$, on définit
\begin{equation}\label{A2CNeg1bis}
v _{N+1} (m, P, e) :=
\sum _{\underline{r}}
v _{N+1} (mb _{\underline{r}} ,\underline{\partial} ^{<\underline{r}>}, e).
\end{equation}
Si $P\in \widetilde{\D} ^{(m)} _{X,N} $, on remarque que $v _{N+1} (m, P, e) =v _{N} (m, P, e)$.
De plus, on vérifie que l'application $v _{N +1} $ est un morphisme de groupes.
Les morphismes $v _N$ induisent alors le morphisme de groupes
  $v$ : $ \M\times \widetilde{\D} ^{(m)} _{X} \times \E ^\#   \rightarrow \NN $.
Nous aurons besoin des lemmes ci-après.

\begin{lemm}
  \label{A2CN-lemm2pre}
    Pour tous $b \in \B _{X } $, $P \in \widetilde{\D} ^{(m)} _{X} $, $ m \in \M$, $e \in \E ^\#$,
\begin{gather}
  \label{A2CN-lemm2preeq2}
v (m, P, e)b = v (m b, P , e)= v (m, b P , e).
\end{gather}
\end{lemm}

\begin{proof}
L'égalité de droite de \ref{A2CN-lemm2preeq2} résulte de \ref{A2CNeg1bis}.
Par additivité, pour vérifier celle de gauche, on en déduit qu'il suffit d'établir
que
$\epsilon :=
 v (m, \underline{\partial} ^{<\underline{k}>} , e )b - v (mb, \underline{\partial} ^{<\underline{k}>} , e )$ est nul.
On procède par récurrence sur $N:=|\underline{k}|$.
On obtient par $\B _X$-linéarité de $\alpha$ :
$v (m, 1,e)b=v  (m b, 1,e)$.
Supposons à présent la formule vraie pour $N -1$.
Par \ref{A2CNeg1} puis par hypothèse de récurrence, on calcule :
\begin{align}
  \label{A2CN-lemm2prealign1}
v (m , 1, e ) \underline{\partial} ^{<\underline{k}>} b=
\sum _{\underline{i} \leq \underline{k}}
(-1) ^{|\underline{i}|}
\left \{ \begin{smallmatrix}  \underline{k} \\  \underline{i} \\ \end{smallmatrix} \right \}
v ( m \underline{\partial} ^{<\underline{k}-\underline{i}>} , \underline{\partial} ^{<\underline{i}>}, e)b
=
  (-1) ^{|\underline{k}|} \epsilon +
\sum _{\underline{i} \leq \underline{k}}
(-1) ^{|\underline{i}|}
\left \{ \begin{smallmatrix}  \underline{k} \\  \underline{i} \\ \end{smallmatrix} \right \}
v ( m \underline{\partial} ^{<\underline{k}-\underline{i}>} b, \underline{\partial} ^{<\underline{i}>}, e).
\end{align}
D'un autre côté, d'après \cite[2.3.5.1]{Be1}, on dispose de la formule :
$\underline{\partial} ^{<\underline{k}>} b =
\sum _{\underline{h} \leq \underline{k}}
\left \{ \begin{smallmatrix}  \underline{k} \\  \underline{h} \\ \end{smallmatrix} \right \}
\underline{\partial} ^{<\underline{k}-\underline{h}>}(b) \underline{\partial} ^{<\underline{h}>}$.
D'où :
\begin{align}
\notag
  v (m , 1, e ) \underline{\partial} ^{<\underline{k}>} b
  &=
\sum _{\underline{h} \leq \underline{k}}
\left \{ \begin{smallmatrix}  \underline{k} \\  \underline{h} \\ \end{smallmatrix} \right \}
v( m \underline{\partial} ^{<\underline{k}-\underline{h}>}(b) ,1,e) \underline{\partial} ^{<\underline{h}>}
\\ \notag
&=
\sum _{\underline{h} \leq \underline{k}}
\sum _{\underline{i} \leq \underline{h}}
(-1) ^{|\underline{i}|}
\left \{ \begin{smallmatrix}  \underline{k} \\  \underline{h} \\ \end{smallmatrix} \right \}
\left \{ \begin{smallmatrix}  \underline{h} \\  \underline{i} \\ \end{smallmatrix} \right \}
v(  m \underline{\partial} ^{<\underline{k}-\underline{h}>}(b)\underline{\partial} ^{<\underline{h}-\underline{i}>} ,
\underline{\partial} ^{<\underline{i}>},e)
& \text{(via \ref{A2CNeg1})}
\\ \notag
&=
\sum _{\underline{i} \leq \underline{h} \leq \underline{k}}
(-1) ^{|\underline{i}|}
\left \{ \begin{smallmatrix}  \underline{k} \\  \underline{i} \\ \end{smallmatrix} \right \}
\left \{ \begin{smallmatrix}  \underline{k} -\underline{i} \\  \underline{k} -\underline{h} \\ \end{smallmatrix} \right \}
v(  m \underline{\partial} ^{<\underline{k}-\underline{h}>}(b)\underline{\partial} ^{<\underline{h}-\underline{i}>} ,
\underline{\partial} ^{<\underline{i}>},e)
\\   \label{A2CN-lemm2prealign2}
 &=
\sum _{\underline{i} \leq \underline{k}}
(-1) ^{|\underline{i}|}
\left \{ \begin{smallmatrix}  \underline{k} \\  \underline{i} \\ \end{smallmatrix} \right \}
v ( m \underline{\partial} ^{<\underline{k}-\underline{i}>} b, \underline{\partial} ^{<\underline{i}>}, e).
&
\text{(à nouveau via \cite[2.3.5.1]{Be1})}
\end{align}
En comparant \ref{A2CN-lemm2prealign1} et \ref{A2CN-lemm2prealign2}, on obtient $\epsilon =0$.
\end{proof}

\begin{lemm}
    Pour tous $b \in \B _{X } $, $P \in \widetilde{\D} ^{(m)} _{X} $, $ m \in \M$, $e \in \E ^\#$,
\begin{gather}
  \label{A2CN-lemm2preeq}
  v (m, P, b e)=  v (m, P b, e).
\end{gather}
\end{lemm}

\begin{proof}

Grâce à la relation \ref{A2CN-lemm2preeq2}, on se ramène au cas où $P$ est de la forme
$\underline{\partial} ^{<\underline{k}>}$, avec $\underline{k} \in \N ^d$.
Il s'agit d'établir que
$\epsilon :=
v (m, \underline{\partial} ^{<\underline{k}>} , b e )- v (m, \underline{\partial} ^{<\underline{k}>}b  , e )$ est nul.
On procède par récurrence sur $N:=|\underline{k}|$. Lorsque $N  =0$, c'est évident.
Avec \cite[1.1.24]{caro_comparaison} et \cite[2.3.5.1]{Be1}, on calcule dans $\M\otimes _{\B _{X}} \widetilde{\D} ^{(m)} _{X}$
(on prend par défaut la structure gauche de $\widetilde{\D} ^{(m)} _{X}$-module à droite) :
 \begin{align}
   \label{A2CN-lemm2prealign11}
   (m b \otimes 1) \underline{\partial} ^{<\underline{k}>}
   & =
\sum _{\underline{h} \leq \underline{k}}
(-1) ^{|\underline{h}|}
\left \{ \begin{smallmatrix}  \underline{k} \\  \underline{h} \\ \end{smallmatrix} \right \}
m b\underline{\partial} ^{<\underline{k}-\underline{h}>} \otimes \underline{\partial} ^{<\underline{h}>},
\\
   (m  \otimes b ) \underline{\partial} ^{<\underline{k}>}
   & =
\sum _{\underline{r} \leq \underline{k}}
(-1) ^{|\underline{r}|}
\left \{ \begin{smallmatrix}  \underline{k} \\  \underline{r} \\ \end{smallmatrix} \right \}
m \underline{\partial} ^{<\underline{k}-\underline{r}>} \otimes \underline{\partial} ^{<\underline{r}>} b
   \label{A2CN-lemm2prealign12}
 =
\sum _{\underline{r} \leq \underline{k},\, \underline{s} \leq \underline{r} }
(-1) ^{|\underline{r}|}
\left \{ \begin{smallmatrix}  \underline{k} \\  \underline{r} \\ \end{smallmatrix} \right \}
\left \{ \begin{smallmatrix}  \underline{r} \\  \underline{s} \\ \end{smallmatrix} \right \}
m \underline{\partial} ^{<\underline{k}-\underline{r}>} \underline{\partial} ^{<\underline{r}-\underline{s}>}( b )
\otimes
\underline{\partial} ^{<\underline{s}>}.
 \end{align}
Comme $m b \otimes 1 = m \otimes b$ et que $\widetilde{\D} ^{(m)} _{X}$ est un $\B _{X}$-module
(pour la structure gauche ou droite)
libre de base les
$\underline{\partial} ^{<\underline{n}>}$ avec $\underline{n}$ parcourant $\N ^d$,
il découle de \ref{A2CN-lemm2prealign11} et \ref{A2CN-lemm2prealign12}
la relation dans $\M\times \widetilde{\D} ^{(m)} _{X}$ :
\begin{equation}
  \label{A2CN-lemm2pre=}
\sum _{\underline{h} \leq \underline{k}}
(         (-1) ^{|\underline{h}|}
\left \{ \begin{smallmatrix}  \underline{k} \\  \underline{h} \\ \end{smallmatrix} \right \}
m b\underline{\partial} ^{<\underline{k}-\underline{h}>} , \underline{\partial} ^{<\underline{h}>})
=
\sum _{\underline{r} \leq \underline{k},\, \underline{s} \leq \underline{r} }
(   (-1) ^{|\underline{r}|}
\left \{ \begin{smallmatrix}  \underline{k} \\  \underline{r} \\ \end{smallmatrix} \right \}
\left \{ \begin{smallmatrix}  \underline{r} \\  \underline{s} \\ \end{smallmatrix} \right \}
m \underline{\partial} ^{<\underline{k}-\underline{r}>} \underline{\partial} ^{<\underline{r}-\underline{s}>}( b )
,
\underline{\partial} ^{<\underline{s}>}  ).
\end{equation}
Par successivement \ref{A2CNeg1}, hypothèse de récurrence, \cite[2.3.5.1]{Be1} et \ref{A2CN-lemm2pre},
\ref{A2CN-lemm2pre=} puis \ref{A2CNeg1},
on vérifie les égalités :
\begin{align}
\notag
  v (m , 1, ae ) \underline{\partial} ^{<\underline{k}>}
  &
= \sum _{\underline{r} \leq \underline{k}}
(-1) ^{|\underline{r}|}
\left \{ \begin{smallmatrix}  \underline{k} \\  \underline{r} \\ \end{smallmatrix} \right \}
v ( m \underline{\partial} ^{<\underline{k}-\underline{r}>} , \underline{\partial} ^{<\underline{r}>} ,b  e)
\\ \notag
& =(-1) ^{|\underline{k}|}  \epsilon + \sum _{\underline{r} \leq \underline{k}}
(-1) ^{|\underline{r}|}
\left \{ \begin{smallmatrix}  \underline{k} \\  \underline{r} \\ \end{smallmatrix} \right \}
v ( m \underline{\partial} ^{<\underline{k}-\underline{r}>} , \underline{\partial} ^{<\underline{r}>} b   ,e)
\\ \notag
& =  (-1) ^{|\underline{k}|} \epsilon +
\sum _{\underline{r} \leq \underline{k},\, \underline{s} \leq \underline{r} }
(-1) ^{|\underline{r}|}
\left \{ \begin{smallmatrix}  \underline{k} \\  \underline{r} \\ \end{smallmatrix} \right \}
\left \{ \begin{smallmatrix}  \underline{r} \\  \underline{s} \\ \end{smallmatrix} \right \}
v ( m \underline{\partial} ^{<\underline{k}-\underline{r}>} \underline{\partial} ^{<\underline{r}-\underline{s}>}( b )
, \underline{\partial} ^{<\underline{s}>} , e)
\\ \notag
& =  (-1) ^{|\underline{k}|} \epsilon +
\sum _{\underline{h} \leq \underline{k}}
(-1) ^{|\underline{h}|}
\left \{ \begin{smallmatrix}  \underline{k} \\  \underline{h} \\ \end{smallmatrix} \right \}
v ( m b\underline{\partial} ^{<\underline{k}-\underline{h}>} , \underline{\partial} ^{<\underline{h}>}, e)
\\ \notag
& =  (-1) ^{|\underline{k}|} \epsilon +
  v (m b, 1, e ) \underline{\partial} ^{<\underline{k}>}
  =  (-1) ^{|\underline{k}|} \epsilon +
  v (m , 1, ae ) \underline{\partial} ^{<\underline{k}>}.
\end{align}
D'où $\epsilon  =0$.
\end{proof}

\begin{lemm}
\label{A2CN-lemm1}
Pour tous $P \in \widetilde{\D} ^{(m)} _{X} $, $ m \in \M$, $e \in \E ^\#$, $\underline{k} \in \N ^d$, la formule suivante est validée :
\begin{equation}
  \label{A2CNeg2}
v  ( m, P, e)\underline{\partial} ^{<\underline{k}>}
=
\sum _{\underline{h} \leq \underline{k}}
(-1) ^{|\underline{h}|}
\left \{\begin{smallmatrix}  \underline{k} \\  \underline{h} \\ \end{smallmatrix} \right \}
v (m \underline{\partial} ^{<\underline{k}-\underline{h}>},\underline{\partial} ^{<\underline{h}>} P , e ).
\end{equation}
\end{lemm}

\begin{proof}
Comme $P $ s'écrit sous la forme
$\sum \underline{\partial} ^{<\underline{r}>} b _{\underline{r}} $, où
$b _{\underline{r}} \in \B _X$,
par additivité et \ref{A2CN-lemm2preeq}, on se ramène au cas où $P=\underline{\partial} ^{<\underline{r}>}$.
 Prouvons le lemme par récurrence sur $N:=|\underline{r}|$.
  Pour $N =0$, cela découle de \ref{A2CNeg1}. 
Supposons l'égalité validée pour $N $, prouvons-le pour $N +1$.

Par
hypothèse de récurrence, on vérifie les égalités ci-dessous :
\begin{align}
\notag
 \bullet \  v (m, 1 , e)\underline{\partial} ^{<\underline{r}>} \underline{\partial} ^{<\underline{k}>}
  =&
(-1) ^{|\underline{r}|}
v (m ,\underline{\partial} ^{<\underline{r}>} , e )
\underline{\partial} ^{<\underline{k}>}
+
\sum _{\underline{s} <\underline{r}}
(-1) ^{|\underline{s}|}
\left \{\begin{smallmatrix}  \underline{r} \\  \underline{s} \\ \end{smallmatrix} \right \}
v (m \underline{\partial} ^{<\underline{r}-\underline{s}>},\underline{\partial} ^{<\underline{s}>} , e )
\underline{\partial} ^{<\underline{k}>},
\\
    \label{A2CNeg3}
=&
(-1) ^{|\underline{r}|}
v (m ,\underline{\partial} ^{<\underline{r}>} , e )
\underline{\partial} ^{<\underline{k}>}
+
\sum _{\underline{h} \leq \underline{k}}
\sum _{\underline{s} <\underline{r}}
(-1) ^{|\underline{h}|+|\underline{s}|}
\left \{\begin{smallmatrix}  \underline{k} \\  \underline{h} \\ \end{smallmatrix} \right \}
\left <\begin{smallmatrix}  \underline{h} + \underline{s} \\  \underline{s} \\ \end{smallmatrix} \right >
\left \{ \begin{smallmatrix}   \underline{r} \\  \underline{s} \\ \end{smallmatrix} \right \}
v (m \underline{\partial} ^{<\underline{r}-\underline{s}>}\underline{\partial} ^{<\underline{k}-\underline{h}>},
\underline{\partial} ^{<\underline{h}+\underline{s}>} , e ),
\\  \label{A2CNeg4}
 \bullet \   v (m, 1 , e)\underline{\partial} ^{<\underline{r}>} \underline{\partial} ^{<\underline{k}>}
  =&
\left < \begin{smallmatrix}   \underline{r} +\underline{k} \\   \underline{k} \\ \end{smallmatrix} \right >
 v (m, 1 , e)\underline{\partial} ^{<\underline{r} +\underline{k} >}
=
\sum _{\underline{l} \leq \underline{r} +\underline{k}}
(-1) ^{|\underline{l} |}
\left < \begin{smallmatrix}  \underline{r} +\underline{k} \\  \underline{k} \\ \end{smallmatrix} \right >
\left \{ \begin{smallmatrix} \underline{r} +\underline{k}\\   \underline{l} \\ \end{smallmatrix} \right \}
v (m \underline{\partial} ^{<\underline{r} +\underline{k}-\underline{l}>},
\underline{\partial} ^{<\underline{l}>} , e ).
\end{align}

 D'un autre côté, on calcule de même
 $(m \otimes 1 )\underline{\partial} ^{<\underline{r}>} \underline{\partial} ^{<\underline{k}>}
\in \M\otimes _{\B _{X}} \widetilde{\D} ^{(m)} _{X}$ (on prend par défaut la structure
 gauche de $\widetilde{\D} ^{(m)} _{X}$-module à droite)
des deux différentes façons :
\small
\begin{align}
  \bullet\  (m \otimes 1 )\underline{\partial} ^{<\underline{r}>} \underline{\partial} ^{<\underline{k}>}
  =&
\sum _{\underline{s} \leq \underline{r}}
(-1) ^{|\underline{s}|}
\left \{ \begin{smallmatrix}   \underline{r} \\   \underline{s} \\ \end{smallmatrix} \right \}
(m \underline{\partial} ^{<\underline{r}-\underline{s}>} \otimes \underline{\partial} ^{<\underline{s}>})
\underline{\partial} ^{<\underline{k}>}
\label{A2CNeg5}
=
\sum _{\underline{h} \leq \underline{k}}
\sum _{\underline{s} \leq \underline{r}}
(-1) ^{|\underline{h}|+|\underline{s}|}
\left \{ \begin{smallmatrix}   \underline{k} \\   \underline{h} \\ \end{smallmatrix} \right \}
\left <\begin{smallmatrix}  \underline{h} + \underline{s} \\  \underline{s} \\ \end{smallmatrix} \right >
\left \{ \begin{smallmatrix}   \underline{r} \\   \underline{s} \\ \end{smallmatrix} \right \}
m \underline{\partial} ^{<\underline{r}-\underline{s}>}\underline{\partial} ^{<\underline{k}-\underline{h}>}
\otimes
\underline{\partial} ^{<\underline{h}+\underline{s}>},
\\
\label{A2CNeg6}
\bullet\
(m \otimes  1 )\underline{\partial} ^{<\underline{r}>} \underline{\partial} ^{<\underline{k}>}
  = &
\left < \begin{smallmatrix}   \underline{r} +\underline{k} \\   \underline{k} \\ \end{smallmatrix} \right >
(m \otimes  1)
\underline{\partial} ^{<\underline{r} +\underline{k} >}
=
\sum _{\underline{l} \leq \underline{r} +\underline{k}}
(-1) ^{|\underline{l} |}
\left < \begin{smallmatrix}   \underline{r} +\underline{k} \\   \underline{k} \\ \end{smallmatrix} \right >
\left \{\begin{smallmatrix} \underline{r} +\underline{k}\\   \underline{l} \\ \end{smallmatrix} \right \}
m \underline{\partial} ^{<\underline{r} +\underline{k}-\underline{l}>}
\otimes
\underline{\partial} ^{<\underline{l}>}.
\end{align}
\normalsize
Il résulte de \ref{A2CNeg5} et \ref{A2CNeg6} l'égalité dans $\M \times \widetilde{\D} ^{(m)} _{X}$ :
\begin{equation}
  \label{A2CNegcup=cup}
  \sum  _{\underline{l} \leq \underline{r} +\underline{k}}
((-1) ^{|\underline{l} |}
\left < \begin{smallmatrix}   \underline{r} +\underline{k} \\   \underline{k} \\ \end{smallmatrix} \right >
\left \{ \begin{smallmatrix} \underline{r} +\underline{k}\\  \underline{l} \\ \end{smallmatrix} \right \}
m \underline{\partial} ^{<\underline{r} +\underline{k}-\underline{l}>}
,
\underline{\partial} ^{<\underline{l}>})
=
\sum _{\underline{h} \leq \underline{k},\ \underline{s} \leq \underline{r}}
((-1) ^{|\underline{h}|+|\underline{s}|}
\left \{ \begin{smallmatrix}  \underline{k} \\  \underline{h} \\ \end{smallmatrix} \right \}
\left <\begin{smallmatrix}  \underline{h} + \underline{s} \\  \underline{s} \\ \end{smallmatrix} \right >
\left \{ \begin{smallmatrix}   \underline{r} \\   \underline{s} \\ \end{smallmatrix} \right \}
m \underline{\partial} ^{<\underline{r}-\underline{s}>}\underline{\partial} ^{<\underline{k}-\underline{h}>}
,
\underline{\partial} ^{<\underline{h}+\underline{s}>}).
\end{equation}
On déduit de \ref{A2CNegcup=cup} et \ref{A2CNeg4} la formule :
\begin{equation}
    \label{A2CNeg7}
  v (m, 1 , e)\underline{\partial} ^{<\underline{r}>} \underline{\partial} ^{<\underline{k}>}
    =
  \sum _{\underline{h} \leq \underline{k}}
\sum _{\underline{s} \leq \underline{r}}
(-1) ^{|\underline{h}|+|\underline{s}|}
\left \{ \begin{smallmatrix}   \underline{k} \\   \underline{h} \\ \end{smallmatrix} \right \}
\left <\begin{smallmatrix}  \underline{h} + \underline{s} \\  \underline{s} \\ \end{smallmatrix} \right >
\left \{ \begin{smallmatrix}   \underline{r} \\   \underline{s} \\ \end{smallmatrix} \right \}
v (m \underline{\partial} ^{<\underline{r}-\underline{s}>}\underline{\partial} ^{<\underline{k}-\underline{h}>},
\underline{\partial} ^{<\underline{h}+\underline{s}>} , e ).
\end{equation}
Il résulte de \ref{A2CNeg3} et \ref{A2CNeg7}
\begin{equation}
(-1) ^{|\underline{r}|}
  v  ( m, \underline{\partial} ^{<\underline{r}>}, e)\underline{\partial} ^{<\underline{k}>}
=
(-1) ^{|\underline{r}|}
\sum _{\underline{h} \leq \underline{k}}
(-1) ^{|\underline{h}|}
\left \{ \begin{smallmatrix}   \underline{k} \\   \underline{h} \\ \end{smallmatrix} \right \}
v (m \underline{\partial} ^{<\underline{k}-\underline{h}>},
\underline{\partial} ^{<\underline{h}>} \underline{\partial} ^{<\underline{r}>} , e ).
\end{equation}
La formule \ref{A2CNeg2} est donc vérifiée pour $P =\underline{\partial} ^{<\underline{r}>}$.
\end{proof}

\begin{lemm}\label{A2CN-lemm2}
  Pour tous $P ^\# \in \widetilde{\D} ^{(m)} _{X ^\#} $, $P \in \widetilde{\D} ^{(m)} _{X} $, $ m \in \M$, $e \in \E ^\#$,
\begin{equation}
  \label{A2CN-lemm2eq}
  v (m, P, P ^\# e)=  v (m, P P ^\#, e).
\end{equation}
\end{lemm}

\begin{proof}
Dans un premier temps, supposons $P =1$, i.e., vérifions l'égalité
$v  ( m, 1, P ^\# e) = v  ( m, P ^\# , e)$.
Par récurrence sur $|\underline{k}|$, établissons d'abord la formule :
$v  ( m, 1, \underline{\partial} _\# ^{<\underline{k}>}  e)
=v  ( m, \underline{\partial} _\# ^{<\underline{k}>} ,  e)$.

En multipliant dans l'égalité \ref{A2CNeg2} par $t _1  ^{k _1} \dots t _d ^{k _d}$
et grâce à \ref{A2CN-lemm2preeq2}, on obtient (avec $P =1$) :
\begin{equation}
  \label{A2CNeg2bis}
v  ( m, 1, e)\overset{^\mathrm{t}}{} \underline{\partial} ^{<\underline{k}>}_\#
=
\sum _{\underline{h} \leq \underline{k}}
\left \{\begin{smallmatrix}  \underline{k} \\  \underline{h} \\ \end{smallmatrix} \right \}
v (m \overset{^\mathrm{t}}{} \underline{\partial} ^{<\underline{k}-\underline{h}>}_\# ,
\underline{\partial} _\# ^{<\underline{h}>}  , e ).
\end{equation}
Or, par $\widetilde{\D} ^{(m)} _{X ^\#}$-linéarité de $\alpha$, il découle de \ref{MotimesEetcegbis} :
\begin{equation}
  \label{MotimesEetcegcons}
 ( \alpha  (m \otimes e) )\overset{^\mathrm{t}}{} \underline{\partial} ^{<\underline{k}>}_\#
=
\sum _{\underline{h} \leq \underline{k}}
\left \{ \begin{smallmatrix}   \underline{k} \\   \underline{h} \\ \end{smallmatrix} \right \}
\alpha (m \overset{^\mathrm{t}}{} \underline{\partial} ^{<\underline{k}-\underline{h}>}_\#
\otimes
\underline{\partial} ^{<\underline{h}>}_\# e).
\end{equation}
Par hypothèse de récurrence, pour tout $\underline{h} < \underline{k}$,
$v (m \overset{^\mathrm{t}}{} \underline{\partial} ^{<\underline{k}-\underline{h}>}_\# ,
\underline{\partial} _\# ^{<\underline{h}>} , e )=
v (m \overset{^\mathrm{t}}{} \underline{\partial} ^{<\underline{k}-\underline{h}>}_\# ,
1 , \underline{\partial} _\# ^{<\underline{h}>} e )=
\alpha (m \overset{^\mathrm{t}}{} \underline{\partial} ^{<\underline{k}-\underline{h}>}_\#
\otimes
\underline{\partial} ^{<\underline{h}>}_\# e)$.
En comparant \ref{A2CNeg2bis} et \ref{MotimesEetcegcons}, on en conclut
que
$v  ( m, \underline{\partial} _\# ^{<\underline{k}>} ,  e)=
\alpha (m \otimes
\underline{\partial} ^{<\underline{k}>}_\# e)
=v  ( m, 1, \underline{\partial} _\# ^{<\underline{k}>}  e)$.

Enfin, si $P ^\# $ est de la forme $\sum _{\underline{k}} b _{\underline{k}} \underline{\partial} _\# ^{<\underline{k}>}$,
où $b _{\underline{k}}  \in \B _X$,
par linéarité de $\alpha$,
ce que l'on vient d'établir,
puis \ref{A2CN-lemm2preeq2}, on obtient :
$v  ( m, 1, P ^\#  e)
= \sum _{\underline{k}} v  ( m b _{\underline{k}},1 ,  \underline{\partial} _\# ^{<\underline{k}>}  e)
=\sum _{\underline{k}} v  ( m b _{\underline{k}},   \underline{\partial}_\#  ^{<\underline{k}>},  e)
= v( m ,  P ^\# ,  e)$.

Traitons à présent le cas général. Par \ref{A2CN-lemm2preeq2}, il suffit de vérifier \ref{A2CN-lemm2eq}
lorsque $P$ est de la forme
$\underline{\partial} ^{<\underline{k}>}$, avec $\underline{k} \in \N ^d$.
Effectuons alors une récurrence sur l'entier $N:=|\underline{k}|$.
Pour $N =0$, il s'agit de ce que l'on vient de prouver.
Supposons le lemme vrai pour $N$ et supposons $|\underline{k}| \leq N +1$.
D'après \ref{A2CNeg1} :
  \begin{equation}
  \label{A2CNeg1-bis}
  v (m, \underline{\partial} ^{<\underline{k}>}, P ^\# e)
=
(-1) ^{|\underline{k}|}
\left (
v  ( m, 1, P ^\# e)\underline{\partial} ^{<\underline{k}>}
-
\sum _{\underline{h} <\underline{k}}
(-1) ^{|\underline{h}|}
\left \{ \begin{smallmatrix}   \underline{k} \\   \underline{h} \\ \end{smallmatrix} \right \}
v (m \underline{\partial} ^{<\underline{k}-\underline{h}>},\underline{\partial} ^{<\underline{h}>} , P ^\# e )
\right ).
\end{equation}
De plus, il dérive de \ref{A2CNeg2} la formule :
  \begin{equation}
  \label{A2CNeg2-bis}
  v (m, \underline{\partial} ^{<\underline{k}>}P ^\#  , e)
=
(-1) ^{|\underline{k}|}
\left (
v  ( m, P ^\# , e)\underline{\partial} ^{<\underline{k}>}
-
\sum _{\underline{h} <\underline{k}}
(-1) ^{|\underline{h}|}
\left \{ \begin{smallmatrix}   \underline{k} \\   \underline{h} \\\end{smallmatrix}\right \}
v (m \underline{\partial} ^{<\underline{k}-\underline{h}>},\underline{\partial} ^{<\underline{h}>} P ^\# , e )
\right ).
\end{equation}
Via \ref{A2CNeg1-bis}, \ref{A2CNeg2-bis} et par hypothèses de récurrence, on établit
$  v (m, \underline{\partial} ^{<\underline{k}>} ,P ^\#  e) =
  v (m, \underline{\partial} ^{<\underline{k}>}P ^\#  , e)$.
\end{proof}
Concluons maintenant la preuve du théorème.
Il dérive de \ref{A2CN-lemm2preeq2} et \ref{A2CN-lemm2} que le morphisme $v$ induit un morphisme
de $\B _X$-modules
$\beta $ :
$\M\otimes _{\B _{X}} ( \widetilde{\D} ^{(m)} _{X} \otimes _{\widetilde{\D} ^{(m)} _{X^\#}} \E ^\# ) \rightarrow \NN$.
Enfin, il résulte des formules \ref{A2CNeg2} et de \cite[1.1.24.1]{caro_comparaison} que $\beta$
est $\widetilde{\D} ^{(m)} _{X}$-linéaire.
\end{proof}

\begin{theo}\label{A2CNbis}
  Soient $\E ^\#$ un $\widetilde{\D} ^{(m)} _{X^\#}$-module à gauche et $\FF $ un $\widetilde{\D} ^{(m)} _{X}$-module à gauche.
  On dispose du morphisme de $\widetilde{\D} ^{(m)} _{X^\#}$-modules à gauche
  $\E ^\#\otimes _{ \B _{X}} \FF
  \rightarrow
(\widetilde{\D} ^{(m)} _{X} \otimes _{\widetilde{\D} ^{(m)} _{X^\#}} \E ^\# )\otimes _{ \B _{X}} \FF$,
envoyant $e \otimes f $ sur $(1 \otimes e )\otimes f$ où
$e \in \E ^\#$, $f \in \FF$.
Le morphisme $\widetilde{\D} ^{(m)} _{X}$-linéaire induit par extension :
\begin{equation}
  \label{A2CNbisiso}
\widetilde{\D} ^{(m)} _{X} \otimes _{\widetilde{\D} ^{(m)} _{X^\#}} (\E ^\#\otimes _{ \B _{X}} \FF)
  \rightarrow
(\widetilde{\D} ^{(m)} _{X} \otimes _{\widetilde{\D} ^{(m)} _{X^\#}} \E ^\# )\otimes _{ \B _{X}} \FF
\end{equation}
est un isomorphisme de $\widetilde{\D} ^{(m)} _{X}$-modules à gauche.
\end{theo}

\begin{proof}
  On procède de façon analogue à \ref{A2CN}.
\end{proof}

\begin{rema}
\label{rema37}
Soient $\E ^\#$ un $\widetilde{\D} ^{(m)} _{X^\#}$-module à gauche,
$\FF $ un $\widetilde{\D} ^{(m)} _{X}$-module à gauche
(resp. un $\widetilde{\D} ^{(m)} _{X}$-bimodule).
  Il découle de \ref{A2CN} l'isomorphisme canonique
  $(\widetilde{\omega} _X \otimes _{ \B _{X}} \E ^\#) \otimes _{\widetilde{\D} ^{(m)} _{X^\#}} \widetilde{\D} ^{(m)} _{X}
  \riso
\widetilde{\omega} _X \otimes _{\B _{X}} ( \widetilde{\D} ^{(m)} _{X} \otimes _{\widetilde{\D} ^{(m)} _{X^\#}} \E ^\# )$.
En lui appliquant le foncteur $- \otimes _{\widetilde{\D} ^{(m)} _{X}} \FF$, cela donne :
$(\widetilde{\omega} _X \otimes _{ \B _{X}} \E ^\#) \otimes _{\widetilde{\D} ^{(m)} _{X^\#}} \FF
  \riso
[\widetilde{\omega} _X \otimes _{\B _{X}} ( \widetilde{\D} ^{(m)} _{X} \otimes _{\widetilde{\D} ^{(m)} _{X^\#}} \E ^\# )] \otimes _{\widetilde{\D} ^{(m)} _{X}} \FF$.
L'isomorphisme de transposition
$\delta$ :
  $\widetilde{\omega} _X \otimes _{\B _X} \widetilde{\D} ^{(m)} _X \riso
  \widetilde{\omega} _X \otimes _{\B _X} \widetilde{\D} ^{(m)} _X $,
  qui échange les deux structures de $\widetilde{\D} ^{(m)} _X $-modules à droite,
  induit par fonctorialité
$[(\widetilde{\omega} _X \otimes _{\B _{X}} \widetilde{\D} ^{(m)} _{X} )\otimes _{\widetilde{\D} ^{(m)} _{X^\#}} \E ^\# ] \otimes _{\widetilde{\D} ^{(m)} _{X}} \FF
\riso
(\widetilde{\omega} _X \otimes _{\B _{X}} \FF )\otimes _{\widetilde{\D} ^{(m)} _{X^\#}} \E ^\# $.
On obtient par composition l'isomorphisme de groupes (resp. de $\widetilde{\D} ^{(m)} _{X}$-modules à droite) :
\begin{equation}
\label{omegaEDF}
  (\widetilde{\omega} _X \otimes _{ \B _{X}} \E ^\#) \otimes _{\widetilde{\D} ^{(m)} _{X^\#}} \FF
  \riso
  (\widetilde{\omega} _X \otimes _{\B _{X}} \FF )\otimes _{\widetilde{\D} ^{(m)} _{X^\#}} \E ^\# .
\end{equation}
\end{rema}

\begin{rema}
  \label{rema38}
Soit $\B _\X$ une $\O _\X$-algèbre commutative vérifiant les hypothèses de \ref{chap-coh-21}.
Le faisceau $\B _\X \smash{\widehat{\otimes}} _{\O _{\X}} \widehat{\D} ^{(m)} _{\X ^\#}$ est cohérent
(voir \ref{cohhat}).
Soient $\E ^\#$ un $\B _\X \smash{\widehat{\otimes}} _{\O _{\X}} \widehat{\D} ^{(m)} _{\X ^\#}$-module à gauche,
$\FF $ un $\B _\X \smash{\widehat{\otimes}} _{\O _{\X}} \widehat{\D} ^{(m)} _{\X }$-module à gauche
(resp. un $\B _\X \smash{\widehat{\otimes}} _{\O _{\X}} \widehat{\D} ^{(m)} _{\X }$-bimodule).

De manière analogue à \ref{dr-gabis},
on vérifie que
les foncteurs $-\otimes _{\B _{\X}} \widetilde{\omega}_{\X ^\# } ^{-1}$
et $\widetilde{\omega}_{\X ^\# }  \otimes _{\B _{\X}} -$
(resp. $-\otimes _{\B _{\X}} \widetilde{\omega}_{\X  } ^{-1}$
et $\widetilde{\omega}_{\X }  \otimes _{\B _{\X}} -$) induisent des équivalences quasi-inverses exactes entre
la catégorie des $\B _\X \smash{\widehat{\otimes}} _{\O _{\X}} \widehat{\D} ^{(m)} _{\X ^\#}$-modules
(resp. cohérents, resp. plats, resp. localement projectifs de type fini)
à gauche et celle des
$\B _\X \smash{\widehat{\otimes}} _{\O _{\X}} \widehat{\D} ^{(m)} _{\X ^\#}$-modules
(resp. cohérents, resp. plats, resp. localement projectifs de type fini) à droite.
On obtient alors par extension (comme pour \ref{A2CN}) l'homomorphisme de
$\B _\X \smash{\widehat{\otimes}} _{\O _{\X}} \widehat{\D} ^{(m)} _{\X }$-modules cohérents à droite :
  $$(\widetilde{\omega} _\X \otimes _{ \B _{\X}} \E ^\#)
  \otimes _{\B _\X \smash{\widehat{\otimes}} _{\O _{\X}} \widehat{\D} ^{(m)} _{\X ^\#}}
  \B _\X \smash{\widehat{\otimes}} _{\O _{\X}} \widehat{\D} ^{(m)} _{\X }
  \rightarrow
\widetilde{\omega} _\X \otimes _{\B _{\X}} ( \B _\X \smash{\widehat{\otimes}} _{\O _{\X}} \widehat{\D} ^{(m)} _{\X }
\otimes _{\B _\X \smash{\widehat{\otimes}} _{\O _{\X}} \widehat{\D} ^{(m)} _{\X ^\#}} \E ^\# ),$$
défini pour $x\in \widetilde{\omega} _\X $, $e\in \E ^\#$ par
$(x \otimes e ) \otimes 1 \mapsto x \otimes ( 1\otimes e )$.
Celui-ci est en fait un isomorphisme puisqu'il l'est modulo $\mathfrak{m} ^{i +1}$.
On en déduit comme pour \ref{omegaEDF} l'isomorphisme canonique
de groupes (resp. de $\B _\X \smash{\widehat{\otimes}} _{\O _{\X}} \widehat{\D} ^{(m)} _{\X }$-modules à droite) :
\begin{equation}
\label{omegaEDFhat}
(\widetilde{\omega} _\X \otimes _{ \B _{\X}} \E ^\#)
  \otimes _{\B _\X \smash{\widehat{\otimes}} _{\O _{\X}} \widehat{\D} ^{(m)} _{\X ^\#}}
  \FF
  \rightarrow
\widetilde{\omega} _\X \otimes _{\B _{\X}} (\FF
\otimes _{\B _\X \smash{\widehat{\otimes}} _{\O _{\X}} \widehat{\D} ^{(m)} _{\X ^\#}} \E ^\# ).
\end{equation}
\end{rema}

\section{Log-isocristaux surconvergents}

\begin{vide}
\label{cohhat}
On garde les notations et hypothèses de \ref{chap-coh-21}.
Nous sommes ainsi dans le contexte de la section \cite[3.3]{Be1} en prenant (avec ses notations) pour
anneau $\D :=\B _\X \otimes _{\O _{\X}} \D ^{(m)} _{\X ^\#}$. On obtient en particulier la cohérence de
son complété $p$-adique, noté $\B _\X \smash{\widehat{\otimes}} _{\O _{\X}} \widehat{\D} ^{(m)} _{\X ^\#}$
ainsi que des théorèmes de type $A$ et $B$ pour les
$\B _\X \smash{\widehat{\otimes}} _{\O _{\X}} \widehat{\D} ^{(m)} _{\X ^\#}$-modules cohérents à gauche ou à droite
(pour plus de précisions, voir \cite[3.3]{Be1}).
De même, il découle de \cite[3.4]{Be1} la cohérence de
$\B _\X \smash{\widehat{\otimes}} _{\O _{\X}} \widehat{\D} ^{(m)} _{\X ^\#}$ ainsi que des théorèmes de type
$A$ et $B$ pour les $\B _\X \smash{\widehat{\otimes}} _{\O _{\X}}\widehat{\D} ^{(m)} _{\X ^\#,\Q}$-modules cohérents à gauche ou à droite.
Le théorème suivant correspond à l'analogue logarithmique du théorème de platitude \cite[3.5.3]{Be1} de Berthelot.
\end{vide}

\begin{theo}
  \label{mm+1platQ}
  Soit $\B _\X$ une $\O _\X$-algèbre vérifiant les conditions de \ref{chap-coh-21} pour $m+1$.
  L'homomorphisme canonique
$\B _\X \smash{\widehat{\otimes}} _{\O _{\X}}\widehat{\D} ^{(m)} _{\X ^\#,\Q}
\rightarrow \B _\X \smash{\widehat{\otimes}} _{\O _{\X}}\widehat{\D} ^{(m+1)} _{\X ^\#,\Q}$
  est plat à droite et à gauche.
\end{theo}
\begin{proof}
La preuve est identique à celle de \cite[3.5.3]{Be1}.
Nous allons toutefois rappeler les points fondamentaux de la preuve de Berthelot.
On se contentera de mettre en exergue les propriétés
fondamentales des faisceaux utilisés (e.g. $\B _\X \smash{\widehat{\otimes}} _{\O _{\X}}\widehat{\D} ^{(m)} _{\X ^\#,\Q}$)
qui permettent de reprendre les calculs
de Berthelot que nous ne referons pas (pour ceux-ci, on se reportera à la preuve de \cite[3.5.3]{Be1}).

On se contente de prouver la platitude à gauche.
  L'assertion est locale. On peut donc supposer $\X$ affine et muni de coordonnées locales logarithmiques
  $t _1,\dots, t _d$. On note $D ^{(m)} :=\Gamma (\X, \B _\X \otimes _{\O _{\X}} \D ^{(m)} _{\X ^\#})$,
  $\widehat{D} ^{(m)}$ son complété $p$-adique et de même pour $m+1$. Il suffit de prouver que
  l'extension $\widehat{D} ^{(m)} _\Q \rightarrow \widehat{D} ^{(m+1)} _\Q$ est plate.
  Avec les notations de \ref{nota11} et en utilisant par exemple l'isomorphisme \ref{lemm-wildetildeDcoh2361},
  $D ^{(m)}$ est un $\Gamma (\X, \B _{\X})$-module libre de base
$\underline{\partial} _\# ^{<\underline{k}>_{(m)} }$ avec $\underline{k} \in \N ^d$,
ce qui donne aussi une description simple
de son complété $p$-adique (de même pour $m+1$).
Avec les arguments de la preuve de \cite[3.5.3]{Be1} (qui fonctionne toujours grâce à \ref{lemm-wildetildeDcoh}),
on peut en outre supposer
$\B _\X$ sans $p$-torsion.

Pour tout $\underline{k} \in \N ^d$, pour tout $i=1,\dots ,d$, posons
$k _i = p^m q ^{(m)} _{k_i} +r ^{(m)} _{k_i}=
p^{m+1} q ^{(m+1)} _{k_i} +r ^{(m+1)} _{k_i}$, avec
$0\leq r ^{(m)} _{k_i} <p^m$, $0\leq r ^{(m+1)} _{k_i} <p^{m+1}$. Alors,
$\underline{\partial} _\# ^{<\underline{k}>_{(m)} }= \frac{q _{\underline{k}} ^{(m)} !}{q _{\underline{k}} ^{(m+1)} !}
\underline{\partial} _\# ^{<\underline{k}>_{(m+1)} }$.
Ainsi, $\widehat{D} ^{(m)}$ et $D ^{(m+1)}$ sont canoniquement inclus dans
$\widehat{D} _\Q ^{(m)}$.
Notons alors
$D'$ le sous-groupe de $\widehat{D} _\Q ^{(m)}$ engendré par $\widehat{D} ^{(m)}$ et $D ^{(m+1)}$.
En reprenant les calculs de Berthelot, on obtient alors
que $D'$ est en fait un sous-anneau de $\widehat{D} _\Q ^{(m)}$, que
$\widehat{D} ^{(m+1)} = \widehat{D} ^{\prime }$ et $\widehat{D} ^{(m)} _\Q= D ^{\prime } _\Q$.
Pour terminer la preuve, il suffit donc d'établir que
$D'$ est noethérien (car cela implique que $\widehat{D}'$ est plat à droite et à gauche sur $D ' $ et donc de même avec l'indice $\Q$).

D'après \ref{lemm-wildetildeDcoh}.\ref{lemm-wildetildeDcohii},
$D ^{(m)}$ est engendré comme $\Gamma(\X,\B _\X)$-algèbre par
les opérateurs $\partial_{\#i} ^{<p^j>_{(m)} }$, où $1\leq i\leq d$, $0\leq j\leq m$,
ces derniers commutant deux à deux (de même pour $m+1$).
On en déduit que $D'$ est engendré en tant que $\widehat{D} ^{(m)}$-module à gauche par les éléments de la forme
$(\underline{\partial} _\# ^{<p^{m+1}>_{(m)} }) ^{\underline{q}}$, pour $\underline{q} \in \N ^d$ (cela a un sens via \ref{prodpartial}).
En utilisant \ref{be1224iv},
on vérifie comme Berthelot que, pour tout $r \in \N $,
pour tout $P \in \widehat{D} ^{(m)}$, on a :
\begin{equation}
  \notag
  [(\partial_{\#i} ^{<p^{m+1}>_{(m)} }) ^r,P] \in \sum _{s <r} \widehat{D} ^{(m)} (\partial_{\#i} ^{<p^{m+1}>_{(m)} }) ^s.
\end{equation}
On en déduit alors de manière identique à Berthelot que $D'$ est noethérien.
\end{proof}

\begin{theo}
\label{Be1-4.3.5}
  Soient $T$ un diviseur de $X_0$, $m' \geq m $, $r$ (resp. $r'$) un multiple de $p^{m+1}$ (resp. $p^{m'+1}$).
Avec les notations de \ref{exB(T)bis}, l'homomorphisme canonique
    $\widehat{\B} _\X  (T,r)\smash{\widehat{\otimes}} _{\O _{\X}}\widehat{\D} ^{(m)} _{\X ^\#,\Q}
\rightarrow
\widehat{\B} _\X  (T,r') \smash{\widehat{\otimes}} _{\O _{\X}}\widehat{\D} ^{(m')} _{\X ^\#,\Q}$
  est plat à droite et à gauche.
\end{theo}
\begin{proof}
Via \ref{mm+1platQ}, cela se vérifie de façon identique à la preuve de \cite[4.3.5]{Be1} (dont
les étapes clés sont les mêmes que celle de \ref{mm+1platQ}).
\end{proof}

\begin{vide}
Soit $T$ un diviseur de $X_0$.
On définit le faisceau $\D ^\dag _{\X ^\#} (\hdag T )$ des
{\og opérateurs différentiels de niveau fini sur $\X ^\#$
à singularités surconvergentes le long de $T$\fg}
en posant
$\D ^\dag _{\X ^\#} (\hdag T) := \smash{\underset{\longrightarrow}{\lim}} _m
\widehat{\B} _\X  (T,p^{m+1})\smash{\widehat{\otimes}} _{\O _{\X}}\widehat{\D} ^{(m)} _{\X ^\#} $.
Lorsque le diviseur $\ZZ$ est vide, on retrouve
$\D ^\dag _{\X } (\hdag T ) $ (voir \cite[4.2.5.3]{Be1}).
\end{vide}

\begin{theo}
\label{platdagtXdiese}
Soit $T$ un diviseur de $X_0$.
Le faisceau
$\D ^\dag _{\X ^\#} (\hdag T ) _{\Q}$ est cohérent.
On bénéficie de théorèmes de type $A$ et $B$ pour les
$\D ^\dag _{\X ^\#} (\hdag T ) _{\Q}$-modules cohérents à droite ou à gauche.
\end{theo}

\begin{proof}
  Par \cite[3.6]{Be1}, cela découle de \ref{Be1-4.3.5}.
\end{proof}

\begin{rema}
  Par contre, on ignore si
$\D ^\dag _{\X }$ (et a fortiori $\D ^\dag _{\X ^\#} (\hdag T )$) est cohérent.
\end{rema}

\begin{theo}
\label{4.3.10Be1}
  Soient $T$ un diviseur de $X_0$, $\U ^\#$ l'ouvert de $\X ^\#$ complémentaire de $T$,
  $j$ : $\U ^\# \subset \X ^\#$ l'inclusion canonique.
L'homomorphisme $\D ^\dag _{\X ^\#} (\hdag T ) _{\Q}
\rightarrow
j _* \D ^\dag _{\U ^\#,\Q} $
est fidèlement plat à droite et à gauche.
\end{theo}
\begin{proof}
  Il s'agit de reprendre la preuve de \cite[4.3.10.2]{Be2}.
\end{proof}

De manière analogue à \cite[4.3.11 et 4.3.12]{Be1},
on déduit de \ref{4.3.10Be1} la proposition ci-après.

\begin{prop}
  \label{4.3.12Be1}
  Avec les notations \ref{4.3.10Be1},
  pour qu'un $\D ^\dag _{\X ^\#} (\hdag T ) _{\Q}$-module cohérent soit nul, il faut et il suffit que sa restriction à $\U ^\#$ soit nulle.
\end{prop}

Par commodité, on s'intéressera dans un premier temps au cas des log isocristaux convergents
puis dans un second temps à celui des log isocristaux surconvergents.
Via l'équivalence de catégories de Berthelot (\cite{Be4} ou \cite[2.2.12]{caro_courbe-nouveau} pour la version publiée), la définition suivante
 correspond à celle de Kedlaya dans \cite[6.3.1]{kedlaya-semistableI}
 (d'après \cite[6.4.1]{kedlaya-semistableI}, cette notion est équivalente à celle de A. Shiho).
 Avec le théorème \ref{theologisoDdag}, nous retrouvons
 la description classique des log-isocristaux convergents en terme de log-$\D$-module arithmétique.
\begin{defi}
\label{defilogiso}
  Soit $\E$ un $\D _{\X ^\#,\Q}$-module cohérent qui soit localement projectif et de type fini sur $\O _{\X,\Q}$.
  On dit que $\E$ est un log-isocristal convergent sur $\X ^\#$
  si la structure de
  $\D _{\Y,\Q}$-module de $\E |_{\Y}$ se prolonge en une structure de $\D ^\dag _{\Y,\Q}$-module cohérent.
\end{defi}

\begin{rema}
  Si on voulait calquer \cite[6.4.1]{kedlaya-semistableI}, on aurait dû remplacer
  {\og la structure de
  $\D _{\Y,\Q}$-module de $\E |_{\Y}$ se prolonge en une structure de $\D ^\dag _{\Y,\Q}$-module cohérent\fg}
  par la condition apparemment (voir les théorèmes qui suivent) plus forte
{\og la structure de
  $\O _{\X} (\hdag Z) _{\Q} \otimes _{\O_{\X,\Q}} \D _{\X,\Q}$-module de $\O _{\X} (\hdag Z) _{\Q}  \otimes _{\O _{\X,\Q}}\E$
  se prolonge en une structure de $\D ^\dag  _{\X} (\hdag Z) _{\Q}$-module cohérent\fg}.
\end{rema}

\begin{lemm}\label{654semistableI}
  On suppose $\X ^\#$ affine et muni de coordonnées logarithmiques locales $t _1, \dots, t _d \in M (\ZZ)$.
 On pose alors $\underline{\partial} _{\#} ^{[\underline{k}]}:=
 \frac{\underline{\partial} _{\#} ^{<\underline{k}> _{(0)}}}{\underline{k}!}=\underline{t} ^{\underline{k}}
 \underline{\partial}  ^{[\underline{k}]}.$
  Soit $\E$ un log-isocristal convergent sur $\X ^\#$.
  Pour pour toute section $e \in \Gamma (\X, \E)$, tout $0\leq \eta < 1$, on a
\begin{equation}
  \label{654semistableIeq}
  \parallel \underline{\partial} _{\#} ^{[\underline{k}]} e \parallel \eta ^{|\underline{k} |} \rightarrow 0\ \text{pour}\
  |\underline{k} | \rightarrow \infty.
\end{equation}
\end{lemm}
\begin{proof}
  Via la formule \cite[Lemme 2.3.3.(c)]{these_montagnon}, cela est une réécriture de \cite[6.3.4]{kedlaya-semistableI}.
\end{proof}

\begin{prop}\label{312be0}
  Soient $\E$ un log-isocristal convergent sur $\X ^\#$ et $m\in  \N $ un entier.
  Il existe alors un $\widehat{\D} ^{(m)} _{\X ^\#}$-module $\overset{\circ}{\E}$, cohérent sur $\O _{\X}$ et
  un isomorphisme $\widehat{\D} ^{(m)} _{\X ^\#,\Q}$-linéaire
  $\overset{\circ}{\E} _\Q \riso \E$.
\end{prop}
\begin{proof}
Supposons $\X ^\#$ muni de coordonnées locales logarithmiques $t _1,\dots, t _d$.
Avec les notations \ref{nota11}, il suffit de reprendre les calculs de la preuve de \cite[3.1.2]{Be0} (ou \cite[4.4.7]{Be1})
en remplaçant \cite[3.0.1.1]{Be0} par \ref{654semistableI},
$\underline{\tau}  ^{\underline{k}}$
par $\underline{\tau} _\# ^{\underline{k}}$ et
$\underline{\partial}  ^{\underline{k}}$
(resp. $\underline{\partial} _\# ^{\underline{k}}$).
\end{proof}

\begin{prop}\label{313Be0}
  Soit $\E$ un $\D ^{(m)} _{\X ^\#}$-module, cohérent en tant que $\O _{\X}$-module.
  \begin{enumerate}
  \item \label{313iBe0} Si $\X$ est affine alors $\E$ est globalement de présentation finie sur $\D ^{(m)} _{\X ^\#}$.
    \item  \label{313iiBe0} Le faisceau $\E$ est cohérent sur $\D ^{(m)} _{\X ^\#}$.
    \item \label{313iiiBe0} L'homomorphisme canonique
\begin{equation}
  \label{313Be0eq}
  \E \rightarrow \widehat{\D} ^{(m)} _{\X ^\#} \otimes _{\D ^{(m)} _{\X ^\#}} \E
\end{equation}
est un isomorphisme.
  \end{enumerate}
\end{prop}

\begin{proof}
On vérifie \ref{313iBe0} en reprenant la preuve de \cite[3.1.3.(i)]{Be0}.
Cela implique aussitôt \ref{313iiBe0}. Traitons à présent \ref{313iiiBe0}.
Comme $\E$ est un $\O _{\X}$-module cohérent, il est canoniquement isomorphe à son complété $p$-adique.
Or, comme $\E$ est un $\D ^{(m)} _{\X ^\#}$-module cohérent (d'après ce que l'on vient de prouver),
son complété $p$-adique est canoniquement isomorphe à $\widehat{\D} ^{(m)} _{\X ^\#} \otimes _{\D ^{(m)} _{\X ^\#}} \E$.
D'où le résultat.
\end{proof}

De manière analogue à \cite[3.1.4]{Be0}, il découle de \ref{312be0} et de \ref{313Be0} la proposition suivante.
\begin{prop}\label{314be0}
  Soient $\E$ un log-isocristal convergent sur $\X ^\#$ et $m\in \N$ un entier.
    \item Les homomorphismes
    \begin{equation}
      \label{314be0eq}
      \E \rightarrow \widehat{\D} ^{(m)} _{\X ^\#,\Q} \otimes _{\D  _{\X ^\#,\Q}}\E,\hspace{1cm}
      \E \rightarrow \D ^{\dag} _{\X ^\#,\Q} \otimes _{\D  _{\X ^\#,\Q}}\E
    \end{equation}
   sont des isomorphismes.
\end{prop}

\begin{theo}\label{theologisoDdag}
Soit $\E$ un $\D _{\X ^\#,\Q}$-module cohérent localement projectif et de type fini sur $\O _{\X,\Q}$.
Le faisceau $\E$ est un log-isocristal convergent sur $\X ^\#$ si et seulement si sa structure
de $\D _{\X ^\#,\Q}$-module se prolonge (de façon unique) en une structure de
$\D ^{\dag} _{\X ^\#,\Q}$-module cohérent.
\end{theo}
\begin{proof}
  Comme $\D ^{\dag} _{\X ^\#,\Q} | \Y \riso \D ^{\dag} _{\Y,\Q}$, la condition est suffisante.
  La réciproque découle de \ref{314be0eq} et de la cohérence de $\D ^{\dag} _{\X ^\#,\Q}$.
\end{proof}

Traitons à présent le cas des log isocristaux surconvergents. Dans la suite de cette section,
$T$ sera un diviseur de $X _0$.
\begin{defi}
\label{defisurcvlog}
Un log isocristal (ou simplement isocristal) sur $\X ^\#$ surconvergent le long de $T$
  est un $\O _{\X} (\hdag T) _{\Q} \otimes _{\O _{\X,\Q} }\D _{\X ^\#,\Q}$-module cohérent $\E$,
  localement projectif et de type fini sur $\O _{\X} (\hdag T) _{\Q}$ et
  tel que
  $\E |\Y$ soit un $\D ^\dag _{\Y} (\hdag T\cap Y) _{\Q}$-module cohérent, i.e.,
  $\E |\Y$ est associé à un isocristal sur $Y_0 \setminus T$ surconvergent le long de $Y_0\cap T$.
\end{defi}

\begin{vide}
\label{dualisocsurcv}
Soit $\E$ un isocristal sur $\X ^\#$ surconvergent le long de $T$.
     Notons $\E ^\vee := \mathcal{H} om _{\O _{\X} (\hdag T) _{\Q}} (\E,\O _{\X} (\hdag T) _{\Q})$.
Comme la catégorie des isocristaux sur $Y _0\setminus T$ surconvergent le long de $Y_0 \cap T$ est stable par dualité,
on en déduit que $\E ^\vee$ est un isocristal sur $\X ^\#$ surconvergent le long de $T$.
De même, la catégorie des isocristaux sur $\X ^\#$ surconvergent le long de $T$
est stable par le bifoncteur $- \otimes _{ \O _{\X} (\hdag T) _{\Q}} -$.

\end{vide}

\begin{lemm}\label{654semistableIsurcv}
  On suppose $\X ^\#$ affine, muni de coordonnées logarithmiques locales $t _1, \dots, t _d \in M (\ZZ)$ et qu'il existe
  un relèvement $f \in \O _{\X}$ d'une équation locale de $T$ dans $X$.
  Soit $\E$ un log-isocristal sur $\X ^\#$ surconvergent le long de $T$.
Notons $\X _K$ l'espace analytique rigide de $\X$ au sens de Raynaud, $\sp $ : $ \X _K \rightarrow \X$ le morphisme de spécialisation,
$U _\lambda := \{ x \in \X _K \, | \, |f(x)| \geq \lambda\}$ et $E := \sp ^* (\E)$.
  Pour pour tout $0\leq \eta < 1$, il existe $0\leq \lambda _\eta<1$ tel que, pour tout
  $\lambda _\eta \leq \lambda < 1$ et toute section $e \in \Gamma (U _\lambda, E)$, on ait
\begin{equation}
  \label{654semistableIeqsurcv}
  \parallel \underline{\partial} _{\#} ^{[\underline{k}]} e \parallel \eta ^{|\underline{k} |} \rightarrow 0\ \text{pour}\
  |\underline{k} | \rightarrow \infty.
\end{equation}
\end{lemm}
\begin{proof}
Par hypothèse, $E | \Y _K$ est un isocristal sur $Y _0\setminus T$ surconvergent le long de $Y_0\cap T$, i.e.,
\ref{654semistableIeqsurcv}
est vrai sans dièses en remplaçant $U _\lambda $ par $U _\lambda \cap \Y _K$.
  On procède alors de manière analogue à \cite[6.3.4]{kedlaya-semistableI}.
\end{proof}

\begin{theo}
\label{isoctheodef}
  Soit $\E$ un log-isocristal sur $\X ^\#$ surconvergent le long de $T$.
  Alors $\E$ est un $\D  ^\dag _{\X ^\#} (\hdag T) _{\Q} $-module cohérent.
\end{theo}
\begin{proof}
  Via la formule \ref{654semistableIeqsurcv}, il s'agit de reprendre la preuve de \cite[4.4.12]{Be1}.
\end{proof}

\begin{rema}
\label{rema-isoctheodef}
  D'après le théorème \ref{isoctheodef}, un log-isocristal sur $\X ^\#$ surconvergent le long de $T$
  est un $\D  ^\dag _{\X ^\#} (\hdag T) _{\Q} $-module cohérent,
  localement projectif et de type fini sur $\O _{\X} (\hdag T) _{\Q}$.
\end{rema}

\begin{prop}
\label{videsurcvlog}
  Posons
  $\D  _{\X ^\#} (\hdag T) _{\Q} := \O _\X (\hdag T) _{\Q} \otimes _{\O _{\X,\Q}} \D  _{\X ^\#,\Q} $.
L'homomorphisme canonique
$      \E \rightarrow \D  ^\dag _{\X ^\#} (\hdag T) _{\Q} \otimes _{\D  _{\X ^\#} (\hdag T) _{\Q} }\E$
  est un isomorphisme.
\end{prop}

\begin{proof}
  En reprenant les arguments de la preuve de \cite[3.1.3.(i)]{Be0}
  (en effet, $\O _\X (\hdag T) _{\Q} $ est à section noethérienne
  sur les ouverts affines et on dispose de théorèmes de type $A$ et $B$
  pour les $\O _\X (\hdag T) _{\Q} $-modules cohérents), on vérifie
  que $\E$ est
  $\D  _{\X ^\#} (\hdag T) _{\Q} $-cohérent.
Ainsi,
$      \E \rightarrow \D  ^\dag _{\X ^\#} (\hdag T) _{\Q} \otimes _{\D  _{\X ^\#} (\hdag T) _{\Q} }\E$
est un homomorphisme de
$\D  ^\dag _{\X ^\#} (\hdag T) _{\Q} $-modules cohérents.
Par \ref{314be0}, cet homomorphisme est un isomorphisme en dehors de $T$.
On conclut ensuite via \ref{4.3.12Be1}.
\end{proof}

\begin{theo}
  Soit $\E$ un log-isocristal sur $\X ^\#$ surconvergent le long de $T$. Le morphisme canonique
\begin{equation}
  \label{resspencerextlogdiag4}
  \D ^\dag _{\X } (\hdag T) _{\Q} \otimes _{\D ^\dag _{\X ^\#} (\hdag T) _{\Q}} ^{\L} \E
  \rightarrow
  \D ^\dag _{\X } (\hdag T) _{\Q} \otimes _{\D ^\dag _{\X ^\#} (\hdag T) _{\Q}} \E
\end{equation}
  est un isomorphisme.
\end{theo}

\begin{proof}
Par \ref{4.3.12Be1}, il suffit de le vérifier lorsque $T$ est vide.
  Il résulte de \ref{corologdiag3q} que
  le morphisme canonique
\begin{equation}
  \label{resspencerextlogdiag3bis}
  \D  _{\X ,\Q }\otimes _{\D  _{\X ^\#,\Q }} ^{\L} \E
  \rightarrow
  \D  _{\X ,\Q }\otimes _{\D  _{\X ^\#,\Q }} \E
\end{equation}
est un isomorphisme.
Puisque les extensions
$\D  _{\X,\Q} \rightarrow \D ^{\dag} _{\X,\Q} $,
$\D  _{\X ^\#,\Q}  \rightarrow \D ^{\dag} _{\X ^\#, \Q}$
sont plates, il en résulte que le morphisme canonique
\begin{equation}
  \label{resspencerextlogdiag5}
  \D ^{\dag} _{\X,\Q} \otimes _{\D ^{\dag} _{\X ^\#,\Q}} ^{\L}
  (\D ^{\dag} _{\X ^\#,\Q} \otimes _{\D  _{\X ^\#,\Q}} \E)
  \rightarrow
  \D ^{\dag} _{\X,\Q} \otimes _{\D ^{\dag} _{\X ^\#,\Q}}
    (\D ^{\dag} _{\X ^\#,\Q} \otimes _{\D  _{\X ^\#,\Q}} \E)
\end{equation}
est un isomorphisme.
On conclut via l'isomorphisme de droite de \ref{314be0eq}.
\end{proof}

\section{Sur l'holonomie des log isocristaux surconvergents}

\begin{nota}\label{Detc0}
  On définit des $\D ^{(m)} _{X ^\#}$-modules à gauche en posant :
\begin{equation}
  \label{OXdiese}
  \O _{X} (Z):= \mathcal{H} om _{\O _{X}} ( \omega _{X} , \omega _{X ^\#}),
  \O _{X} (-Z):= \mathcal{H} om _{\O _{X}} ( \omega _{X} ^\# , \omega _{X }).
\end{equation}
En tant que $\O _X$-module, le faisceau $\O _{X} (Z)$ correspond à l'$\O _{X}$-module localement engendré
par les inverses d'une équation locale de $Z$ dans $X$, ce qui justifie la notation.
On aurait aussi pu remarquer
$\O _{X} (Z)$ est un sous-$\D ^{(m)} _{X ^\#}$-module à gauche de $j _* \O _Y$.
Via \cite[1.1.73]{Be2}, on calcule que ces deux structures de $\D ^{(m)} _{X ^\#}$-module
sur $\O _{X} (Z)$ sont identiques.

Pour tout entier $n \in \N$, on en déduit des $\D ^{(m)} _{X ^\#}$-modules à gauche en posant :
$\O _{X} (nZ):= \O _{X} (Z) ^{\otimes n}$ et $\O _{X} (-nZ):= \O _{X} (-Z) ^{\otimes n}$,
où $\otimes n$ signifie que l'on tensorise $n$-fois en tant que $\O _X$-module.
  En évaluant deux fois, on obtient (voir \ref{evalDlin}) l'isomorphisme $\D ^{(m)} _{X ^\#}$-linéaire :
  $\omega _{X ^\#} \otimes _{\O _X} \O _{X } (-Z)\otimes _{\O _X} \O _{X } (Z) \riso \omega _{X ^\#} $.
  D'où :
  $\O _{X } (-Z)\otimes _{\O _X} \O _{X } (Z) \riso \O _X $. Pour
  tous $n, n' \in \Z$, les isomorphismes canoniques
  $\O _{X } (n Z)\otimes _{\O _X} \O _{X } (n'Z) \riso \O _{X } ((n +n')Z) $
  sont donc $\D ^{(m)} _{X ^\#}$-linéaires.

Si $\E$ (resp. $\M$) est un $\D  ^{(m)} _{X^\# } $-module à gauche (resp. à droite) et $n \in \Z$,
on définit un $\D ^{(m)} _{X ^\#} $-module à gauche (resp. à droite) en posant
$\E (n Z) := \O _{X } (n Z) \otimes _{\O _{X}} \E  $
(resp. $\M (n Z) := \M  \otimes _{\O _{X}} \O _{X } (n Z)$).

D'après \ref{tranpDotimesB},
on dispose de l'isomorphe de $\D ^{(m)} _{X ^\#}$-bimodules dit de transposition
$\gamma _{\O _X (Z)}$ : $\D ^{(m)} _{X ^\#} \otimes _{\O _{X}} \O _{X } (n Z)
\riso
\O _{X } (n Z) \otimes _{\O _{X}}  \D ^{(m)} _{X ^\#}$.
Par \ref{MotimesEetcegpre}, on vérifie alors la formule
$$\gamma _{\O _X (nZ)} (\underline{\partial} ^{<\underline{k}>}_\# \otimes e)
=
\sum _{\underline{h} \leq \underline{k}}
\left \{ \begin{smallmatrix}   \underline{k} \\   \underline{h} \\ \end{smallmatrix} \right \}
\underline{\partial} ^{<\underline{k}-\underline{h}>}_\# e
\otimes
\underline{\partial} ^{<\underline{h}>}_\# .$$
Avec \cite[2.2.4.(iv)]{Be1}, il en résulte que le composé
$\D ^{(m)} _{X ^\#} \otimes _{\O _{X}} \O _{X } (n Z)
\underset{\gamma _{\O _X (nZ)}}{\riso}
\O _{X } (n Z)\otimes _{\O _{X}} \D ^{(m)} _{X ^\#}
\subset j_* \D ^{(m)} _{Y}$ est égal à l'inclusion canonique
$\D ^{(m)} _{X ^\#} \otimes _{\O _{X}} \O _{X } (n Z)
\subset j_* \D ^{(m)} _{Y}$.
On notera alors sans ambiguïté
$\D ^{(m)} _{X ^\#} (n Z)$ pour
$\D ^{(m)} _{X ^\#} \otimes _{\O _{X}} \O _{X } (n Z)$ ou
$\O _{X } (n Z)\otimes _{\O _{X}} \D ^{(m)} _{X ^\#}$.
On remarque que ni $\O _{X } (n Z)$ ni $\D ^{(m)} _{X ^\#} (n Z)$ sont des faisceaux d'anneaux.
Par \ref{HomD}, on bénéficie de l'isomorphisme canonique :
$\E (n Z) \riso \D ^{(m)} _{X ^\#} (n Z) \otimes _{\D ^{(m)} _{X ^\#}} \E$
et
$\M (n Z) \riso \M\otimes _{\D ^{(m)} _{X ^\#}} \D ^{(m)} _{X ^\#} (n Z) $.
Il en découle les isomorphismes :
\begin{equation}
  \label{Zdroiteàgauche}
\M (n Z)  \otimes _{\D ^{(m)} _{X ^\#}} \E
\riso
\M \otimes _{\D ^{(m)} _{X ^\#}} \D ^{(m)} _{X ^\#}(n Z)  \otimes _{\D ^{(m)} _{X ^\#}}  \E
\riso
\M \otimes _{\D ^{(m)} _{X ^\#}} \E (n Z).
\end{equation}
\end{nota}

\begin{lemm}\label{lemmomeoomegadiese}
Soient $\E$ un $\D ^{(m)}  _{X^\# } $-module à gauche et
$\M$ un $\D ^{(m)} _{X^\#}$-module à droite.
On dispose des isomorphismes canoniques $\D ^{(m)} _{X ^\#}$-linéaires suivants :
\begin{gather}
\label{lemmomeoomegadiese-iso1}
  \mathrm{ev}\ :\ \omega _{X} \otimes _{\O _{X}}  \O _{X} (Z) \riso \omega _{X ^\#}, \
  \mathrm{ev}\ :\ \omega _{X ^\#} \otimes _{\O _{X}}  \O _{X} (-Z) \riso \omega _{X },
  \\
  \label{lemmomeoomegadiese-iso2}
  \mathrm{ev}\otimes \mathrm{Id}\ :\
  \omega _{X } \otimes _{\O _{X}} \E (Z )\riso
    \omega _{X ^\#} \otimes _{\O _{X}} \E, \
    \mathrm{ev}\otimes \mathrm{Id}\ :\
  \omega _{X^\# } \otimes _{\O _{X}} \E (-Z )\riso
    \omega _{X } \otimes _{\O _{X}} \E,
  \\
  \label{lemmomeoomegadiese-iso3}
    \E (-Z)
  \riso
  (\omega _{X }\otimes _{\O _{X}} \E) \otimes _{\O _{X}} \omega ^{-1} _{X ^\#}, \
  \E(Z)
  \riso
  (\omega _{X ^\#}\otimes _{\O _{X}} \E) \otimes _{\O _{X}} \omega ^{-1} _{X },
    \\
  \label{lemmomeoomegadiese-iso4}
  \M (Z)
  \riso
  \omega _{X ^\#}\otimes _{\O _{X}} (\M \otimes _{\O _{X}} \omega ^{-1} _{X}), \
  \M (-Z)
  \riso
  \omega _{X }\otimes _{\O _{X}} (\M \otimes _{\O _{X}} \omega ^{-1} _{X ^\#}).
\end{gather}
\end{lemm}

\begin{proof}
La $\D ^{(m)} _{X ^\#}$-linéarité de \ref{lemmomeoomegadiese-iso1} découle de \ref{evalDlin}.
Par \ref{comm-assoc}, il en dérive \ref{lemmomeoomegadiese-iso2}.
Via \ref{dr-ga}, il en résulte les autres isomorphismes de $\D ^{(m)} _{X ^\#}$-modules.
\end{proof}

\begin{nota}
On désigne par $\widetilde{\D} _{\X ^\#}$ l'un des faisceaux d'anneaux
$\D ^{(0)} _{\X ^\#}$, $\D _{\X ^\#,\Q}$,
$\widehat{\D} ^{(m)} _{\X ^\#}$,
$\widehat{\D} ^{(m)} _{\X ^\#,\Q}$,
$\D ^\dag _{\X ^\# ,\Q}$.
De même en enlevant les dièses.
\end{nota}

\begin{vide}
\label{transdaghat}
On dispose pour tout entier $n $ de l'isomorphisme canonique
$\D _{\X ^\#}$-bimodules dit de transposition
$\gamma _{\O _\X (\ZZ)}$ : $\widetilde{\D} _{\X ^\#} \otimes _{\O _{\X}} \O _{\X } (n \ZZ)
\riso
\O _{\X } (n \ZZ) \otimes _{\O _{\X}}  \widetilde{\D} _{\X ^\#}$.
En effet, on le sait déjà lorsque $\widetilde{\D} _{\X ^\#} =\D ^{(0)} _{\X ^\#}$.
Les autres cas s'en déduisent par tensorisation par $\Q$ sur $\Z$, complétion $p$-adique et
passage à la limite inductive sur le niveau.
Soient $\E$ (resp. $\M$) un $\widetilde{\D} _{\X ^\#}$-module à gauche (resp. à droite).
Avec les notations analogues à \ref{Detc0},
en reprenant la construction de \ref{Zdroiteàgauche}, on obtient
l'isomorphisme canonique fonctoriel $\E$ et $\M$ :
\begin{equation}
\label{Zdroiteàgauchebis}
\M (n \ZZ)  \otimes _{\widetilde{\D} _{\X ^\#}} \E
\riso
\M \otimes _{\widetilde{\D} _{\X ^\#}} \E (n \ZZ).
\end{equation}

De même,
 le lemme \ref{lemmomeoomegadiese} s'étend naturellement en remplaçant
{\og $\D ^{(m)}  _{X^\# } $\fg} par {\og $\widetilde{\D} _{\X ^\# }$\fg}.
Les références relatives à \ref{lemmomeoomegadiese} pourront abusivement concerner ces extensions.
\end{vide}

\begin{vide}
\label{dimcohofinie}
Par \ref{resspencerlog}, on vérifie avec les arguments habituels (e.g. \cite[4.4.3]{Be2})
que $\D ^{(0)} _{X _0 ^\#}$, $\D ^{(0)} _{\X ^\#}$ et $\D _{\X ^\#,\Q}$ sont de dimension homologique finie.
De plus, en reprenant le début de la preuve de \cite[4.4.4]{Be2} et en y remplaçant
\cite[4.4.3]{Be2} par \cite[5.3.1]{these_montagnon}, on vérifie que si
$\X$ est affine alors l'anneau $\Gamma (\X,\widehat{\D} ^{(m)} _{\X ^\#})$ est de dimension homologique finie.
Il en résulte que
$\Gamma (\X,\widehat{\D} ^{(m)} _{\X ^\#,\Q})$
et
$\Gamma (\X,\D ^\dag _{\X ^\#,\Q})$ sont de dimension homologique finie
lorsque $\X$ est affine (car un $\Gamma (\X,\D ^{\dag} _{\X ^\#,\Q})$-module cohérent provient par extension
d'un $\Gamma (\X,\widehat{\D} ^{(m)} _{\X ^\#,\Q})$-module cohérent qui lui provient
d'un $\Gamma (\X,\widehat{\D} ^{(m)} _{\X ^\#})$-module cohérent, de plus ces extensions sont plates).
Via les théorèmes de type $A$, il en résulte que les faisceaux
$\widehat{\D} ^{(m)} _{\X ^\#}$, $\widehat{\D} ^{(m)} _{\X ^\#,\Q}$,
$\D ^\dag _{\X ^\#,\Q}$ sont de dimension homologique finie.
Comme $\widetilde{\D} _{\X^\#}$ est en outre cohérent, on obtient ainsi
$D ^\mathrm{b} _{\mathrm{coh}} (\overset{^*}{} \widetilde{\D}  _{\X ^\#}) =
D _{\mathrm{parf}} (\overset{^*}{} \widetilde{\D}  _{\X ^\#}) $.

Lorsque $i\geq 1$, comme $S _i$ n'est pas régulier, les faisceaux
$\D ^{(m)} _{X _i ^\#}$ ne sont pas de dimension homologique finie.
J'ignore ce qu'il en est lorsque $m\not =0$ de $\D ^{(m)} _{\X ^\#}$
et, lorsque $\ZZ$ et $T$ sont non vides,
de $\D ^\dag _{\X ^\#} (\hdag T) _{\Q}$ (lorsque $\ZZ$
est vide, c'est bien le cas d'après \cite{huyghe_finitude_coho}).
\end{vide}

\begin{defi}
\label{defidual}
  Soient
  $\E \in D ^\mathrm{b} _{\mathrm{coh}} (\overset{^g}{} \widetilde{\D}  _{\X ^\#}) $,
$\M  \in D ^\mathrm{b} _{\mathrm{coh}} ( \widetilde{\D}  _{\X ^\#} \overset{^d}{}  ) $.
  On définit
  {\it les duaux $\widetilde{\D}  _{\X ^\#}$-linéaires} de $\E$ et de $\M$ en posant
\begin{gather}
  \DD _{\X ^\#} (\E) =
\R \mathcal{H} om _{ \widetilde{\D}  _{\X ^\#}}
(\E, \widetilde{\D}  _{\X ^\#} ) \otimes _{\O _{\X}} \omega ^{-1} _{\X ^\#} [d _X],
\\
\DD _{\X ^\#} (\M) =
\omega  _{\X ^\#} \otimes _{\O _{\X}}  \R \mathcal{H} om _{ \widetilde{\D}  _{\X ^\#}}
(\M, \widetilde{\D}  _{\X ^\#} )  [d _X].
\end{gather}

Par \ref{dimcohofinie},
ces foncteurs duaux stabilisent donc
$D ^\mathrm{b} _{\mathrm{coh}} (\overset{^*}{}  \widetilde{\D}  _{\X ^\#}) $.
De plus, on vérifie comme Virrion (voir \cite{virrion}) l'isomorphisme de bidualité
$ \DD _{\X ^\#} \circ \DD _{\X ^\#} (\E) \riso \E$ (de même pour $\M$).
Nous verrons cependant que l'isomorphisme de dualité relative est mise en défaut pour
le morphisme canonique $\X ^\# \rightarrow \X$ (voir \ref{dualrelafaux}).

\end{defi}

\begin{vide}
\label{videcoroequiomeg}
  Il résulte des équivalences de catégories \ref{dr-ga} que, pour tous
  $\E \in D  (\overset{^g}{} \widetilde{\D}  _{\X ^\#}) $,
  $\M \in D ^+ ( \widetilde{\D}  _{\X ^\#}\overset{^d}{}) $,
  on a
  \begin{gather}
    \label{coroequiomeg}
    \R \mathcal{H} om _{\widetilde{\D}  _{\X ^\#}} ( \omega _{\X ^\#} \otimes _{\O _\X} \E , \M)
    \riso
    \R \mathcal{H} om _{\widetilde{\D}  _{\X ^\#}} ( \E , \M \otimes _{\O _\X} \omega _{\X ^\#} ^{-1}),\\
    \label{coroequiomegbis}
    \R \mathcal{H} om _{\widetilde{\D}  _{\X ^\#}} ( \M \otimes _{\O _\X} \omega _{\X ^\#} ^{-1}, \E)
    \riso
    \R \mathcal{H} om _{\widetilde{\D}  _{\X ^\#}} ( \M , \omega _{\X ^\#} \otimes _{\O _\X} \E).
  \end{gather}
D'où :
$\DD _{\X ^\#} (\omega _{\X ^\#} \otimes _{\O _\X} \E ) \riso
\omega _{\X ^\#} \otimes _{\O _\X}  \DD _{\X ^\#} (\E )$
et
$\DD _{\X ^\#} (\M \otimes _{\O _\X} \omega _{\X ^\#} ^{-1})
\riso
\DD _{\X ^\#} (\M )\otimes _{\O _\X} \omega _{\X ^\#} ^{-1}$.
On dispose de même des isomorphismes \ref{coroequiomeg} et \ref{coroequiomegbis} où
{\og $\omega _{\X ^\#}$\fg} est remplacé par {\og $\omega _{\X }$\fg}.
\end{vide}

\begin{nota}
  On note $u$ : $\X ^\# \rightarrow \X$ le morphisme canonique et
$\widetilde{\D}   _{\X \leftarrow \X ^\#} :=\widetilde{\D}  _{\X} \otimes _{\O _{\X}} \O _{\X } (\ZZ)$
vu comme $(\widetilde{\D}  _{\X} , \widetilde{\D} _{\X ^\#} )$-bimodule.
Par \ref{lemmomeoomegadiese-iso4},
$\widetilde{\D}   _{\X \leftarrow \X ^\#}$ est canoniquement isomorphe
à
$\omega _{\X ^\#} \overset{\mathrm{d}}{\otimes} _{\O _{\X}} (\widetilde{\D}  _{\X}  \otimes _{\O _{\X}} \omega ^{-1} _{\X })$, ce qui justifie la notation.
On pourrait aussi désigner par
$\widetilde{\D}   _{\X ^\# \rightarrow \X } :=\widetilde{\D}  _{\X} $
vu comme $(\widetilde{\D} _{\X ^\#}, \widetilde{\D}  _{\X}  )$-bimodule, mais cela ne sert qu'à alourdir les notations.
\end{nota}

\begin{defi}
Pour tous
$\E \in D ^\mathrm{b} _{\mathrm{coh}} (\overset{^g}{} \widetilde{\D}  _{\X ^\#}) $,
$\M  \in D ^\mathrm{b} _{\mathrm{coh}} ( \widetilde{\D} _{\X ^\# } \overset{^d}{}  ) $,
on définit respectivement l'{\it image directe} par $u $ de $\E$ et $\M$ en posant :
\begin{gather}
  \label{defimdirdag}
    u ^\mathrm{g} _+ (\E)  :=
    \widetilde{\D}   _{\X \leftarrow \X ^\#} \otimes ^\L _{ \widetilde{\D}  _{\X ^\# }} \E,
    \hspace{1cm}
    u ^\mathrm{d} _+ (\M)  :=
    \M \otimes ^\L _{ \widetilde{\D}  _{\X ^\#}}\widetilde{\D}   _{\X} .
\end{gather}
Si aucune confusion n'est à craindre, on écrira $u _+$ pour $u ^\mathrm{g} _+ $
ou $u ^\mathrm{d} _+$.
\end{defi}

\begin{lemm}\label{D=wDw-1}
  Soient $\E \in D ^- (\overset{^g}{} \widetilde{\D}  _{\X ^\#}) $,
$\M \in D ^- ( \widetilde{\D} _{\X ^\# } \overset{^d}{}  ) $.
On dispose d'un isomorphisme canonique :
$$ \M \otimes _{\widetilde{\D} _{\X ^\# } } ^\L \E
\riso
(\omega _{\X ^\#} \otimes _{\O _{\X}} \E )
\otimes _{\widetilde{\D} _{\X ^\# } } ^\L
(\M\otimes _{\O _{\X}} \omega _{\X ^\#} ^{-1} ).$$
\end{lemm}
\begin{proof}
Avec \ref{dr-ga}, on vérifie que l'on est bien dans le contexte de \cite[I.2.2]{virrion}.
\end{proof}

\begin{prop}
\label{u+dcommg}
Pour tous
$\E \in D ^\mathrm{b} _{\mathrm{coh}} (\overset{^g}{} \widetilde{\D}  _{\X ^\#}) $,
$\M \in D ^\mathrm{b} _{\mathrm{coh}} ( \widetilde{\D} _{\X ^\# } \overset{^d}{}  ) $,
on bénéficie des isomorphismes canoniques :
\begin{gather}
  \label{u+dcommgeq1}
  u _+ ^\mathrm{d} ( \M ) \otimes _{\O _{\X}} \omega _{\X} ^{-1}
  \riso
  u _+ ^\mathrm{g} ( \M \otimes _{\O _{\X}} \omega _{\X ^\# } ^{-1}),
  \hspace{1cm}
   \omega _{\X } \otimes _{\O _{\X}} u _+ ^\mathrm{g} ( \E)
  \riso
  u _+ ^\mathrm{d} ( \omega _{\X ^\#} \otimes _{\O _{\X}} \E).
\end{gather}
De plus, $u _+ ^\mathrm{g} ( \E)\in D ^\mathrm{b} _{\mathrm{coh}} (\overset{^g}{} \widetilde{\D}  _{\X ^\#}) $ et
$u _+ ^\mathrm{d} ( \M ) \in D ^\mathrm{b} _{\mathrm{coh}} ( \widetilde{\D} _{\X ^\# } \overset{^d}{}  ) $.
\end{prop}

\begin{proof}
  On construit l'isomorphisme de gauche de \ref{u+dcommgeq1} comme suit :
\begin{gather}
\notag  u _+ ^\mathrm{d} ( \M ) \otimes _{\O _{\X}} \omega _{\X } ^{-1}
  \underset{\ref{D=wDw-1}}{\riso}
(\omega _{\X ^\# }
\smash{\overset{\mathrm{g}}{\otimes}} _{\O _{\X}}
( \widetilde{\D} _{\X} \otimes _{\O _{\X}} \omega _{\X  } ^{-1}))
\otimes ^\L _{ \widetilde{\D}  _{\X ^\#}} (\M \otimes _{\O _{\X}} \omega _{\X ^\# } ^{-1})
\\ \notag
  \underset{\delta}{\riso}
(\omega _{\X ^\# }
\smash{\overset{\mathrm{d}}{\otimes}} _{\O _{\X}}
( \widetilde{\D} _{\X} \otimes _{\O _{\X}} \omega _{\X  } ^{-1}))
\otimes ^\L _{ \widetilde{\D}  _{\X ^\#}}
(\M \otimes _{\O _{\X}} \omega _{\X ^\# } ^{-1})
\underset{\ref{lemmomeoomegadiese-iso4}}{\riso}
  u _+ ^\mathrm{g} ( \M \otimes _{\O _{\X}} \omega _{\X ^\# } ^{-1}),
\end{gather}
les symboles {\og g \fg} et {\og d \fg} signifiant respectivement
que pour calculer le produit tensoriel on choisit la structure gauche et droite
de $\widetilde{\D}  _{\X}$-module à gauche de $\widetilde{\D}  _{\X}\otimes _{\O _{\X}} \omega ^{-1} _{\X}$.
On en déduit par passage de gauche à droite (i.e., via les équivalences de catégories de \ref{dr-ga})
l'isomorphisme de droite de \ref{u+dcommgeq1}.
Concernant la dernière assertion, il s'agit d'établir la préservation de la perfection (voir \ref{dimcohofinie}).
Le cas des modules à droite est immédiat.
Comme les structures tordues préservent l'exactitude et
la projectivité locale de type fini (cela découle de \ref{dr-gabis}),
elles préservent aussi les complexes parfaits.
Le cas à gauche résulte alors de \ref{u+dcommgeq1}.
\end{proof}

\begin{defi}
Grâce à \ref{u+dcommg},
pour tout $\G \in D ^\mathrm{b} _{\mathrm{coh}} (\overset{^*}{} \widetilde{\D}  _{\X ^\#}) $,
l'{\it image directe extraordinaire} par $u $ de $\G$ est définie en posant :
\begin{gather}
  \label{defimdirdagext}
    u  _! (\G)  := \DD   _{\X }  \circ  u  _+ \circ  \DD   _{\X ^\#} (\G)\in D ^\mathrm{b} _{\mathrm{coh}} (\overset{^*}{} \widetilde{\D}  _{\X ^\#}) .
\end{gather}
Pour préciser qu'il s'agit de module à gauche ou à droite,
on écrira $u ^\mathrm{g}  _!$ ou $u ^\mathrm{d}  _!$ pour $u _!$.
\end{defi}

\begin{prop}
\label{u+dcommgbis}
Pour tous
$\E \in D ^\mathrm{b} _{\mathrm{coh}} (\overset{^g}{} \widetilde{\D}  _{\X ^\#}) $,
$\M \in D ^\mathrm{b} _{\mathrm{coh}} ( \widetilde{\D} _{\X ^\# } \overset{^d}{}  ) $,
on bénéficie des isomorphismes canoniques :
\begin{gather}
  u _! ^\mathrm{g} ( \M \otimes _{\O _{\X}} \omega _{\X ^\# } ^{-1})
  \riso
  u _! ^\mathrm{d} ( \M ) \otimes _{\O _{\X}} \omega _{\X} ^{-1},
    \label{u+dcommgeq2}
  \hspace{1cm}
  \omega _{\X } \otimes _{\O _{\X}} u _! ^\mathrm{g} ( \E)
  \riso
  u _! ^\mathrm{d} ( \omega _{\X ^\#} \otimes _{\O _{\X}} \E).
\end{gather}
\end{prop}

\begin{proof}
Cela provient par composition de \ref{videcoroequiomeg} et \ref{u+dcommg}.
\end{proof}

\begin{prop}\label{dualrellog}
Pour tous
$\E \in D ^\mathrm{b} _{\mathrm{coh}} (\overset{^g}{} \widetilde{\D}  _{\X ^\#}) $,
$\M \in D ^\mathrm{b} _{\mathrm{coh}} ( \widetilde{\D} _{\X ^\# } \overset{^d}{}  ) $,
  on dispose des isomorphismes canoniques :
\begin{gather}
  \label{dualrellogeq1}
  u _!  (\E) \riso  \widetilde{\D} _{\X } \otimes ^\L _{\widetilde{\D}  _{\X ^\# }} \E,\\
    \label{dualrellogeq2}
u _!  (\M) \riso \M  \otimes ^\L _{ \widetilde{\D}  _{\X ^\# }} (\O _{\X} (-\ZZ) \otimes _{\O _{\X}} \widetilde{\D}   _{\X }) .
\end{gather}
\end{prop}

\begin{proof}
Par définition,
\begin{equation}
  \label{dualrellogiso1}
\DD _{\X} \circ u _+ (\E) \riso
\R \mathcal{H} om _{ \widetilde{\D}   _{\X }}  (
(\widetilde{\D}  _{\X} \otimes _{\O _{\X}} \O _{\X } (\ZZ))\otimes ^\L _{ \widetilde{\D}  _{\X ^\# }} \E,
\widetilde{\D}   _{\X } \otimes _{\O _{\X}} \omega ^{-1} _{\X} ) [d ].
\end{equation}
En notant
$\delta $ l'isomorphisme de transposition
$\widetilde{\D}  _{\X}\otimes _{\O _{\X}} \omega ^{-1} _{\X} \riso
\widetilde{\D}  _{\X}\otimes _{\O _{\X}} \omega ^{-1} _{\X}$ (voir \cite[1.3]{Be2}),
on obtient les isomorphismes de $(\widetilde{\D}  _{\X} , \widetilde{\D} _{\X ^\#} )$-bimodules :
\begin{equation}
  \widetilde{\D}  _{\X} \otimes _{\O _{\X}} \O _{\X} (\ZZ)
  \underset{\ref{lemmomeoomegadiese-iso4}}{\riso}
  \omega _{\X ^\#}\otimes ^{\mathrm{d}} _{\O _{\X}} (\widetilde{\D}  _{\X}\otimes _{\O _{\X}} \omega ^{-1} _{\X})
  \underset{\delta}{\riso}
  (\omega _{\X ^\#}\otimes  _{\O _{\X}} \widetilde{\D}  _{\X} ) \otimes _{\O _{\X}} \omega ^{-1} _{\X},
\end{equation}
le symbole {\og d \fg} signifiant que pour calculer le produit tensoriel on choisit la structure droite
de $\widetilde{\D}  _{\X}$-module à gauche de $\widetilde{\D}  _{\X}\otimes _{\O _{\X}} \omega ^{-1} _{\X}$.
Avec \ref{dualrellogiso1} et la version sans dièse de \ref{coroequiomegbis},
on en déduit le premier isomorphisme :
\begin{gather}
\label{dualrellogiso2pre}
\notag
\DD _{\X} \circ u _+ (\E)
\riso
\R \mathcal{H} om _{ \widetilde{\D}   _{\X }}  (
(\omega _{\X ^\#} \otimes _{\O _{\X}} \widetilde{\D}  _{\X}  )
\smash{\overset{\mathrm{g}}{\otimes}} ^{\L} _{ \widetilde{\D}  _{\X ^\# }} \E,
\omega _{\X } \otimes _{\O _{\X}} \widetilde{\D}   _{\X } \otimes _{\O _{\X}} \omega ^{-1} _{\X} ) [d ]
\\
\notag \label{dualrellogiso2}
\underset{\delta}{\riso}
\R \mathcal{H} om _{ \widetilde{\D}   _{\X }}  (
(\omega _{\X ^\#} \otimes _{\O _{\X}} \widetilde{\D}  _{\X}  )
\smash{\overset{\mathrm{g}}{\otimes}} ^{\L} _{ \widetilde{\D}  _{\X ^\# }}\E , \widetilde{\D}   _{\X } ) [d ]
\\
\notag \label{dualrellogiso3}
\riso
\R \mathcal{H} om _{ \widetilde{\D}   _{\X ^\# }}  (
(\omega _{\X ^\#} \otimes _{\O _{\X}} \widetilde{\D}  _{\X ^\#}  )
\smash{\overset{\mathrm{g}}{\otimes}} ^{\L} _{ \widetilde{\D}  _{\X ^\# }} \E , \widetilde{\D}   _{\X } )  [d ]
\\
\notag \label{dualrellogiso4}
\underset{\delta _\#}{\riso}
\R \mathcal{H} om _{ \widetilde{\D}   _{\X ^\# }}
(\omega _{\X ^\#} \otimes _{\O _{\X}} \E , \widetilde{\D}   _{\X } ) [d ]
\riso
\widetilde{\D}   _{\X } \otimes ^{\L} _{\widetilde{\D}   _{\X ^\# }}
\R \mathcal{H} om _{ \widetilde{\D}   _{\X ^\# }}
(\omega _{\X ^\#} \otimes _{\O _{\X}} \E , \widetilde{\D}   _{\X ^\# } ) [d ]
\\
\notag
\underset{\delta \circ \ref{coroequiomeg}}{\riso}
\widetilde{\D}   _{\X } \otimes ^{\L} _{\widetilde{\D}   _{\X ^\# }}
\DD _{\X ^\#} (\E ).
\end{gather}
D'où \ref{dualrellogeq1} par bidualité.
Afin d'établir \ref{dualrellogeq2}, on peut supposer
$\M =\omega _{\X ^\#} \otimes _{\O _{\X}} \E $.
On obtient :
\begin{gather}
  u _! (\M) \underset{\ref{u+dcommgeq2}}{\riso} \omega _{\X } \otimes _{\O _{\X}} u _! (\E )
\underset{\ref{dualrellogeq1}}{\riso}
(\omega _{\X } \otimes _{\O _{\X}} \widetilde{\D}  _{\X } )\otimes ^\L _{\widetilde{\D}  _{\X ^\# }} \E
\underset{\delta}{\riso}
(\omega _{\X } \otimes _{\O _{\X}} \widetilde{\D}  _{\X } )
\smash{\overset{\mathrm{g}}{\otimes}} ^\L _{\widetilde{\D}  _{\X ^\# }} \E
\\
\underset{\ref{D=wDw-1}}{\riso}
(\omega _{\X ^\#} \otimes _{\O _{\X}} \E)
\otimes ^\L _{\widetilde{\D}  _{\X ^\# }}
((\omega _{\X } \otimes _{\O _{\X}} \widetilde{\D}  _{\X } ) \smash{\overset{\mathrm{g}}{\otimes}} _{\O _{\X}} \omega _{\X ^\#} ^{-1})
\underset{\ref{lemmomeoomegadiese-iso3}}{\riso}
\M
\otimes ^\L _{\widetilde{\D}  _{\X ^\# }} (\O _{\X } (-\ZZ) \otimes _{\O _{\X}} \widetilde{\D} _{\X } ) .
\end{gather}
\end{proof}

\begin{rema}\label{dualrelafaux}
  La proposition \ref{dualrellog} illustre le fait que l'isomorphisme de dualité relative
  est inexact
pour un morphisme de log-$\V$-schémas formels lisses
dont le morphisme induit (en oubliant les monoïdes) de $\V$-schémas formels lisses est propre.
\end{rema}

\begin{lemm}
  \label{extf+log}
  Soient ${\E} \in D ^\mathrm{b} _{\mathrm{coh}} (\overset{^g}{} \widetilde{\D}  _{\X ^\#}) $ et
  ${\M} \in D ^\mathrm{b} _{\mathrm{coh}} (\widetilde{\D}  _{\X ^\#} \overset{^d}{} ) $.
  On dispose des isomorphismes canoniques :
  \begin{gather}
    \label{extf+logeq}
u _+  ( {\M} ) \otimes  _{\widetilde{\D}  _{\X }} \D ^\dag _{\X  ,\Q}
\riso
u _+  ( {\M} \otimes  _{\widetilde{\D}  _{\X ^\#}} \D ^\dag _{\X ^\# ,\Q} ),
\\
    \label{extf+logeq2}
\D ^\dag _{\X  ,\Q} \otimes  _{\widetilde{\D}  _{\X }} u _+  ( {\E} )
\riso
u _+  (\D ^\dag _{\X ^\# ,\Q}  \otimes  _{\widetilde{\D}  _{\X ^\#}} {\E} ).
  \end{gather}
  De même en remplaçant l'image directe par l'image directe extraordinaire.
\end{lemm}
\begin{proof}
Comme le foncteur dual commute à l'extension des scalaires (voir par example \cite{virrion}),
il suffit de traiter le cas de l'image directe.
L'isomorphisme \ref{extf+logeq} est immédiat tandis que \ref{extf+logeq2} se construit comme suit :
\begin{gather}
  \notag
\D ^\dag _{\X  ,\Q} \otimes  _{\widetilde{\D}  _{\X }} u _+  ( {\E} )
=
\D ^\dag _{\X  ,\Q} \otimes  _{\widetilde{\D}  _{\X }}
((\widetilde{\D}  _{\X} \otimes _{\O _{\X}} \O _{\X } (\ZZ)) \otimes ^\L _{ \widetilde{\D}  _{\X ^\# }} \E)
\underset{\ref{HomD2}}{\riso}
(\D ^\dag _{\X  ,\Q} \otimes _{\O _{\X}} \O _{\X } (\ZZ) ) \otimes ^\L _{ \widetilde{\D}  _{\X ^\# }} \E
\riso
\\
\notag
\riso
(\D ^\dag _{\X  ,\Q} \otimes _{\O _{\X}} \O _{\X } (\ZZ) )
\otimes ^\L _{\D ^\dag _{\X ^\# ,\Q} }
\D ^\dag _{\X ^\# ,\Q}
\otimes ^\L _{ \widetilde{\D}  _{\X ^\# }} \E
=
u _+  (\D ^\dag _{\X ^\# ,\Q}  \otimes  _{\widetilde{\D}  _{\X ^\#}} {\E} ).
\end{gather}

\end{proof}

\begin{prop}
\label{u+u!}
  On dispose,
  pour tous
$\E \in D ^\mathrm{b} _{\mathrm{coh}} (\overset{^*}{} \widetilde{\D}  _{\X ^\#}) $, de l'isomorphisme canonique :
\begin{equation}
  \label{u+u!1}
  u _+  (\E) \riso u _!  (\E (\ZZ)).
\end{equation}
\end{prop}

\begin{proof}
D'après \cite[3.6.2.(ii)]{Be1} (resp. \cite[3.4.5]{Be1}),
un $\D ^\dag _{\X ^\# ,\Q}$-module cohérent (resp. $\widehat{\D} ^{(m)} _{\X ^\#,\Q}$-module cohérent)
provient par extension
d'un $\widehat{\D} ^{(m)} _{\X ^\#,\Q}$-module cohérent
(resp. $\widehat{\D} ^{(m)} _{\X ^\#}$-module cohérent).
Par \ref{extf+log}, il suffit alors de traiter le cas où $\widetilde{\D} = \D ^{(0)}$
ou celui où $\widetilde{\D} = \widehat{\D} ^{(m)}$. La preuve du second cas étant la même
(on remplace \ref{omegaEDF} par \ref{omegaEDFhat}), contentons-nous d'étudier le premier.
Avec \ref{u+dcommgeq1} et \ref{u+dcommgeq2}, il suffit de traiter le cas à gauche (i.e. $* =g$).
  Soient $\PP$ une résolution de $\D ^{(0)}  _{\X }$ par des $\D ^{(0)}  _{\X }$-bimodules plats
et $\PP ^\#$ une résolution de $\E$ par des $\D ^{(0)}  _{\X ^\#}$-modules à gauche plats.
On dispose alors des isomorphismes :
\begin{gather}
\notag
  (\omega _\X \otimes _{ \O _{\X}} \E) \otimes ^\L _{\D ^{(0)}  _{\X ^\#}} \D ^{(0)}  _{\X }
  \liso
  (\omega _\X \otimes _{ \O _{\X}} \PP ^\#) \otimes _{\D ^{(0)}  _{\X ^\#}} \PP
  \overset{\ref{omegaEDF}}{\riso}
(\omega _\X \otimes _{ \O _{\X}} \PP)
\smash{\overset{\mathrm{g}}{\otimes}}  _{ \D ^{(0)}  _{\X ^\# }}   \PP ^\#
\\ \notag
\riso
(\omega _\X \otimes _{ \O _{\X}} \D ^{(0)}  _{\X })
\smash{\overset{\mathrm{g}}{\otimes}}  _{ \D ^{(0)}  _{\X ^\# }}   \PP ^\#
\overset{\delta}{\riso}
\omega _\X \otimes _{ \O _{\X}} (\D ^{(0)}  _{\X } \otimes _{\D ^{(0)}  _{\X ^\#}} \PP ^\#)
\riso
\omega _\X \otimes _{ \O _{\X}} (\D ^{(0)}  _{\X } \otimes ^\L _{\D ^{(0)}  _{\X ^\#}} \E).
\end{gather}
Ainsi, $u _+ ( \omega _\X \otimes _{\O _\X } \E) \riso \omega _\X \otimes _{\O _\X} u _! (\E)$.
En lui appliquant
$-\otimes _{\O _\X } \omega _\X  ^{-1}$,
il en résulte, via \ref{lemmomeoomegadiese}, \ref{u+dcommgeq1},
l'isomorphisme :
$  u _+  (\E (-\ZZ)) \riso u _!  (\E )$.

\end{proof}

\begin{nota}
\label{notaDDT}
Soit $T$ un diviseur de $X_0$.
    De manière analogue à \ref{defidual}, on définit le dual
  $\D  _{\X ^\#} (\hdag T) _{\Q}$-linéaire (resp. $\D ^\dag  _{\X ^\#} (\hdag T) _{\Q}$-linéaire)
  des complexes parfaits de $\D  _{\X ^\#} (\hdag T) _{\Q}$-modules
  (resp. $\D ^\dag  _{\X ^\#} (\hdag T) _{\Q}$-modules)
  que l'on notera $\DD _{\X ^\#,T}$  (resp. $\DD ^\dag _{\X ^\#,T}$ ou si aucune confusion
  n'est à craindre $\DD _{\X ^\#,T}$).

Remarquons que comme on ne sait pas si $\D ^\dag  _{\X ^\#} (\hdag T) _{\Q}$ est de dimension homologique finie,
pour utiliser les isomorphismes standards concernant les faisceaux d'homomorphismes,
il faut manipuler les complexes parfaits à la place des complexes à cohomologie bornée et cohérente.
Par exemple, pour obtenir l'isomorphisme \ref{dualrellogeq1}, nous avons utilisé la perfection
de $\omega _{\X ^\#} \otimes _{\O _{\X}} \E $. La proposition suivante donne un exemple de
tels complexes.
\end{nota}

\begin{prop}
\label{isocparf}
  Soient $T$ un diviseur de $X_0$, $\E $ un log isocristal sur $\X ^\#$ surconvergent le long de $T$.
  Alors $\E \in D  _{\mathrm{parf}} (\overset{^g}{} \D  _{\X ^\#} (\hdag T) _{\Q}) $,
$\E \in D  _{\mathrm{parf}} (\overset{^g}{} \D ^\dag  _{\X ^\#} (\hdag T) _{\Q}) $.
\end{prop}

\begin{proof}
Nous avons vu au cours de la preuve de \ref{videsurcvlog} que
$\E$ est $\D  _{\X ^\#} (\hdag T) _{\Q} $-cohérent.
D'après \cite[3.6.2]{Be1}, il existe un entier $m_0$ suffisamment grand tel que $\E$ provienne par extension d'un
$\widehat{\B} ^{(m_0)} _\X  \otimes _{\O _{\X}} \D  _{\X ^\#,\Q}$-module cohérent $\E ^{(m_0)}$.
La proposition étant locale, on peut supposer $\E ^{(m_0)}$ muni d'une bonne filtration.
D'après \ref{theofirstspencer} et \ref{remaQ},
pour $s$ assez grand, la première suite de Spencer
$Sp ^\bullet _{s,\widehat{\B} _\X  ^{(m_0)} \otimes _{\O _{\X}} \D  _{\X ^\#,\Q}} (\E ^{(m_0)})$
est exacte.
Comme l'extension $\widehat{\B} _\X  ^{(m_0)} \otimes _{\O _{\X}} \D  _{\X ^\#,\Q}
\rightarrow \D  _{\X ^\#} (\hdag T) _{\Q} $ est plate,
il en résulte
que la suite $\D  _{\X ^\#} (\hdag T) _{\Q} \otimes _{\widehat{\B} _\X ^{(m_0)} \otimes _{\O _{\X}} \D  _{\X ^\#,\Q}}
Sp ^\bullet _{s,\widehat{\B} _\X  ^{(m_0)} \otimes _{\O _{\X}} \D  _{\X ^\#,\Q}} (\E ^{(m_0)})$
est exacte.
Comme $\E$ est localement projectif de type fini sur
$\O _{\X} (\hdag T) _\Q$,
on remarque que cette suite donne une résolution finie de $\E$ par
des $\D  _{\X ^\#} (\hdag T) _{\Q} $-modules localement projectif de type fini.
Donc, $\E \in D  _{\mathrm{parf}} (\overset{^g}{} \D  _{\X ^\#} (\hdag T) _{\Q}) $.
Puisque l'extension $\D  _{\X ^\#} (\hdag T) _{\Q} \rightarrow \D ^\dag   _{\X ^\#} (\hdag T) _{\Q}$
est plate,
avec \ref{videsurcvlog}, on en déduit que
$\E \in D  _{\mathrm{parf}} (\overset{^g}{} \D ^\dag  _{\X ^\#} (\hdag T) _{\Q}) $.
\end{proof}

\begin{defi}
\label{defhol}
Soit $T$ un diviseur de $X_0$.
En s'inspirant de \cite[III.4.2]{virrion},
si $\E$ est un $\D ^\dag  _{\X ^\#} (\hdag T) _{\Q}$-module à gauche cohérent,
on dira que $\E$ est {\og $\D ^\dag  _{\X ^\#} (\hdag T) _{\Q}$-holonome \fg}
si, pour tout $l \not = 0$,
$\mathcal{H} ^l (\DD _{\X ^\#,T} (\E)) = 0$.
De même pour les $\D ^\dag  _{\X ^\#} (\hdag T) _{\Q}$-modules à droite cohérents.

Lorsque la log-structure est triviale, $T$ est vide et $\E$ est muni d'une structure de
Frobenius, nous retrouvons la notion d'holonomie de Berthelot (d'après \cite[III.4.2]{virrion}).
\end{defi}

\begin{lemm}
\label{omegareso}
Soit $T$ un diviseur de $X_0$.
On désigne par $\widetilde{\D} _{\X ^\#}$ l'un des faisceaux d'anneaux
$\D ^{(0)} _{\X ^\#}$,
$\widehat{\D} ^{(0)} _{\X ^\#}$
(resp. $\D _{\X ^\#,\Q}$, $\widehat{\D} ^{(m)} _{\X ^\#,\Q}$, $\D ^\dag _{\X ^\# ,\Q}$,
resp. $\D  _{\X ^\#} (\hdag T) _{\Q}$, $\D ^\dag  _{\X ^\#} (\hdag T) _{\Q}$).
Notons $\omega _{\X ^\#} (\hdag T) :=\omega _{\X ^\#} \otimes _{\O _{\X}} \O _{\X} (\hdag T)$.

\begin{enumerate}
  \item \label{omegareso1}
  Le morphisme canonique $\omega _{\X ^\#} \otimes _{\O _{\X}} \widetilde{\D} _{\X ^\#}  \rightarrow  \omega _{\X ^\#} $
  (resp. $\omega _{\X ^\#} \otimes _{\O _{\X}} \widetilde{\D} _{\X ^\#} \rightarrow  \omega _{\X ^\#,\Q} $,
  resp. $\omega _{\X ^\#} \otimes _{\O _{\X}} \widetilde{\D} _{\X ^\#} \rightarrow  \omega _{\X ^\#} (\hdag T) _{\Q}$)
  induit un quasi-isomorphisme
  $\Omega _{\X ^\#} ^\bullet \otimes _{\O _{\X}} \widetilde{\D} _{\X ^\#}   [d]
  \riso \omega _{\X ^\#} $
  (resp. $\Omega _{\X ^\#} ^\bullet \otimes _{\O _{\X}} \widetilde{\D} _{\X ^\#}   [d]
  \riso \omega _{\X ^\#,\Q} $,
  resp. $\Omega _{\X ^\#} ^\bullet \otimes _{\O _{\X}} \widetilde{\D} _{\X ^\#}   [d]
  \riso \omega _{\X ^\#} (\hdag T) _{\Q}$).

\item \label{DO=O}
  En considérant $\O _{\X}$ (resp. $\O _{\X,\Q}$, resp. $\O _{\X} (\hdag T) _{\Q}$) muni de sa structure canonique de
 $\widetilde{\D} _{\X ^\#}$-module, on dispose de l'isomorphisme canonique : $\DD _{\X ^\#} ( \O _{\X}) \riso \O _{\X}$
  (resp. $\DD _{\X ^\#} ( \O _{\X,\Q}) \riso \O _{\X,\Q}$,
  resp. $\DD _{\X ^\#,T} ( \O _{\X} (\hdag T) _{\Q}) \riso \O _{\X} (\hdag T) _{\Q}$).
\end{enumerate}

\end{lemm}

\begin{proof}
Traitons d'abord le cas non-respectif.
  De manière analogue à \cite[4.1.1]{Be2}, on obtient un quasi-isomorphisme :
$\Omega _{\X ^\#} ^\bullet \otimes _{\O _{\X}} \D ^{(0)} _{\X ^\#} [d]
  \riso \omega _{\X ^\#} $. En lui appliquant le foncteur exact
  $- \otimes _{\D ^{(0)} _{\X ^\#}} \widetilde{\D} _{\X ^\#}$, on en déduit (\ref{omegareso1}) via l'isomorphisme
  canonique
$\omega _{\X ^\#} \otimes _{\D ^{(0)} _{\X ^\#}} \widetilde{\D} _{\X ^\#}\riso \omega _{\X ^\#} $.

Il découle de \ref{theofirstspencer} que le complexe de Spencer
\begin{equation}
\label{sspencerappOpre}
  0 \rightarrow
\D ^{(0)} _{\X ^\#} \otimes _{\O _{\X}} \wedge ^{d} \mathcal{T} _{\X ^\#}
\overset{\epsilon}{\rightarrow}
\cdots
\overset{\epsilon}{\rightarrow}
\D ^{(0)} _{\X ^\#} \otimes _{\O _{\X}} \wedge ^{1} \mathcal{T} _{\X ^\#}
\overset{\epsilon}{\rightarrow}
\D ^{(0)} _{\X ^\#}
\rightarrow
\O _{\X} \rightarrow
0
\end{equation}
est exacte.
En lui appliquant le foncteur exact
$\widetilde{\D} _{\X ^\#}\otimes _{\D ^{(0)} _{\X ^\#}} -$,
on obtient la suite exacte :
\begin{equation}
\label{sspencerappO}
  0 \rightarrow
\widetilde{\D} _{\X ^\#}\otimes _{\O _{\X}} \wedge ^{d} \mathcal{T} _{\X ^\#}
\overset{\epsilon}{\rightarrow}
\cdots
\overset{\epsilon}{\rightarrow}
\widetilde{\D} _{\X ^\#}\otimes _{\O _{\X}} \wedge ^{1} \mathcal{T} _{\X ^\#}
\overset{\epsilon}{\rightarrow}
\widetilde{\D} _{\X ^\#}
\rightarrow
\O _{\X}
\rightarrow
0.
\end{equation}
Cela implique :
$\R \mathcal{H} om _{ \widetilde{\D}   _{\X }}
(\O _{\X}, \widetilde{\D}   _{\X }  )[d]
\riso
\Omega _{\X ^\#} ^\bullet \otimes _{\O _{\X}} \widetilde{\D}  _{\X ^\#} [d]
\riso \omega _{\X ^\#}.$
On conclut par passage de gauche à droite.

Abordons à présent les autres cas.
On déduit des isomorphismes $\omega _{\X ^\#} \otimes _{\D ^{(m)} _{\X ^\#}} \widehat{\D} ^{(m)} _{\X ^\#}
\riso \omega _{\X ^\#} $ pour $m$ variable le suivant
$\omega _{\X ^\#} \otimes _{\D ^{(0)} _{\X ^\#}} \widehat{\D} ^{(m)} _{\X ^\#,\Q}
\riso \omega _{\X ^\#,\Q} $. De même pour les autres anneaux.
On traite alors les autres cas de même que le premier.
\end{proof}

\begin{theo}
\label{compduauxisoc}
  Soient $T$ un diviseur de $X_0$, $\E $ un log isocristal sur $\X ^\#$ surconvergent le long de $T$.

Avec les notations de \ref{dualisocsurcv},
on dispose de l'isomorphisme $\D _{\X ^\#} (\hdag T) _{\Q}$-linéaire :
   $ \E ^\vee
    \riso
    \DD _{\X ^\#,T} (\E).$
Il en résulte l'isomorphisme $\D ^\dag  _{\X ^\#} (\hdag T) _{\Q}$-linéaire :
$ \E ^\vee     \riso     \DD ^\dag _{\X ^\#,T} (\E).$

Le faisceau $\E$ est $\D ^\dag _{\X ^\# } (\hdag T) _{\Q}$-holonome (voir la définition \ref{defhol}).
\end{theo}
\begin{proof}
Avec \ref{omegareso}.\ref{DO=O},
on établit le premier isomorphisme de manière analogue à \cite[2.2.1]{caro_comparaison}.
Or, $\E ^\vee $ est toujours un log isocristal sur $\X ^\#$ surconvergent le long de $T$.
Comme le foncteur dual commute à l'extension des scalaires (voir par example \cite{virrion}),
l'isomorphisme canonique $\D  ^\dag _{\X ^\#} (\hdag T) _{\Q} \otimes _{\D  _{\X ^\#} (\hdag T) _{\Q} }\E$ et celui où
$\E$ est remplacé par $\E ^\vee$ (\ref{videsurcvlog})
nous permettent de conclure.
\end{proof}

\begin{defi}
\label{notauT+}
Soient $T$ un diviseur de $X _0$ et
$\FF \in D  _{\mathrm{parf}} (\overset{^g}{} \D ^\dag  _{\X ^\#} (\hdag T) _{\Q})$.
On définit de manière analogue à \ref{defimdirdag} l'image directe (à singularités surconvergentes le long de $T$)
de $\FF$ par $u$ en posant
$u  _{T,+} (\FF)  :=
    \D ^\dag    _{\X \leftarrow \X ^\#} (\hdag T) _\Q  \otimes
    ^\L _{ \D ^\dag    _{\X ^\#} (\hdag T) _\Q} \FF$,
où
$\D ^\dag    _{\X \leftarrow \X ^\#} (\hdag T) _\Q:=
\D ^\dag  _{\X } (\hdag T) _{\Q} \otimes _{\O _{\X}} \O _{\X } (\ZZ)$
vu comme $(\D ^\dag    _{\X } (\hdag T) _\Q ,\D ^\dag    _{\X ^\#} (\hdag T) _\Q )$-bimodule.

On définit l'image directe extraordinaire (à singularités surconvergentes le long de $T$) de $\FF$ par $u$
en posant
$u  _{T,!} (\FF)  := \DD _{\X ,T}  \circ u  _{T,+} \circ \DD _{\X ^\#,T} (\FF)$.
Lorsque le diviseur $T$ est vide, on omet de l'indiquer dans les opérations cohomologiques.

\end{defi}

\begin{lemm}
  \label{u+u!Tparf}
Soient $T$ un diviseur de $X _0$ et
$\FF \in D  _{\mathrm{parf}} (\overset{^g}{} \D ^\dag  _{\X ^\#} (\hdag T) _{\Q})$.
Alors, $u  _{T,+} (\FF) \in D  _{\mathrm{parf}} (\overset{^g}{} \D ^\dag  _{\X } (\hdag T) _{\Q})$,
$u  _{T,!} (\FF) \in D  _{\mathrm{parf}} (\overset{^g}{} \D ^\dag  _{\X } (\hdag T) _{\Q})$.
\end{lemm}

\begin{proof}
Comme la perfection est stable par dualité, il suffit de traiter le cas de l'image directe.
  De manière analogue à \ref{u+dcommg}, on établit l'isomorphisme canonique :
  $$u _{T,+}  ( \FF) \riso
    ((\omega _{\X ^\#} \otimes _{\O _{\X}} \FF)
\otimes  ^\L _{ \D ^\dag    _{\X ^\#} (\hdag T) _\Q}
\D ^\dag  _{\X } (\hdag T) _{\Q} )
  \otimes _{\O _{\X}} \omega _{\X } ^{-1}.$$
On déduit de \ref{dr-ga} que le foncteur
$\omega _{\X ^\#} \otimes _{\O _{\X}} -$ (resp. $-\otimes _{\O _{\X}} \omega _{\X } ^{-1}$)
préserve la $\D ^\dag    _{\X ^\#} (\hdag T) _\Q$-perfection
(resp. $\D ^\dag    _{\X } (\hdag T) _\Q$-perfection).
D'où
$u _{T,+}  ( \FF) \in
D  _{\mathrm{parf}} (\D ^\dag    _{\X } (\hdag T) _\Q)$.
\end{proof}

Le théorème qui suit signifie que l'holonomie d'un log-isocristal surconvergent est préservée par image directe
et image directe extraordinaire par $u$.

\begin{theo}\label{u=u!=0}
  Soient $T$ un diviseur de $X_0$, $\E $ un log isocristal sur $\X ^\#$ surconvergent le long de $T$.
\begin{enumerate}
  \item
  Pour tout $l \not = 0$,
\begin{equation} \label{u+u!=01}
  \mathcal{H} ^l (u _{T,+} (\E)) =0, \hspace{1cm} \mathcal{H} ^l (u _{T,!} (\E)) =0.
\end{equation}
\item \label{u+u!=02} On dispose des isomorphismes
$u_{T,+} (\E) \riso u_{T,!} (\E (\ZZ)) \riso
\D ^\dag  _{\X } (\hdag T) _{\Q} \otimes _{\D ^\dag  _{\X ^\#} (\hdag T) _{\Q}} \E (\ZZ)$.
\item
Les faisceaux $u _{T,+} (\E)$ et $u _{T,!} ( \E)$ sont $\D ^\dag _{\X} (\hdag T) _{ \Q}$-holonomes (voir \ref{defhol}).
\end{enumerate}
\end{theo}

\begin{proof}
Par \ref{4.3.12Be1} et \ref{u+u!Tparf}, pour établir \ref{u+u!=01}, il suffit de traiter le cas où $T$ est vide.
Pour éviter les confusions, notons $\G$ le faisceau $\E $ vu comme $\D_{\X ^\#,\Q}$-module à gauche.
Par \ref{dualrellogeq1} et \ref{corologdiag3q},
pour tout $l \not = 0$,
$\mathcal{H} ^l (u _! (\G)) =0$.
Il en résulte via \ref{u+u!} que, pour tout $l \not = 0$, $\mathcal{H} ^l (u _+ (\G)) =0$.
On conclut alors la première assertion grâce à
\ref{314be0} et \ref{extf+log}.

Comme $\E \in D  _{\mathrm{parf}} (\overset{^g}{} \D ^\dag  _{\X ^\#} (\hdag T) _{\Q})$ (voir \ref{isocparf}),
en reprenant la preuve de \ref{dualrellog},
on obtient l'isomorphisme :
$\DD _{\X,T} \circ u _{T,+} (\E) \riso
\D ^\dag  _{\X } (\hdag T) _{\Q} \otimes ^\L  _{\D ^\dag  _{\X ^\#} (\hdag T) _{\Q}}
\DD _{\X ^\#,T} (\E )$. Comme celui-ci est encore valable pour
$\DD _{\X ^\#,T} (\E )$ à la place de $\E$ et que $\DD _{\X ^\#,T}\circ\DD _{\X ^\#,T} (\E )\riso \E$,
il en découle :
$u_{T,!} (\E ) \riso
\D ^\dag  _{\X } (\hdag T) _{\Q} \otimes ^\L _{\D ^\dag  _{\X ^\#} (\hdag T) _{\Q}} \E $
(et donc pour $\E (\ZZ)$ à la place de $\E$).
D'où le deuxième isomorphisme de \ref{u+u!=02} via \ref{u+u!=01}.
Le premier s'établit de manière analogue à celle de \ref{u+u!}.

Passons à la dernière assertion.
Il découle de \ref{u+u!=01}
et du théorème \ref{compduauxisoc} que,
pour tout entier $l \not = 0$,
$\mathcal{H} ^l (u _{T,!} \DD _{\X ^\#} (\E)) =0$.
Comme
$\DD _{\X} \circ u _{T,+} (\E)\riso u _{T,!} \circ \DD _{\X ^\#} (\E)$,
il en résulte
que $u _{T,+}  (\E)$ est $\D ^\dag _{\X} (\hdag T) _{ \Q}$-holonome.
L'holonomie se préservant par dualité (voir \cite{virrion} lorsque $T$ est vide mais le cas général s'en déduit grâce à
\cite[4.3.12.(ii)]{Be1}), il en dérive celle de $u _{T,!}  (\E)$.
\end{proof}

\section{Comparaison entre complexes de de Rham non logarithmique et logarithmique}

\begin{nota}
 Soient $T$ un diviseur de $X_0$, $\U ^\#$ l'ouvert de $\X ^\#$ complémentaire de $T$,
 $\E $ un log isocristal sur $\X ^\#$ surconvergent le long de $T$.
 Dans cette section, $Z$ désigne $Z _0$ et $f$ le morphisme structural $\X \rightarrow \S$.
\end{nota}

\begin{lemm}\label{md(eof)=}
  Soient $\M\in D ^- ( \D ^\dag _{\X ^\#} (\hdag T) _{\Q} \overset{^d}{})$,
  $\FF \in D ^- (\overset{^g}{} \D ^\dag _{\X ^\#} (\hdag T) _{\Q})$.
On dispose de
l'isomorphisme canonique
\begin{equation}
  \label{md(eof)=iso}
  \M \otimes ^\L _{\D ^\dag _{\X ^\#} (\hdag T) _{\Q}} (\O _{\X, \Q} (\ZZ) \otimes _{\O _{\X,\Q}}  \FF)
  \riso
(\M \otimes  _{\O _{\X,\Q}} \O _{\X, \Q} (\ZZ))\otimes _{\D ^\dag _{\X ^\#} (\hdag T) _{\Q}} ^\L \FF.
\end{equation}
\end{lemm}

\begin{proof}
De manière analogue à \ref{transdaghat},
on dispose de l'isomorphisme de transposition
$\D ^\dag _{\X ^\#} (\hdag T) _{\Q} \otimes _{\O _{\X,\Q}} \O _{\X, \Q} (\ZZ)
\riso
\O _{\X, \Q} (\ZZ) \otimes _{\O _{\X,\Q}} \D ^\dag _{\X ^\#} (\hdag T) _{\Q}$.
En résolvant $\M$ et $\FF$ platement, on construit l'isomorphisme \ref{md(eof)=iso} comme
celui de \ref{Zdroiteàgauchebis}.
\end{proof}

\begin{theo}\label{theo-compDR}
On bénéficie de l'isomorphisme canonique dans $D (f ^{-1} \O _{\S})$ :
$$ \Omega ^\bullet _{\X ^\#,\Q} \otimes  _{\O _{\X,\Q} }\E
\riso
\Omega ^\bullet _{\X,\Q} \otimes  _{\O _{\X,\Q} }  u _{T,+} (\E) .$$
\end{theo}

\begin{proof}
On dispose des isomorphismes :
\begin{gather} \notag
\Omega ^\bullet _{\X ^\#,\Q} \otimes  _{\O _{\X,\Q} }\E
\riso
( \Omega ^\bullet _{\X ^\#,\Q} \otimes  _{\O _{\X,\Q} } \D ^\dag _{\X ^\#} (\hdag T) _{\Q} )\otimes _{\D ^\dag _{\X ^\#} (\hdag T) _{\Q}}\E
\\ \notag
\underset{\ref{omegareso}.\ref{omegareso1}}{\riso}
\omega  _{\X ^\#,\Q} (\hdag T) \otimes ^\L _{\D ^\dag _{\X ^\#} (\hdag T) _{\Q}}\E [-d]
\underset{\ref{lemmomeoomegadiese-iso2}}{\riso}
(\omega  _{\X ,\Q} (\hdag T)
\otimes _{\O _{\X, \Q}} \O _{\X, \Q} (\ZZ))
\otimes ^\L _{\D ^\dag _{\X ^\#} (\hdag T) _{\Q}}\E [-d]
\\ \notag
\underset{\ref{md(eof)=iso}}{\riso}
\omega  _{\X ,\Q} (\hdag T)
\otimes ^\L _{\D ^\dag _{\X ^\#} (\hdag T) _{\Q}}
\E (\ZZ)
[-d]
\underset{\ref{u=u!=0}.\ref{u+u!=02}}{\riso}
\omega  _{\X ,\Q} (\hdag T)
\otimes ^\L _{\D ^\dag _{\X} (\hdag T)_{\Q}}
 u _{T,+} (\E )
[-d]
\\
\underset{\ref{omegareso}.\ref{omegareso1}}{\riso}
\Omega ^\bullet _{\X,\Q} \otimes  _{\O _{\X,\Q} }  u _{T,+} (\E).
\end{gather}
\end{proof}

\begin{rema}
On obtient une seconde preuve de \ref{theo-compDR} via les isomorphismes ci-dessous :
\begin{gather}
\notag
   \R \mathcal{H} om _{\D ^\dag _{\X ^\#} (\hdag T) _{\Q}} ( \O _{\X} (\hdag T) _{\Q}, \E )
\riso
\R \mathcal{H} om _{\D ^\dag _{\X ^\#} (\hdag T) _{\Q}} (\DD _{\X ^\# ,T} ( \E ), \DD _{\X ^\# ,T} (\O _{\X} (\hdag T) _{\Q}) )
\\ \notag
\underset{\ref{omegareso}.\ref{DO=O}}{\riso}
\R \mathcal{H} om _{\D ^\dag _{\X ^\#} (\hdag T) _{\Q}} (\DD _{\X ^\# ,T} ( \E ), \O _{\X} (\hdag T) _{\Q} )
\underset{\ref{u=u!=0}.\ref{u+u!=02}}{\riso}
\R \mathcal{H} om _{\D ^\dag _{\X} (\hdag T)_{\Q}} (u _{T,!} \circ \DD _{\X ^\# ,T} ( \E ), \O _{\X} (\hdag T) _{\Q} )
\\ \notag
\riso
\R \mathcal{H} om _{\D ^\dag _{\X} (\hdag T)_{\Q}} (\DD _{\X ,T}  (\O _{\X} (\hdag T) _{\Q})  ,
\DD _{\X ,T} \circ   u _{T,!} \circ \DD _{\X ^\# ,T} ( \E ))
\\ \notag
\underset{\ref{omegareso}.\ref{DO=O}}{\riso}
\R \mathcal{H} om _{\D ^\dag _{\X} (\hdag T)_{\Q}} (\O _{\X} (\hdag T) _{\Q}  ,
\DD _{\X ,T} \circ   u _{T,!} \circ \DD _{\X ^\# ,T} ( \E ))
\riso
 \R \mathcal{H} om _{\D ^\dag _{\X} (\hdag T)_{\Q}} ( \O _{\X} (\hdag T) _{\Q}, u _{T,+} (\E) ),
\end{gather}
le dernier isomorphisme résultant de l'isomorphisme de bidualité.
\end{rema}

\begin{lemm}
Le morphisme canonique
$\O _{\X} (\hdag Z \cup T ) _{\Q} \otimes _{\O _{\X} (\hdag T ) _{\Q} } \E \rightarrow
\D ^\dag _{\X } (\hdag Z \cup T) _{\Q} \otimes  _{\D ^\dag _{\X ^\#} (\hdag T ) _{\Q} } \E $
est un isomorphisme.
Ainsi, $ \E (\hdag Z):=\D ^\dag _{\X } (\hdag Z \cup T) _{\Q} \otimes ^\L _{\D ^\dag _{\X ^\#} (\hdag T ) _{\Q} } \E $
est l'isocristal sur $Y _0\cap U_0 $ surconvergent le long de $T \cup Z$ associé à $\E$ (via l'équivalence de
catégories de Berthelot énoncée dans \cite[2.2.12]{caro_courbe-nouveau}).
\end{lemm}

\begin{proof}
L'isocristal surconvergent sur $Y _0\cap U _0$ surconvergent le long de $T \cup Z$ associé à $\E$ est isomorphe à
$\O _{\X} (\hdag Z \cup T ) _{\Q} \otimes _{\O _{\X} (\hdag T ) _{\Q} } \E $.
De plus, $\D ^\dag _{\X } (\hdag Z \cup T) _{\Q} \otimes  _{\D ^\dag _{\X ^\#} (\hdag T ) _{\Q} } \E$
est un $\D ^\dag _{\X } (\hdag Z \cup T) _{\Q}$-module cohérent dont la restriction sur
$\Y\cap \U $ est isomorphe à $\E  |\Y\cap \U $, qui est $\O _{\Y\cap \U ,\Q}$-cohérent.
D'après un théorème de Berthelot (voir \cite[2.2.12]{caro_courbe-nouveau}), il en résulte que
$\D ^\dag _{\X } (\hdag Z \cup T) _{\Q} \otimes  _{\D ^\dag _{\X ^\#} (\hdag T ) _{\Q} } \E$
est $\O _{\X} (\hdag Z \cup T) _{\Q}$-cohérent. Comme
$\O _{\X} (\hdag Z \cup T ) _{\Q} \otimes _{\O _{\X} (\hdag T ) _{\Q} } \E \rightarrow
\D ^\dag _{\X } (\hdag Z \cup T) _{\Q} \otimes  _{\D ^\dag _{\X ^\#} (\hdag T ) _{\Q} } \E $
est un morphisme de
$\O _{\X} (\hdag Z \cup T) _{\Q}$-modules cohérents qui est un isomorphisme sur $\Y\cap \U$, il dérive
de \cite[4.3.12]{Be1} que celui-ci est un isomorphisme.
\end{proof}

\begin{vide}
\label{notarho}
De manière analogue à \ref{Zdroiteàgauchebis}, on dispose de l'isomorphisme canonique
$$\D ^\dag _{\X } (\hdag T ) _{\Q} \otimes  _{\D ^\dag _{\X ^\#} (\hdag T ) _{\Q} } \E(\ZZ)
\riso
\D ^\dag _{\X } (\hdag T ) _{\Q} (\ZZ) \otimes  _{\D ^\dag _{\X ^\#} (\hdag T ) _{\Q} } \E.$$
  Il en résulte ensuite par extension l'homomorphisme
$\D ^\dag _{\X } (\hdag T ) _{\Q} \otimes  _{\D ^\dag _{\X ^\#} (\hdag T ) _{\Q} } \E(\ZZ)
\rightarrow \E(\hdag Z) $. Par \ref{u=u!=0}.\ref{u+u!=02}, ce dernier est canoniquement isomorphisme
à un homomorphisme de la forme
$\rho _{\E}$ : $u_{T,+} (\E) \rightarrow\E(\hdag Z) $.
\end{vide}

Lorsque $\ZZ$ est vide, on remarque que $\rho _{\E}$ est l'identité de $\E$.
Plus généralement, le théorème \ref{theo-rhoisocv} et le lemme \ref{rhoisoext} ci-dessous donnent des exemples
de cas où l'homomorphisme $\rho _{\E}$ est un isomorphisme.
Quoiqu'il existe des contre-exemples
(voir \ref{contre-ex}), nous terminerons par une conjecture à ce sujet.
\'Enonçons d'abord les conséquences immédiates du fait que $\rho _{\E}$ soit un isomorphisme :

\begin{prop}
  Si l'homomorphisme $\rho _{\E}$ : $u _{T,+} (\E) \rightarrow \E (\hdag Z) $ est un isomorphisme,
  alors $\E (\hdag Z) $ est un $\D ^\dag _{\X} (\hdag T) _{\Q}$-module holonome
  et on dispose de l'isomorphisme canonique dans $D (f ^{-1} \O _{\S})$ :
  $$ \Omega ^\bullet _{\X ^\#,\Q} \otimes  _{\O _{\X,\Q} }\E
\riso
\Omega ^\bullet _{\X,\Q} \otimes  _{\O _{\X,\Q} }   \E (\hdag Z)  .$$
\end{prop}

Avant d'établir \ref{theo-rhoisocv},
nous aurons besoin des trois lemmes suivants.
\begin{lemm}
  \label{theo-rhoisoO}
Le morphisme canonique
$\rho _{\O _{\X} (\hdag T) _{\Q}}$ :
  $u _{T,+} (\O _{\X} (\hdag T) _{\Q}) \rightarrow \O _{\X} (\hdag T \cup Z) _\Q $
  est un isomorphisme.
\end{lemm}

\begin{proof}

Comme $\rho _{\O _{\X} (\hdag T) _{\Q}}$
est un morphisme de $\D ^\dag _{\X } (\hdag T ) _{\Q}$-modules cohérents
(pour le second terme, voir \cite{Becohdiff}),
par \cite[4.3.12]{Be1}, on se ramène à supposer $T$ vide.
  L'assertion étant locale,
  supposons qu'il existe des coordonnées locales logarithmiques
  $t _1,\dots, t _d$
  telles que $\ZZ = V ( t_1 \dots t _s)$ (rappelons que par convention, $t _{s+1},\dots, t _d$ sont inversibles).
  Avec les notations de \ref{nota11},
  on obtient la suite exacte :
\begin{equation}
  \label{se}
  ( \D _{\X ^\#,\Q} ) ^d \overset{\psi}{\rightarrow} \D _{\X ^\#,\Q} \overset{\phi}{\rightarrow} \O _{\X } (\ZZ)_\Q \rightarrow 0,
\end{equation}
où $\phi (P) = P \cdot (1/t_1\cdots t_s)$ et
$\psi ( P _1,\dots, P _d)= \sum _{i=1} ^s P _i \partial _i t _i + \sum _{i=s+1} ^d P _i \partial _i$.
Il résulte alors de la suite exacte \cite[4.3.2.1]{Be0} que
$\D ^\dag _{\X ,\Q} \otimes _{\D _{\X ^\#,\Q}} \O _{\X } (\ZZ)_\Q \riso \O _{\X } (\hdag Z)_\Q $.
On conclut grâce à \ref{314be0eq}.

\end{proof}

\begin{lemm}
  \label{lemmtheo-rhoisocv}
Soit $\B _\X$ une $\O _\X$-algèbre commutative munie d'une structure de $\D ^{(m)} _{\X^\#}$-module à gauche compatible
à sa structure de $\O _\X$-algèbre vérifiant les conditions \ref{chap-coh-21}.
On pose comme d'habitude $\B _\X (\ZZ)=\B _\X \otimes _{\O _{\X}} \O _{\X } (\ZZ)$.

Le faisceau $\B _\X  \widehat{\otimes} _{\O _{\X}} \widehat{\D} _{\X } ^{(m)}
\otimes _{\B _\X  \widehat{\otimes} _{\O _{\X}} \widehat{\D} _{\X ^\#} ^{(m)}}
\B _\X  (\ZZ)$
est alors $\B _\X  \widehat{\otimes} _{\O _{\X}} \widehat{\D} _{\X } ^{(m)}$-cohérent.
Plus précisément, s'il existe des coordonnées locales logarithmiques
  $t _1,\dots, t _d$
  telles que $\ZZ = V ( t_1 \dots t _s)$,
  alors
  $\B _\X  \widehat{\otimes} _{\O _{\X}} \widehat{\D} _{\X ^\#} ^{(m)}/
\I \riso \B _\X (\ZZ)$,
où $\I$ est
l'idéal à gauche engendré par
$\overset{^\mathrm{t}}{}\partial _{\#,i} ^{<p^j>_{(m)}}$, avec $i=1,\dots, s$, $j=1,\dots, m$
et par
$\partial _{i} ^{<p^j>_{(m)}}$, avec $i=s+1,\dots, d$, $j=1,\dots, m$.
\end{lemm}

\begin{proof}
  L'assertion étant locale,
  supposons qu'il existe des coordonnées locales logarithmiques
  $t _1,\dots, t _d$
  telles que $\ZZ = V ( t_1 \dots t _s)$.
  Le morphisme canonique $\B _\X  \widehat{\otimes} _{\O _{\X}} \widehat{\D} _{\X ^\#} ^{(m)}
  \rightarrow \B _\X (\ZZ)$ défini par
  $P \in \B _\X  \widehat{\otimes} _{\O _{\X}} \widehat{\D} _{\X ^\#} ^{(m)} \mapsto P \cdot \frac{1}{t_1 \cdots t_s}$
  induit l'isomorphisme
$\B _\X  \widehat{\otimes} _{\O _{\X}} \widehat{\D} _{\X ^\#} ^{(m)}/
\I \riso \B _\X (\ZZ)$.
En effet, un élément $P$ de
$\B _\X  \widehat{\otimes} _{\O _{\X}} \widehat{\D} _{\X ^\#} ^{(m)}$ s'écrit de manière unique
sous la forme
$\sum _{\underline{k}\geq 0} b _{\underline{k}}
\overset{^\mathrm{t}}{} \partial ^{<k_1>}_{\#,1}
\cdots
\overset{^\mathrm{t}}{} \partial ^{<k_s>}_{\#,s}
\partial ^{<k_{s+1}>}_{s+1} \cdots \partial ^{<k_{d}>}_{d}$, où
$b _{\underline{k}} \in \B _\X$ tend vers $0$ lorsque $|\underline{k}|$ tend vers l'infini.
On calcule que $P \cdot \frac{1}{t_1 \cdots t_s}=0$ si et seulement si
$b _0 =0$. Enfin, de manière analogue à \ref{lemm-wildetildeDcoh}.\ref{lemm-wildetildeDcohii},
on vérifie que cette idéal est engendré par les éléments décrits ci-dessous.

\end{proof}

\begin{lemm}
    \label{lemmtheo-rhoisocvbis}
    Avec les notations et hypothèses de \ref{lemmtheo-rhoisocv},
    soient $\E$ un $\B _\X  \widehat{\otimes} _{\O _{\X}} \widehat{\D} _{\X } ^{(m)}$-module qui soit
    cohérent sur $\B _\X $ et
    $\FF$ un $\B _\X  \widehat{\otimes} _{\O _{\X}} \widehat{\D} _{\X } ^{(m)}$-module cohérent.
    Le faisceau $\E \otimes _{\B _\X} \FF$ est alors muni d'une structure de
    $\B _\X  \widehat{\otimes} _{\O _{\X}} \widehat{\D} _{\X } ^{(m)}$-module cohérent.
\end{lemm}

\begin{proof}
On dispose des isomorphismes $\B _\X  \otimes _{\O _{\X}} \D _{\X } ^{(m)}$-linéaires :
\begin{gather}
  \notag
  \E \otimes _{\B _\X} \FF
  \underset{\ref{HomD1}}{\riso}
  (\E \otimes _{\B _\X} (\B _\X  \otimes _{\O _{\X}} \D _{\X } ^{(m)}))
  \otimes _{\B _\X  \otimes _{\O _{\X}} \D _{\X } ^{(m)}}\FF
  \riso
  \\
  \label{isocohhat}
  (\E \otimes _{\B _\X} (\B _\X  \otimes _{\O _{\X}} \D _{\X } ^{(m)}))
  \otimes _{\B _\X  \otimes _{\O _{\X}} \D _{\X } ^{(m)}}
\B _\X  \widehat{\otimes} _{\O _{\X}} \widehat{\D} _{\X } ^{(m)}
\otimes _{\B _\X  \widehat{\otimes} _{\O _{\X}} \widehat{\D} _{\X } ^{(m)}}
  \FF
  \underset{\ref{HomD1}}{\riso}
  (\E \otimes _{\B _\X} (\B _\X  \widehat{\otimes} _{\O _{\X}} \widehat{\D} _{\X } ^{(m)}))
  \otimes _{\B _\X  \widehat{\otimes} _{\O _{\X}} \widehat{\D} _{\X } ^{(m)}}\FF.
\end{gather}
Comme $\E \otimes _{\B _\X} (\B _\X  \widehat{\otimes} _{\O _{\X}} \widehat{\D} _{\X } ^{(m)})$
(resp. $(\B _\X  \widehat{\otimes} _{\O _{\X}} \widehat{\D} _{\X } ^{(m)})\otimes _{\B _\X}\E $)
est un $\B _\X  \widehat{\otimes} _{\O _{\X}} \widehat{\D} _{\X } ^{(m)}$-module à droite (resp. à gauche)
cohérent, il est $p$-adiquement séparé et complet. On en déduit par complétion
l'isomorphisme de $\B _\X  \widehat{\otimes} _{\O _{\X}} \widehat{\D} _{\X } ^{(m)}$-bimodules de transposition :
$(\B _\X  \widehat{\otimes} _{\O _{\X}} \widehat{\D} _{\X } ^{(m)})\otimes _{\B _\X}\E
\riso
\E \otimes _{\B _\X} (\B _\X  \widehat{\otimes} _{\O _{\X}} \widehat{\D} _{\X } ^{(m)})$.
Cela implique que
$\E \otimes _{\B _\X} (\B _\X  \widehat{\otimes} _{\O _{\X}} \widehat{\D} _{\X } ^{(m)})$ est
un $\B _\X  \widehat{\otimes} _{\O _{\X}} \widehat{\D} _{\X } ^{(m)}$-module à gauche cohérent.
Via les théorèmes de type $A$ pour les
$\B _\X  \widehat{\otimes} _{\O _{\X}} \widehat{\D} _{\X } ^{(m)}$-modules cohérents,
on en déduit que
$(\E \otimes _{\B _\X} (\B _\X  \widehat{\otimes} _{\O _{\X}} \widehat{\D} _{\X } ^{(m)}))
  \otimes _{\B _\X  \widehat{\otimes} _{\O _{\X}} \widehat{\D} _{\X } ^{(m)}}\FF$
  est un $\B _\X  \widehat{\otimes} _{\O _{\X}} \widehat{\D} _{\X } ^{(m)}$-module cohérent.
On conclut alors avec \ref{isocohhat}.
\end{proof}

\begin{theo}
\label{theo-rhoisocv}
On suppose que $\E$ est en fait un isocristal sur $\X $ surconvergent le long de $T$, i.e.,
un $\D ^\dag _{\X} (\hdag T) _{ \Q}$-module cohérent $\O _{\X} (\hdag T \cup Z) _\Q$-cohérent.

Le morphisme canonique $\rho _{\E}$ :
  $u _{T,+} (\E) \rightarrow \E (\hdag Z )  $
  est alors un isomorphisme.
  Par \ref{u=u!=0}, l'isocristal $\E (\hdag Z )  $ sur $\X $ surconvergent le long de $T \cup Z$
  est donc un $\D ^\dag _{\X} (\hdag T) _{ \Q}$-module holonome.
\end{theo}
\begin{proof}

D'après \cite[4.4.5]{Be1} (et avec la remarque \cite[4.4.6]{Be1}),
on peut construire une suite
croissante d'entiers $(n _m) _{m\in \N}$ avec $n _m \geq m$
telle qu'il existe un
$\widehat{\B} _\X ^{(n _0)} (T ) _{\Q}$-module cohérent $\E ^{(0)}$ tel
que $\E ^{(m)} :=\widehat{\B} _\X ^{(n _m)} (T ) _{\Q} \otimes _{\widehat{\B} _\X ^{(n _0)} (T ) _{\Q}} \E ^{(0)}$
soit muni d'une structure canonique de
$\widehat{\B} _\X ^{(n_m)} (T ) \widehat{\otimes} _{\O _{\X}} \widehat{\D} _{\X ,\Q} ^{(m)}$-module topologiquement
nilpotent induisant l'isomorphisme
$\D ^\dag _{\X} (\hdag T) _{\Q}$-linéaire :
$\E \riso \smash{\underset{\longrightarrow}{\lim}} _m \E ^{(m)}$.

Par \cite[4.4.7]{Be1}, il existe un
$\widehat{\B} _\X ^{(n_m)} (T ) \widehat{\otimes} _{\O _{\X}} \widehat{\D} _{\X} ^{(m)}$-module cohérent $\smash{\overset{\circ}{\E}}  ^{(m)}$,
$\widehat{\B} _\X ^{(m)} (T )$-cohérent
tel que $\smash{\overset{\circ}{\E}} _\Q ^{(m)} \riso \E ^{(m)}$.
Posons :
$\widehat{\B} _\X ^{(n_m)} (T ,\ZZ):= \widehat{\B} _\X ^{(n_m)} (T ) \otimes _{\O _{\X}} \O _{\X } (\ZZ)$,
$\widehat{\D} _{\X} ^{(m)}(T):=
\widehat{\B} _\X ^{(n_m)} (T ) \widehat{\otimes} _{\O _{\X}} \widehat{\D} _{\X} ^{(m)}$,
$\widehat{\D} _{\X  ^\#} ^{(m)}(T):=
\widehat{\B} _\X ^{(n_m)} (T ) \widehat{\otimes} _{\O _{\X}} \widehat{\D} _{\X ^\#} ^{(m)}$.
Comme
$\smash{\overset{\circ}{\E}}  ^{(m)}$ est $\widehat{\B} _\X ^{(n_m)} (T )$-cohérent,
par \ref{lemmtheo-rhoisocv} et \ref{lemmtheo-rhoisocvbis},
le faisceau
$(\widehat{\D} _{\X} ^{(m)}(T) \otimes _{\widehat{\D} _{\X ^\#} ^{(m)}(T)}
\widehat{\B} _\X ^{(n_m)} (T ,\ZZ) )\otimes _{\widehat{\B} _\X ^{(n_m)} (T )} \smash{\overset{\circ}{\E}}  ^{(m)}$
est muni d'une structure canonique de
$\widehat{\D} _{\X} ^{(m)}(T) $-module cohérent.
On obtient alors par extension le morphisme canonique :
\begin{equation}
\label{ACCNbishat}
\widehat{\D} _{\X} ^{(m)}(T) \otimes _{\widehat{\D} _{\X ^\#} ^{(m)}(T)}
(\widehat{\B} _\X ^{(n_m)} (T ,\ZZ) \otimes _{\widehat{\B} _\X ^{(n_m)} (T )} \smash{\overset{\circ}{\E}}  ^{(m)})
\rightarrow
(\widehat{\D} _{\X} ^{(m)}(T) \otimes _{\widehat{\D} _{\X ^\#} ^{(m)}(T)}
\widehat{\B} _\X ^{(n_m)} (T ,\ZZ) )\otimes _{\widehat{\B} _\X ^{(n_m)} (T )} \smash{\overset{\circ}{\E}}  ^{(m)},
\end{equation}
qui vérifie $1 \otimes (x\otimes y) \mapsto (1 \otimes x )\otimes y $, où
$x \in \widehat{\B} _\X ^{(n_m)} (T ,\ZZ) $,
$y \in  \smash{\overset{\circ}{\E}}  ^{(m)}$.
De même que pour \ref{313Be0} (ou \cite[3.1.3]{Be1}), on vérifie que
$\widehat{\B} _\X ^{(n_m)} (T ,\ZZ) \otimes _{\widehat{\B} _\X ^{(n_m)} (T )} \smash{\overset{\circ}{\E}}  ^{(m)}$
est
$\widehat{\D} _{\X ^\#} ^{(m)}(T)$-cohérent.
Le terme de gauche de \ref{ACCNbishat} est donc, comme celui de droite,
$p$-adiquement séparé et complet.
Or, par \ref{A2CNbis},
\ref{ACCNbishat} est un isomorphisme modulo $\mathfrak{m} ^{i+1}$ pour tout entier $i\geq0$.
Cela implique que \ref{ACCNbishat} est un isomorphisme.
D'où en tensorisant par $\Q$ :
\begin{equation}
\label{ACCNbishatQ}
\widehat{\D} _{\X} ^{(m)}(T) _{\Q} \otimes _{\widehat{\D} _{\X ^\#} ^{(m)}(T)_{\Q}}
\E  ^{(m)} (\ZZ)
\riso
(\widehat{\D} _{\X} ^{(m)}(T) _{\Q} \otimes _{\widehat{\D} _{\X ^\#} ^{(m)}(T) _{\Q}}
\widehat{\B} _\X ^{(n_m)} (T ,\ZZ) _{\Q} )\otimes _{\widehat{\B} _\X ^{(n_m)} (T ) _{\Q}}
\E  ^{(m)}.
\end{equation}
Pour terminer la preuve, nous aurons besoin du lemme suivant.

\begin{lemm}
\label{theo-rhoisoO-lemm1}
L'homomorphisme canonique
\begin{equation}
\label{theo-rhoisoO-lemm1iso}
\E  ^{(m)} (\ZZ)  \rightarrow
  \widehat{\D} _{\X ^\#} ^{(m)}(T)_{\Q} \otimes _{\widehat{\D} _{\X ^\#} ^{(0)}(T)_{\Q}}
 \E  ^{(0)} (\ZZ)
\end{equation}
est un isomorphisme.
\end{lemm}

\begin{proof}
On procède comme pour \cite[4.4.10]{Be1} : de manière analogue à \cite[4.4.9]{Be1}, l'homomorphisme canonique :
\begin{equation}
  \notag \label{theo-rhoisoO-lemm1eq1}
  \smash{\overset{\circ}{\E}}  ^{(m)} \rightarrow
\widehat{\B} _\X ^{(n_m)} (T ) \widehat{\otimes} _{\O _{\X}} \smash{\widehat{\D}} _{\X ^\#} ^{(m)}
\otimes _{\widehat{\B} _\X ^{(n_m)} (T ) \otimes _{\O _{\X}} \D _{\X^\#} ^{(m)}}
\smash{\overset{\circ}{\E}}  ^{(m)}
\end{equation}
est un isomorphisme et de même en remplaçant {\og $\D _{\X ^\#} ^{(m)}$ \fg} par
{\og $\D _{\X ^\#} ^{(0)}$ \fg}.
Il en résulte l'isomorphisme :
\begin{equation}
  \label{theo-rhoisoO-lemm1eq2}
  \E  ^{(m)}
  \riso
\widehat{\B} _\X ^{(n_m)} (T ) \widehat{\otimes} _{\O _{\X}} \smash{\widehat{\D}} _{\X ^\#,\Q} ^{(m)}
\otimes _{\widehat{\B} _\X ^{(n_m)} (T ) \widehat{\otimes} _{\O _{\X}} \smash{\widehat{\D}} _{\X ^\#,\Q} ^{(0)}}
\E  ^{(m)}
\end{equation}
On établit comme dans la preuve de \cite[4.4.8]{Be1} que l'homomorphisme canonique
\begin{equation}
  \label{theo-rhoisoO-lemm1eq3}
  \smash{\overset{\circ}{\E}}  ^{(m)} \rightarrow
\widehat{\B} _\X ^{(n_m)} (T ) \widehat{\otimes} _{\O _{\X}} \smash{\widehat{\D}} _{\X ^\#} ^{(0)}
\otimes _{\widehat{\B} _\X ^{(n_0)} (T ) \widehat{\otimes} _{\O _{\X}} \smash{\widehat{\D}} _{\X^\#} ^{(0)}}
\smash{\overset{\circ}{\E}}  ^{(0)}
\end{equation}
est un isomorphisme.
En tensorisant \ref{theo-rhoisoO-lemm1eq3} par $\Q$, on conclut avec \ref{theo-rhoisoO-lemm1eq2}.
\end{proof}
Revenons à présent à la preuve du théorème.
On déduit de \ref{theo-rhoisoO-lemm1} les isomorphismes suivants :
\begin{gather}
  \notag
  \smash{\underset{\underset{m}{\longrightarrow}}{\lim}} \,
\widehat{\smash{\widehat{\D}}} _{\X} ^{(m)}(T) _{\Q} \otimes _{\widehat{\smash{\widehat{\D}}} _{\X ^\#} ^{(m)}(T)_{\Q}}
\E  ^{(m)} (\ZZ)
\underset{\ref{theo-rhoisoO-lemm1iso}}{\riso}
\smash{\underset{\underset{m}{\longrightarrow}}{\lim}} \,
  \widehat{\D} _{\X } ^{(m)}(T)_{\Q} \otimes _{\widehat{\D} _{\X ^\#} ^{(0)}(T)_{\Q}}
 \E  ^{(0)} (\ZZ)
\\ \label{theo-rhoisoO-lemm1eq4}
\riso
\D ^\dag _{\X } (\hdag T ) _{\Q}
\otimes _{\widehat{\D} _{\X ^\#} ^{(0)}(T)_{\Q}}
 \E  ^{(0)} (\ZZ)
\underset{\ref{theo-rhoisoO-lemm1iso}}{\riso}
\D ^\dag _{\X } (\hdag T ) _{\Q}
\otimes  _{\D ^\dag _{\X ^\#} (\hdag T ) _{\Q} } \E(\ZZ).
\end{gather}
De même, comme
$\E ^{(m)} =\widehat{\B} _\X ^{(n _m)} (T ) _{\Q} \otimes _{\widehat{\B} _\X ^{(n _0)} (T ) _{\Q}} \E ^{(0)}$
et en utilisant \ref{theo-rhoisoO-lemm1eq4} dans le cas où $\E = \O _{\X } (\hdag T ) _{\Q}$,
on obtient l'isomorphisme :
\begin{equation}
\label{theo-rhoisoO-lemm1eq5}
\smash{\underset{\underset{m}{\longrightarrow}}{\lim}} \,
  (\widehat{\D} _{\X} ^{(m)}(T) _{\Q} \otimes _{\widehat{\D} _{\X ^\#} ^{(m)}(T) _{\Q}}
\widehat{\B} _\X ^{(n_m)} (T ,\ZZ) _{\Q} )\otimes _{\widehat{\B} _\X ^{(n_m)} (T ) _{\Q}}
\E  ^{(m)}
\riso
(\D ^\dag _{\X } (\hdag T ) _{\Q}
\otimes  _{\D ^\dag _{\X ^\#} (\hdag T ) _{\Q} } \O _{\X } (\hdag T ) _{\Q} (\ZZ))
\otimes _{\O _{\X } (\hdag T ) _{\Q}}
\E.
\end{equation}
Il résulte de \ref{ACCNbishatQ}, \ref{theo-rhoisoO-lemm1eq4}, \ref{theo-rhoisoO-lemm1eq5} l'isomorphisme :
\begin{equation}
  \notag
  \D ^\dag _{\X } (\hdag T ) _{\Q}
\otimes  _{\D ^\dag _{\X ^\#} (\hdag T ) _{\Q} } \E(\ZZ)
\riso
  (\D ^\dag _{\X } (\hdag T ) _{\Q}
\otimes  _{\D ^\dag _{\X ^\#} (\hdag T ) _{\Q} } \O _{\X } (\hdag T ) _{\Q} (\ZZ))
\otimes _{\O _{\X } (\hdag T ) _{\Q}}
\E.
\end{equation}
Or, d'après \ref{theo-rhoisoO},
$\D ^\dag _{\X } (\hdag T ) _{\Q}
\otimes  _{\D ^\dag _{\X ^\#} (\hdag T ) _{\Q} } \O _{\X } (\hdag T ) _{\Q} (\ZZ) \riso
u_{T,+} (\O _{\X } (\hdag T ) _{\Q}) \riso \O _{\X } (\hdag T \cup Z) _{\Q}$. D'où le résultat.
\end{proof}

\begin{lemm}
\label{rhoisoext}
Soit $0\rightarrow \E ' \rightarrow \E \rightarrow \E''\rightarrow 0$ une suite exacte de
log-isocristaux sur $\X ^\#$ surconvergents le long de $T$.
Si $\rho _{\E'}$ et $\rho _{\E''}$ sont des isomorphismes alors $\rho _{\E}$ l'est aussi.
\end{lemm}

\begin{proof}
Cela résulte du lemme des cinq.
\end{proof}

\begin{rema}
\label{contre-ex}
  Si on ne fait aucune hypothèse sur les exposants, l'homomorphisme $\rho$ n'est pas toujours un isomorphisme.
  Voici deux contre-exemples :

  $\bullet$ Lorsque $\ZZ$ est non vide,
  on vérifie que $\rho _{\O _{\X} (\hdag T) _{\Q} (-\ZZ)}$ n'est pas un isomorphisme.
  En effet, par \cite[4.3.12]{Be1}, on se ramène au cas où $T$ est vide.
  Supposons qu'il existe des coordonnées locales logarithmiques
  $t _1,\dots, t _d$
  telles que $\ZZ = V ( t_1 \dots t _s)$ (rappelons que par convention, $t _{s+1},\dots, t _d$ sont inversibles).
On calcule $\D _{\X ^\#,\Q} / \D _{\X ^\#,\Q} (\partial _{\#,1}, \dots, \partial _{\#,d})
\riso \O _{\X,\Q} $. Par \ref{314be0eq}, on en déduit :
$\D ^\dag _{\X ,\Q} / \D ^\dag _{\X ,\Q} (t _1\partial _{1}, \dots, t_d\partial _{d})
\riso \D ^\dag _{\X ,\Q}  \otimes _{\D ^\dag _{\X ^\# ,\Q} } \O _{\X,\Q} \riso u _{T,+} (\O _{\X,\Q} (-\ZZ))$.
Or,
$$\D ^\dag _{\X ,\Q} / \D ^\dag _{\X ,\Q} (\partial _{1} t _1, \dots, \partial _{s}t _s, \partial _{s+1},\dots, \partial _{d})
\underset{\cite[4.3.2.1]{Be0}}{\riso}
\O _{\X} (\hdag Z) _{\Q}
\riso \O _{\X} (\hdag Z) _{\Q} \otimes _{\O _{\X,\Q}}\O _{\X,\Q} (-\ZZ) \riso (\O _{\X,\Q} (-\ZZ) ) (\hdag Z).$$
Lorsque $s \geq 1$,
on conclut alors en remarquant $\D ^\dag _{\X ,\Q} (t _1\partial _{1}, \dots, t_d\partial _{d}) \not =
\D ^\dag _{\X ,\Q} (\partial _{1} t _1, \dots, \partial _{s}t _s, \partial _{s+1},\dots, \partial _{d})$.

  $\bullet$ Lorsque $\X$ est propre et $T$ est vide, le fait que
  $\rho _{\E}$ soit un isomorphisme implique que la cohomologie rigide de l'isocristal surconvergent $\E (\hdag Z)$ serait
  toujours de dimension finie. Il suffit alors de regarder l'isocristal surconvergent (qui provient d'un log isocristal convergent)
  décrit par Berthelot dans la dernière remarque de \cite{Be1}
  pour constater que ce n'est pas toujours le cas.
\end{rema}

\begin{conj}
\label{conj}
On suppose que
\begin{itemize}
  \item toutes les différences des exposants le long des composantes
irréductibles de $Z$ ne sont pas des entiers $p$-adiques de Liouville,
  \item tous les exposants le long des composantes
irréductibles de $Z$ ne sont pas des entiers $p$-adiques de Liouville
et ne sont pas strictement positifs.
\end{itemize}

Le morphisme $\rho _{\E}$ : $u _{T,+} (\E) \rightarrow \E (\hdag Z) $ est alors un isomorphisme.
\end{conj}

\bibliographystyle{smfalpha}
\bibliography{bib1}

\noindent Daniel Caro\\
Arithmétique et Géométrie algébrique, Bât. 425\\
Université Paris-Sud\\
91405 Orsay Cedex\\
France.\\
email: daniel.caro@math.u-psud.fr

\end{document}